\documentclass[11pt,reqno,tbtags]{amsart} 
\usepackage[utf8]{inputenc}

\usepackage[a4paper,width=150mm,top=25mm,bottom=25mm]{geometry}
\usepackage{mathtools}
\usepackage{suffix}
\usepackage{enumerate}
\usepackage{enumitem}
\usepackage{listings}

\makeatletter
\newcommand*{\rom}[1]{\expandafter\@slowromancap\romannumeral #1@}
\makeatother
\usepackage{pgf, tikz}
\usepackage{subcaption}
\usetikzlibrary{arrows, automata}
\usepackage{float}

\usepackage{parskip}
\usepackage{amsmath,amsthm,amssymb}
\numberwithin{equation}{section}

\allowdisplaybreaks
\usepackage{bbm}
\usepackage[makeroom]{cancel}

\usepackage{xcolor}
\definecolor{coolblack}{rgb}{0.0, 0.18, 0.39}

\usepackage[breaklinks=true]{hyperref}
\hypersetup{
colorlinks=true,
linkcolor=blue,
urlcolor=blue,
citecolor=black
}

\title{The number of descendants in a preferential attachment graph}
\author{Svante Janson, Tiffany Y.\ Y.\ Lo}
\thanks{Supported by the Knut and Alice Wallenberg Foundation, 
 Ragnar Söderberg Foundation, 
the Swedish Research Council  (VR), 
and Sverker Lerheden Foundation. 
}
\address{Department of Mathematics, Uppsala University, PO Box 480,
SE-751~06 Uppsala, Sweden}
\email{svante.janson@math.uu.se }
\address{Department of Mathematics, Stockholm University, SE-106 91 Stockholm, Sweden}
\email{tiffany.y.y.lo@math.su.se}
\date{18 December, 2024}

\theoremstyle{plain}
\newtheorem{theorem}{Theorem}[section]

\newtheorem{lemma}[theorem]{Lemma}

\theoremstyle{definition}

\newtheorem{definition}[theorem]{Definition}

\newtheorem{remark}[theorem]{Remark}

\newtheorem{problem}[theorem]{Problem}

\newcommand{\al}{\alpha}
\newcommand{\IP}{\mathbbm{P}}
\newcommand\E{\operatorname{\mathbb E}} 

\newcommand{\G}{\Gamma}

\newcommand\numberthis{\addtocounter{equation}{1}\tag{\theequation}}
\newcommand{\wt}{\widetilde}

\newcommand{\bB}{\mathbf{B}}

\newcommand{\wh}{\widehat}
\newcommand{\tone}{\mathbf{1}}

\newcommand{\cF}{\mathcal{F}}

\newcommand{\var}{\mathrm{Var}}

\newcommand{\toinf}{\to\infty}
\renewcommand{\leq}{\leqslant}
\renewcommand{\geq}{\geqslant}
\renewcommand{\phi}{\varphi}
\newcommand{\eps}{\varepsilon}

\newenvironment{romenumerate}[1][-10pt]{
\addtolength{\leftmargini}{#1}\begin{enumerate}
 }{\end{enumerate}}

\renewcommand{\le}{\leq}
\renewcommand{\ge}{\geq}
\newcommand{\refT}[1]{Theorem~\ref{#1}}
\newcommand{\refTs}[1]{Theorems~\ref{#1}}

\newcommand{\refL}[1]{Lemma~\ref{#1}}
\newcommand{\refLs}[1]{Lemmas~\ref{#1}}
\newcommand{\refR}[1]{Remark~\ref{#1}}

\newcommand{\refS}[1]{Section~\ref{#1}}
\newcommand{\refSs}[1]{Sections~\ref{#1}}
\newcommand{\refApp}[1]{Appendix~\ref{#1}}

\newcommand{\refD}[1]{Definition~\ref{#1}}

\newcommand{\refStep}[1]{Step~\ref{#1}}
\newcommand{\refSteps}[1]{Steps~\ref{#1}}
\newcommand\ga{\alpha}
\newcommand\gb{\beta}
\newcommand\gd{\delta}
\newcommand\gD{\Delta}

\newcommand\GF{\Phi}

\newcommand\gG{\Gamma}
\newcommand\gk{\varkappa}
\newcommand\kk{\kappa}
\newcommand\gl{\lambda}

\newcommand\go{\omega}

\newcommand\gth{\theta}

\newcommand\cA{\mathcal A}
\newcommand\cB{\mathcal B}

\newcommand\cU{\mathcal U}
\newcommand\cX{\mathcal X}
\newcommand\tcB{\widetilde{\mathcal B}}
\newcommand\tU{\widetilde{U}}
\newcommand\sU{\mathsf{U}}
\newcommand\hA{\widehat{A}}
\newcommand\xD{\widehat{D}}

\newcommand\hH{\widehat{H}}
\newcommand\hP{\widehat{P}}
\newcommand\hT{\widehat{T}}
\newcommand\hV{\widehat{V}}

\newcommand\set[1]{\ensuremath{\{#1\}}}

\newcommand\xpar[1]{(#1)}
\newcommand\bigpar[1]{\bigl(#1\bigr)}
\newcommand\Bigpar[1]{\Bigl(#1\Bigr)}
\newcommand\biggpar[1]{\biggl(#1\biggr)}
\newcommand\lrpar[1]{\left(#1\right)}
\newcommand\bigsqpar[1]{\bigl[#1\bigr]}
\newcommand\sqpar[1]{[#1]}

\newcommand\abs[1]{\lvert#1\rvert}
\newcommand\bigabs[1]{\bigl\lvert#1\bigr\rvert}

\newcommand\biggabs[1]{\biggl\lvert#1\biggr\rvert}

\newcommand\downto{\searrow}
\newcommand\upto{\nearrow}
\newcommand{\tend}{\longrightarrow}
\newcommand\dto{\overset{\mathrm{d}}{\tend}}
\newcommand\pto{\overset{\mathrm{p}}{\tend}}
\newcommand\asto{\overset{\mathrm{a.s.}}{\tend}}
\newcommand\ktoo{\ensuremath{{k\to\infty}}}
\newcommand\ntoo{\ensuremath{{n\to\infty}}}

\newcommand\ttoo{\ensuremath{{t\to\infty}}}

\newcommand\Bin{\operatorname{Bin}}

\renewcommand\P{\IP}
\newcommand\Var{\operatorname{Var}}

\newcommand\bbR{\mathbb R}

\newcommand\qw{^{-1}}

\newcommand\qq{^{1/2}}

\newcommand\intoo{\int_0^\infty}
\newcommand\dd{\,\mathrm{d}}
\newcommand\ddx{\mathrm{d}}

\newcommand\eqd{\overset{\mathrm{d}}{=}}

\newcommand\lhs{left-hand side}
\newcommand\rhs{right-hand side}
\newcommand\hcY{\widehat{\mathcal Y}}
\newcommand\cY{\mathcal{Y}}
\newcommand\nn{^{(n)}}
\newcommand\xfrac[2]{#1/#2}
\newcommand\whp{w.h.p.}
\newcounter{steps}
\newcommand\stepp{\par\noindent\refstepcounter{steps}%
  \emph{Step \arabic{steps}. }\noindent}
\newcommand\steppx[1]{\par\noindent\refstepcounter{steps}%
  \emph{Step \arabic{steps}. #1}\noindent}
\newcommand\resetsteps{\setcounter{steps}{0}}
\newcommand\oi{\ensuremath{[0,1]}}
\newcommand\nxoo{_{n=1}^\infty}
\newcommand\Beta{\mathrm{Beta}}
\newcommand\GAMMA{\mathrm{Gamma}}
\newcommand\Phix{\widehat\Psi}
\newcommand\xM{\mathfrak M}
\newcommand\tM{\widetilde M}
\newcommand\gln{\gl_n}
\newcommand\tgb{\tilde\beta}
\newcommand\op{o_{\mathrm p}}

\newcommand\Mx{M_*}
\newcommand\Mxx{\Mx}
\newcommand\bignorm[1]{\bigl\lVert#1\bigr\rVert}

\newcommand\lrnorm[1]{\left\lVert#1\right\rVert}

\newcommand\MM{\widehat M}

\begingroup
  \divide\count255 by 60
  \multiply\count255 by -60
  \advance\count255 by \time
  \ifnum \count255 < 10 \xdef\klockan{\the\count1.0\the\count255}
  \else\xdef\klockan{\the\count1.\the\count255}\fi
\endgroup

\DeclarePairedDelimiter\ceil{\lceil}{\rceil}
\DeclarePairedDelimiter\floor{\lfloor}{\rfloor}

\def\given{\typeout{Command 'given' should only be used within bracket command}}
\newcounter{@bracketlevel}
\def\@bracketfactory#1#2#3#4#5#6{
\expandafter\def\csname#1\endcsname##1{%
\addtocounter{@bracketlevel}{1}%
\global\expandafter\let\csname @middummy\alph{@bracketlevel}\endcsname\given%
\global\def\given{\mskip#5\csname#4\endcsname\vert\mskip#6}\csname#4l\endcsname#2##1\csname#4r\endcsname#3%
\global\expandafter\let\expandafter\given\csname @middummy\alph{@bracketlevel}\endcsname
\addtocounter{@bracketlevel}{-1}}%
}
\def\bracketfactory#1#2#3{%
\@bracketfactory{#1}{#2}{#3}{relax}{1mu plus 0.25mu minus 0.25mu}{0.6mu plus 0.15mu minus 0.15mu}
\@bracketfactory{b#1}{#2}{#3}{big}{1mu plus 0.25mu minus 0.25mu}{0.6mu plus 0.15mu minus 0.15mu}
\@bracketfactory{bb#1}{#2}{#3}{Big}{2.4mu plus 0.8mu minus 0.8mu}{1.8mu plus 0.6mu minus 0.6mu}
\@bracketfactory{bbb#1}{#2}{#3}{bigg}{3.2mu plus 1mu minus 1mu}{2.4mu plus 0.75mu minus 0.75mu}
\@bracketfactory{bbbb#1}{#2}{#3}{Bigg}{4mu plus 1mu minus 1mu}{3mu plus 0.75mu minus 0.75mu}
}
\bracketfactory{clc}{\lbrace}{\rbrace}
\bracketfactory{clr}{(}{)}
\bracketfactory{cls}{[}{]}
\bracketfactory{abs}{\lvert}{\rvert}
\bracketfactory{norm}{\Vert}{\Vert}
\bracketfactory{floor}{\lfloor}{\rfloor}
\bracketfactory{ceil}{\lceil}{\rceil}
\bracketfactory{angle}{\langle}{\rangle}

\begin{document}
\begin{abstract}
We study the number $X^{(n)}$ of vertices that can be reached from the last
added vertex $n$ via a directed path (the descendants) in the standard
preferential attachment graph. In this model, vertices are sequentially
added, each born with outdegree $m\ge 2$; the endpoint of each outgoing edge
is chosen among previously added vertices with probability proportional to
the current degree of the vertex
plus some number $\rho$.

We show that $X^{(n)}/n^\nu$ converges in distribution as $n\to\infty$,
where $\nu$ depends on both $m$ and $\rho$, and the limiting distribution is
given by a product of a constant factor and the $(1-\nu)$-th power of a
$\GAMMA(m/(m-1),1)$ variable. The proof uses a P\'olya urn representation of
preferential attachment graphs, and the arguments of Janson (2024) where the
same problem was studied in uniform attachment graphs.   
Further results, including convergence of all moments and analogues for the version with possible self-loops are  provided. 
\end{abstract}

\maketitle
\section{Introduction}
Preferential attachment models have emerged as a popular class of random
 graphs since it was proposed in \cite{BA1997} as an explanation for the power-law
 degree sequences observed in real-world networks.
There are several versions of these models, differing in minor details,
see e.g.\ \cite{vdh2017};
we will use the version defined below, which is  
the sequential model in \cite{Berger2014}. In this version,
 self-loops are not allowed but multiple edges are possible. 
The graph is often treated as undirected, but  we regard it as directed,
with all edges directed from the 
younger vertex (with larger label) to the older vertex (with smaller label).

\begin{definition}[Preferential attachment graph]\label{de:pa}
  Fix an integer $m\geq 2$ and a real number $\rho>-m$, 
and let $(G_n)_{n\geq 1}$ be the sequence of random graphs  
that are generated as follows;
$G_n$ has $n$  vertices with labels in $[n]:=\{1,\dots,n\}$. 
The initial graph $G_1$ consists of a  single vertex (labelled 1) with no
edges.
Given $G_{n-1}$, we
construct $G_{n}$ from $G_{n-1}$ by adding the new vertex with label $n$, 
and sequentially attaching $m$ edges between vertex~$n$ and
at most $m$ vertices in $G_{n-1}$ as follows. 
Let {$d_j(n)$} be the degree of vertex $j$ in $G_n$.
If $n\ge2$, each outgoing edge of vertex $n$ is attached to vertex
$j\in[n-1]$ with probability proportional to 
$\rho$ + the current degree of vertex~$j$. 
(In particular,
if $n=2$, we add $m$ edges from vertex~2 to vertex 1.)
This means that
the first {outgoing} edge of
  vertex $n$ is attached to vertex $j\in[n-1]$ with probability
  \begin{align}\label{eq:pa1}
      \frac{d_j(n-1)+\rho}{2m(n-2)+(n-1)\rho};
  \end{align}
  noting that $\sum^{n-1}_{k=1}d_k(n-1)=2m(n-2)$ and 
$d_j(n-1)+\rho\ge m+\rho >0$.
Furthermore,  given that the first $1\leq k\leq m-1$ outgoing
edges of vertex $n$ have been added to the graph, the $(k+1)$th edge of
vertex $n$ is attached to vertex $j\in{[n-1]}$ with probability
  \begin{align}\label{eq:pa2}
      \frac{d_j(n-1)+\sum^k_{\ell=1}\tone[n\overset{\ell}{\to} j]+\rho}{2m(n-2)+k+(n-1)\rho},
  \end{align}
  where $n\overset{\ell}{{\to}} j$ is shorthand for the event that the
    $\ell$-th outgoing edge of vertex $n$ is attached to vertex $j$.
  The resulting graph $G_n$ is a preferential attachment graph with
  $n$ vertices with parameters~$m$ and $\rho$, and we denote its law by
  $\mathrm{PA}(n,m,{\rho})$. 
\end{definition}

The formulation of the sequential model in \cite{Berger2014}
is somewhat different, but is easily seen to be equivalent.
Note also that \cite{Berger2014} assume (in our notation)
$\rho\ge 0$, but in the formulation above, only $\rho>-m$  is needed.
The definition above is valid also for $m=1$ (in which case
the graph is a tree), but we do not consider this case in the present paper;
see Remark \ref{Rm=1} below for a further discussion.

Since \cite{Bollobas2001}
 proved that the degree sequence of a certain class of
preferential attachment models indeed has a power-law behaviour, many
other properties of the model above and its variants have been investigated
over the last two decades. These results include for example, vertex degrees, distance
and local weak convergences; and
we refer to the books \cite{vdh2017,vdh2024} for a comprehensive
overview.

In this paper, we study the number of vertices that can be reached from the lastly added vertex $n$ via a directed path in the preferential attachment graph. We refer to these vertices (including vertex $n$)
as the \emph{descendants} of $n$ and their count as $X^{(n)}$, even though
all of them (apart from vertex $n$ itself) are added to $G_n$ before $n$. The
problem was first considered 
in \cite[Exercise 7.2.2.3-371 and 372]{Knuth}
for a uniform attachment graph, where each vertex 
has $m\ge 2$ outgoing edges and the
endpoints of these edges are chosen uniformly among the
existing vertices. 
(\cite{Knuth} uses drawing without replacement, thus
avoiding multiple edges, but as shown in \cite{Janson2023}, this makes no
difference asymptotically.)
This uniform attachment version is studied in \cite{Janson2023},
where it is shown that as $n\to\infty$,
if $\nu=(m-1)/m$, then
$X^{(n)}/n^{\nu}$ converges in distribution, and the limiting distribution
is given by a product of a constant factor and the $(1-\nu)$-th power of a
$\GAMMA(m/(m-1),1)$ variable.
The main result of the present paper is that 
for the preferential
attachment graph defined above,  $X^{(n)}$ behaves similarly, but
with a different exponent $\nu$ which furthermore depends on both $m$ and
$\rho$.

As in previous works such as \cite{Berger2014, Mori2003, PPR2017}, the
analysis in this work is hinged on a connection between P\'olya urns and the
preferential attachment mechanism. We use, in particular, the P\'olya urn representation of \cite{Berger2014} that was originally devised to study the local weak limit of preferential attachment graphs. As we show later, this representation result enables us to adapt the framework of \cite{Janson2023} to study the problem in the preferential attachment setting. 

We state our main results in the next subsection.

\subsection{Main results}
The parameters $m\ge2$ and $\rho>-m$ are fixed throughout the paper.
We define
\begin{align}\label{de:nu}
    \nu := 
\frac{(m-1)(m+\rho)}{m(m+\rho+1)}
\in(0,1)
.\end{align}
The proofs of the results below are developed in \refSs{se:pu}--\ref{Smom},
and as by-products of the proofs, we also prove some results on the
structure of the subgraph of descendants of $n$. 
In \refS{Sloop} we show that the following results hold also
for a preferential attachment model with possible self-loops.

\begin{theorem}\label{Tmain}
As \ntoo,
  \begin{align}\label{tmain}
    n^{-\nu} X
\dto 
\frac{\G\bigpar{\frac{(m-1)(m+\rho)}{m(m+\rho+1)}}
      \G\bigpar{\frac{m+\rho}{m(m+\rho+1)}+1}}
 {\G\bigpar{\frac{m+\rho}{m+\rho+1}}} 
\bbclr{\frac{(m+\rho+1)(m-1)}{2m+\rho}\xi_1 }^{1-\nu},
\end{align}
where $\xi_1\in\GAMMA(m/(m-1),1)$.
\end{theorem}

\begin{theorem}\label{Tmom}
  All moments converge in \eqref{tmain}. In other words, for any $p>0$,
 as \ntoo, 
\begin{align}\label{tmom}
\E[X^p]/n^{p\nu}
&\to
\lrpar{\frac{\G\bigpar{\frac{(m-1)(m+\rho)}{m(m+\rho+1)}}
      \G\bigpar{\frac{m+\rho}{m(m+\rho+1)}+1}}
 {\G\bigpar{\frac{m+\rho}{m+\rho+1}}} 
\lrpar{\frac{(m+\rho+1)(m-1)}{2m+\rho}}^{1-\nu}}^p
\notag\\&
\hskip4em\cdot
\frac{\gG(p(1-\nu)+\frac{m}{m-1})}{\gG(\frac{m}{m-1})}    
.  \end{align}
\end{theorem}

\begin{remark}
  In the special case $\rho=0$, \eqref{de:nu} and \eqref{tmain} simplify to
$\nu=(m-1)/(m+1)$ and 
\begin{align}\label{nov6}
    n^{-\nu} X
\dto \frac{1}{m+1}
\frac{\G\bigpar{\frac{m-1}{m+1}}
      \G\bigpar{\frac{1}{m+1}}}
 {\G\bigpar{\frac{m}{m+1}}} 
\bbclr{\frac{m^2-1}{2m}\xi_1 }^{2/(m+1)}.
\end{align}
If we specialize further to
the  case $m=2$ and $\rho=0$, we get $\nu=1/3$, 
and \eqref{tmain}  simplifies further to
  \begin{align}\label{nov7}
    n^{-1/3} X
\dto 
\frac{\G\bigpar{\frac13}^2}
 {2^{4/3}3^{1/3}
\G\bigpar{\frac23}} 
\xi_1 ^{2/3}
=
\frac{3^{1/6}\G\bigpar{\frac13}^3}
 {2^{7/3}\pi}
\xi_1 ^{2/3}
,\end{align}
with $\xi_1\in\GAMMA(2,1)$ and
    \begin{align}\label{RX}
\frac{\G\bigpar{\frac13}^2}
 {2^{4/3}3^{1/3}
\G\bigpar{\frac23}} 
\doteq  1.45833.
    \end{align}
In this case, \eqref{tmom} yields, for example,
  \begin{align}\label{tmom20}
\E[X]/ n^{1/3}
\to 
\frac{\G\bigpar{\frac13}^2}
 {2^{4/3}3^{1/3}
\G\bigpar{\frac23}} 
\gG\bigpar{2+\tfrac23}
=
\frac{5\,\G\bigpar{\frac13}^2}
 {2^{1/3}3^{7/3}} 
\doteq
2.19416 
.\end{align}
\end{remark}

\begin{remark}\label{Rm=1}
Definition \ref{de:pa} is valid also for $m=1$, and then defines a random
tree; such  preferential attachment trees have been studied by many authors.
In this case, $X\nn$ equals 1 + the depth of vertex $n$, and it is known
that
$X\nn$ grows like $\log n$,
in contrast to the case $m\ge2$ studied in the present paper,
where we show that $X\nn$ grows as a power of $n$.
More precisely,
as \ntoo,
\begin{align}\label{m=1}
  X\nn/\log n \pto \frac{1+\rho}{2+\rho},
\end{align}
and precise results are known 
on the exact distribution, 
Poisson approximation, 
and a central limit theorem, see
\cite{Dobrow-Smythe1996},
\cite[Theorem 6]{PP2007},
and \cite[Theorem 3]{Kuba-Wagner2010}.
(Papers on preferential attachment trees usually use a slightly different
definition, where the attachment probabilities depend on the outdegree
rather than the degree as in \eqref{de:pa}; apart from a shift in the
parameter $\rho$, this makes a difference only at the root. 
This minor difference ought not to affect asymptotic result;
for $X\nn$ this follows rigorously by the bijection in \cite{Kuba-Wagner2010}
which yields both exact and asymptotic results,
and in particular \eqref{m=1}, by straightforward calculations for both
versions of the definition.)
\end{remark}

We mention also an open problem, which we have not studied, where the same
methods might be useful.
\begin{problem}
Study the asymptotic behaviour of $\max\{X^{(n+1)},\dots,X^{(n+i)}\}$
for a fixed $i\ge2$, 
in both uniform and preferential attachment graphs.
Perhaps also do the same for $i=i(n)$ growing with $n$ at some rate.
\end{problem}

\subsection{Notation} \label{Snot}
As above,
$k\overset{\ell}{\to}i$ (where $1\le i <k\le n$ and $\ell\in[m]$)
denotes that in $G_n$
the $\ell$-th outgoing edge of vertex $k$ is attached to vertex $i$.
We say that
vertex $i$ is a \emph{child} of vertex $k$
if there is such an edge. 

As usual, empty sums are 0, and empty products are 1.

Convergence in distribution, in probability, and a.s.\ (almost surely) are
denoted by $\dto$, $\pto$, and $\asto$, respectively. 
Equality in distribution is
denoted by $\eqd$, and w.h.p. (with high probability) is short for 
``with probability tending one as $n\toinf$''.  

We frequently use two standard probability distributions. 
The $\GAMMA(a,b)$ distribution, with $a,b>0$,
has density $\G(a)^{-1}b^{-a} x^{a-1} e^{-x/b}$ on $(0,\infty)$.
The Beta$(a,b)$ distribution,  with $a,b>0$, has density
$\frac{\gG(a+b)}{\gG(a)\gG(b)}x^{a-1}(1-x)^{b-1}$ on $(0,1)$.

Most quantities defined below depend on $n$. We sometimes indicate
this by a superscript ${}\nn$, but usually we omit this to simplify the
notation. 
We may in proofs sometimes tacitly assume that $n$ is large enough.

$C[a,b]$, $C[0,\infty)$ and $C(0,\infty)$ denote the spaces of continuous
functions on the indicated intervals, equipped with the topology of
uniform convergence on compact subsets.
These spaces are complete separable metric spaces.
Note that a sequence of random functions in $C[0,\infty)$ or $C(0,\infty)$
converges (a.s., in probability, or in distribution)
if and only if it converges in the same sense in $C[a,b]$ for each compact
interval $[a,b]$ in $[0,\infty)$ or $(0,\infty)$, respectively.
(For $C[0,\infty)$ it is obviously equivalent to consider intervals $[0,b]$
only.) The case $C[0,\infty)$ is treated in detail in \cite{Whitt1970}; the
case $C(0,\infty)$ 
is similar.

$C$ denotes positive constants (not depending on $n$) that may vary from one
occasion to another. 
The constants may depend on the parameters $m$ and $\rho$;
we 
indicate dependence on other parameters (if any) by
writing e.g.\ $C_a$.

\section{P\'olya urn representation}\label{se:pu}
We shall use a celebrated result of \cite{Berger2014},
which states that the dynamics of the preferential attachment graph can be
encoded in a collection of classical P\'olya urns; see also \cite[Chapter 5]{vdh2024} for more details. In a classical P\'olya
urn with initially $a$ red balls and $b$ black balls, a ball is randomly
sampled from the urn at each {step}, and is then returned to the urn with
another ball of the same colour. 
(The ``numbers'' of balls are not necessarily integers; any positive real
numbers are allowed.)
In the preferential attachment graph,
for each $i\ge2$,
  the weight of vertex $i$, defined as the degree + $\rho$,
and the total weight of the first $i-1$ vertices
  evolve like the numbers of red and black balls in a classical P\'olya
  urn. The initial numbers of red and black balls are {$a=m+\rho$ and $b=(2i-3)m+(i-1)\rho$},
  which are the weights of vertex $i$ and the first $i-1$ vertices before
  the edges of vertex $i+1$ are added to the graph. When one of the first
$i$ vertices is chosen as a recipient of a newly added edge, the number of
red balls in the urn increases by one if vertex $i$ is the recipient;
otherwise we add a new black ball to the urn. 
It is well-known, for example as a consequence of exchangeability and de
Finetti's  theorem, that the proportion of red balls a.s.\ converges to 
a random number $\beta\in \mathrm{Beta}(a,b)$,
and that conditioned on $\beta$,
the indicators that a red ball is chosen at each step are distributed as
conditionally independent Bernoulli variables with parameter
$\beta$. Consequently, by conditioning on suitable beta variables,
  the preferential attachment graph can instead be generated using
  independent steps.

The model and the theorem below are easy variations of their counterparts in \cite[Section 2.2]{Berger2014}. The only difference is that $\rho$ is allowed to be negative here.

\begin{definition}[P\'olya urn representation, \cite{Berger2014}]\label{de:PUR}
Given the integer $m\geq 2$ {and the real number $\rho>-m$}, let
$(B_j)^{\infty}_{j =1}$ be independent random variables such that $B_1=1$ and
\begin{align}\label{de:betas}
    B_j\in\mathrm{Beta}(m+\rho, (2j-3)m+(j-1)\rho),
\qquad j\ge2.
\end{align}
Given $(B_j)^{\infty}_{j =1}$, 
construct for each $n\ge1$ 
a (directed) graph $G_n$ on $n$ vertices (labelled by $[n]$)
such that each vertex $2\leq k\leq n$ has $m$ outgoing edges, 
and the recipient of each outgoing edge of $k$ is $i\in[k-1]$ 
with probability 
\begin{align}\label{pb}
  B_i\prod_{j=i+1}^{k-1}(1-B_j),
\end{align}
with the endpoints of all edges in $G_n$ chosen (conditionally)  independently.
The law of $G_n$ is denoted by $\mathrm{PU}(n,m,{\rho}),$ where PU is short for
P\'olya Urn. 
\end{definition}

\begin{remark}\label{Rstop}
  The probabilities \eqref{pb} can be interpreted as follows, which will be
  useful below:
Given $(B_j)_{j=1}^\infty$, each edge from vertex $k$ tries to land at $k-1$, $k-2$,
\dots{} successively; at each vertex $j$ it stops with probability $B_j$,
and otherwise it continues to the next vertex. (All random choices are
independent, given $(B_j)_{j=1}^\infty$.)
\end{remark}

\begin{remark}\label{RBerger}
  The construction in \cite{Berger2014} is actually
  formulated in the following somewhat different way, 
which  obviously is equivalent; we will use this version too below.
Define 
\begin{align}\label{de:S}
    S_{n,j}= \prod^{n-1}_{i=j+1} (1-B_i)\qquad\text{for $0\leq j\leq n-1$}.
\end{align}
(In particular,  $S_{n,0}=0$ and $S_{n,n-1}=1$.)
Conditioned on $(B_j)^{n-1}_{j=2}$, let $(U_{k,\ell})^{n,m}_{k=2,\ell=1}$ be
independent random variables with 
\begin{equation}\label{de:Uij}
    U_{k,\ell}\in \sU[0,S_{n,k-1}).
\end{equation}
For each vertex $2\le k\le n$, add the $m$ outgoing edges such that
\begin{align}\label{kil}
  k \overset{\ell}\to i 
\iff
U_{k,\ell} \in [S_{n,i-1},S_{n,i}), \qquad \ell\in[m],\; i\in[k-1].
\end{align}
Note also that a natural way to achieve \eqref{de:Uij} is to
let $(\widetilde U_{k,\ell})^{n-1,m}_{k=2,\ell=1}$ be independent
$\sU[0,1]$ variables,
independent of $(B_i)^{n-1}_{i=2}$, 
and set
\begin{align}\label{eq:sU}
    U_{k,\ell} := S_{n,k-1} \widetilde U_{k,\ell}.
\end{align}
\end{remark}

\begin{theorem}[\cite{Berger2014}, Theorem 2.1]
    For all integers $n\geq2$, $m\geq 2$ and real $\rho>-m$, $\mathrm{PA}(n,m, {\rho})=\mathrm{PU}(n,m,{\rho})$.
\end{theorem}

In view of this theorem, it is enough to consider the P\'olya urn
representation instead of the preferential attachment graph. We shall do so
in the subsequent analysis and always have $G_n\in \mathrm{PU}(n,m,
{\rho})$.

\begin{remark}\label{Runiform}
  The uniform directed acyclic graph studied in \cite{Janson2023},
where each new edge from $k$ is attached uniformly to a vertex in $[k-1]$,
can be seen as the limit as $\rho\to\infty$ of the construction above;
it can be constructed by the same procedure, except that we let $B_j:=1/j$
(deterministically). This may help seeing the similarities and differences 
in the arguments below and in \cite{Janson2023}.
Not surprisingly, 
formally taking the limit $\rho\to\infty$ in \eqref{tmain}
yields 
the main result of \cite{Janson2023}.
\end{remark}

\begin{remark}
  Unless we say otherwise, we use the same sequence $(B_i)_{i=1}^\infty$ for
  every $n$. (But see \refS{Sconv} for an exception.) 
\end{remark}

\section{Preliminaries}\label{Sprel}

For convenience, we define
the positive constants
\begin{align}
    \theta&:=2m+\rho,\label{de:theta}
\\\label{de:chi}
    \chi &:= \frac{m+\rho}{2m+\rho}
= \frac{m+\rho}{\theta},
\end{align}
noting that if $\rho=0$, then $\chi=1/2$ for any $m$.

Recall that if $B\in\Beta(a,b)$, it follows by evaluating a beta integral that the moments are given by
\begin{align}\label{betamom}
  \E B^s = \frac{\gG(a+b)\gG(a+s)}{\gG(a)\gG(a+b+s)}
=\frac{\gG(a+s)/\gG(a)}{\gG(a+b+s)/\gG(a+b)},
\qquad s>0.
\end{align}
Recall also that for any fixed real (or complex) $a$ and $b$, 
and $x>0$ (with $x+a\notin\set{0,-1,-2,\dots}$,
\begin{align}\label{gg}
\frac{ \gG(x+a)}{\gG(x+b)} = x^{a-b}\bigpar{1+O\bigpar{x\qw}},
\end{align}
which follows readily from Stirling's formula;
see also \cite[5.11.13]{NIST}.

Similarly to the definition of $S_{n,j}$ in \eqref{de:S},
we also define 
\begin{align}\label{de:phi}
    \Phi_k =\prod^k_{j=1} (1+(m-1)B_j), \quad \text{for $k\geq 0$}
.\end{align}

We collect here some simple results for these variables that will be used
later.

\begin{lemma}\label{LB1}
    For $2\leq i<\infty$, let $B_i$ be as in \eqref{de:betas}; and for $1\le
    i<\infty$, let $\Phi_i$ be as in \eqref{de:phi}. 
We then have for $2\le i<\infty$, 
\begin{gather*}
        \E( B_i) = \frac{m+\rho}{\theta i-2m}
=\frac{\chi}{i-2m/\gth}
=\frac{\chi}{i}+O\Bigpar{\frac{1}{i^2}},
\numberthis\label{eq:Bmean}\\
        \E(B_i^2) = \frac{(m+\rho+1)(m+\rho)}{(\theta i-2m+1)(\theta i-2m)} 
 =O\Bigpar{\frac{1}{i^2}}.
\numberthis \label{eq:Ebetasq}
    \end{gather*}
    Furthermore, for $2\leq j\leq k<\infty$, 
    \begin{align*}
        \prod^k_{i=j}\E(1+(m-1)B_i)
          &= \frac{\G\bclr{k+1+[(m-1)(m+\rho)-2m]/\theta}\G(j-2m/\theta)}{\G(k+1-2m/\theta)\G\bclr{j+[(m-1)(m+\rho)-2m]/\theta}}\notag\\
        &= \bbclr{\frac{k}{j}}^{(m-1)\chi} \bclr{1+O\bigpar{j^{-1}}}\numberthis \label{eq:mprodbeta}
    \end{align*}
    and
    \begin{align*}
          \E \Phi_k = \frac{m\cdot \G\bclr{2-2m/\theta}}{\G\big(2+[(m-1)(m+\rho)-2m]/\theta\big)} k^{(m-1)\chi} \big(1+O\bigpar{k^{-1}}\big).
        \numberthis\label{eq:meanphik}
    \end{align*}
    Finally, there is a positive constant $C$  such that, for $2\leq j\leq k<\infty$,
    \begin{align}\label{eq:trunphi2bd}
     \prod^k_{i=j}\E(1+(m-1)B_i)^2 \leq C\bbclr{\frac{k}{j}}^{2(m-1)\chi} 
    \end{align}
    and, for $2\le k<\infty$,
    \begin{align}\label{eq:phibd}
        \E(\Phi_k^2)\leq Ck^{2(m-1)\chi}, 
\quad \E(\Phi^{-1}_k)\leq Ck^{-(m-1)\chi},
\quad  \E(\Phi^{-2}_k)\leq Ck^{-2(m-1)\chi}.
    \end{align}
\end{lemma}

\begin{remark}
    If $m=2$ and $\rho=0$, then $\theta=4$ and so in \eqref{eq:meanphik},
    \begin{align}
        \frac{m\cdot \G\bclr{2-2m/\theta}}{\G\big(2+[(m-1)(m+\rho)-2m]/\theta\big)} = \frac{2}{\G(3/2)} = \frac{4}{\sqrt{\pi}}.
    \end{align}
\end{remark}

\begin{proof}
    The equalities in \eqref{eq:Bmean} and \eqref{eq:Ebetasq} follow from 
\eqref{betamom}, recalling \eqref{de:theta}--\eqref{de:chi}.
For  \eqref{eq:mprodbeta},
    we use \eqref{eq:Bmean} {to obtain}
\begin{align}\label{eq:Eprodbeta}
        \prod^k_{i=j}\E(1+(m-1)B_i)  &= \prod^k_{i=j}\frac{i+[(m-1)(m+\rho)-2m]/\theta}{i-2m/\theta}\notag\\
       &= \frac{\G\bclr{k+1+[(m-1)(m+\rho)-2m]/\theta}\G(j-2m/\theta)}{\G(k+1-2m/\theta)\G\bclr{j+[(m-1)(m+\rho)-2m]/\theta}},
    \end{align}
    and so \eqref{eq:mprodbeta} follows 
from \eqref{gg} and \eqref{de:chi}.
The formula  \eqref{eq:meanphik} follows 
similarly by taking $j=2$ in \eqref{eq:Eprodbeta} and using $B_1=1$.
    
To prove \eqref{eq:trunphi2bd},  we write for $i\geq 2$,
using \eqref{eq:Bmean}--\eqref{eq:Ebetasq},
\begin{align*}
    \E(1+(m-1)B_i)^2 &= 1 + 2(m-1)\E B_i + (m-1)^2\E B_i^2\\
&=1+2(m-1)\frac{\chi}{i}+O\bigpar{i^{-2}}\\
    &=: 1 + y_i.
\end{align*}
Taking logarithms,
    \begin{align*}
        \log\biggpar{\prod^{{k}}_{i={j}} (1+y_i)} &= \sum^{k}_{i={j}} \log(1+y_i) 
        \leq \sum^{k}_{i=j} {y_i}
        = 2(m-1)\chi\sum^{k}_{i=j} \big(i^{-1}+O(i^{-2})\big)\\
        &= 2(m-1)\chi\log\bbclr{\frac{k}{j}} + O(j^{-1}).\numberthis\label{eq:Ebetasq2}
    \end{align*}
     This implies the inequality in \eqref{eq:trunphi2bd}. The bound on 
     $\E(\Phi^2_k)$ in \eqref{eq:phibd} follows from the definition in
     \eqref{de:phi} and applying \eqref{eq:trunphi2bd} with $j=2$. The upper
     bound on $\E(\Phi_k^{-2})$ in \eqref{eq:phibd} can be proved similarly,
     where we can use $(1+x)^{-2}\leq 1-2x+3x^2$ for $x\geq 0$, and thus
\begin{align}\label{dx0}
    \E\bigsqpar{(1+(m-1)B_i)^{-2}} &
\leq 1-2(m-1)\E B_i+3(m-1)^2\E B_i^2 
\notag\\&
= 1- 2(m-1)\chi i^{-1} + O(i^{-2}),
\end{align}
together with $\log(1-x)\leq -x$. 
Finally, by the Cauchy--Schwarz inequality
and the just proven $\E(\Phi_k^{-2})\leq Ck^{-2(m-1)\chi}$, 
\begin{align}
    \E(\Phi_k^{-1})\leq \sqrt{\E(\Phi_k^{-2})}\leq Ck^{-(m-1)\chi},
\end{align}
which completes the proof of all three inequalities claimed in \eqref{eq:phibd}. 
     \end{proof}

\subsection{An infinite product}


\begin{lemma}  \label{LB2}
The infinite product
    \begin{equation}\label{lb2a}
        \beta:=\prod^\infty_{k=1}\frac{1+(m-1)B_k}{\E(1+(m-1)B_k)}
=\lim_\ktoo\frac{\Phi_k}{\E\Phi_k}
    \end{equation}
exists a.s.\ and in $L^p$ for every $p<\infty$.
Furthermore, $\E\beta=1$ and $\beta>0$ a.s.

We have also, as $k\to\infty$,
\begin{align}\label{lb2b}
  k^{-(m-1)\chi}\Phi_k\asto 
\tgb:= \frac{m\cdot \G\bclr{2-2m/\theta}}{\G\big(2+[(m-1)(m+\rho)-2m]/\theta\big)}\gb.
\end{align}
\end{lemma}

\begin{proof}
Define for $k\ge1$
    \begin{align}\label{eq:betamg}
        \tM_k := \frac{\Phi_k}{\E (\Phi_k)} 
= \prod^k_{i=1} \frac{1+(m-1)B_i}{\E(1+(m-1)B_i)}.
    \end{align}
This is a product of independent random variables with mean 1, and thus a
martingale.
For every fixed integer $r>1$, we have by the binomial theorem, 
$|B_k|\le1$,
and \eqref{eq:Bmean}--\eqref{eq:Ebetasq},
\begin{align}\label{freja}
  \E(1+(m-1)B_k)^r
&= \sum_{j=0}^r\binom{r}j{(m-1)^j}\E B_k^j
= 1 + r(m-1) \E B_k + O(\E B_k^2) \notag\\
&= 1 + r(m-1)\E B_k + O(k^{-2})\notag\\
&= (1 + (m-1)\E B_k)^r + O(k^{-2}).
\end{align}
Hence, for every $k\ge1$,
\begin{align}
  \E \tM_k^r = \prod_{i=1}^k\frac{\E(1+(m-1)B_i)^r}{(\E(1+(m-1)B_i))^r}
= \prod_{i=1}^k\bigpar{1+O(i^{-2})} \le C_r
\end{align}
and thus the martingale $\tM_k$ is $L^r$-bounded;
consequently it converges
in $L^r$, and thus in $L^p$ for all real $0<p\le r$.
Since $r$ is arbitrary, this holds for all $p>0$.

In particular, $\tM_k\to\gb$ in $L^1$, which shows that 
$\E\gb=\lim_\ktoo\E\tM_k=1$.

The event $\set{\gb=0}$ is independent of any finite number of
$B_1,B_2,\dots$, and is thus a tail event. The Kolmogorov zero-one law, 
see e.g.\ \cite[Theorem 1.5.1]{Gut}, 
thus shows that $\P(\gb=0)=0$ or 1, but $\P(\gb=0)=1$ is impossible since
$\E\gb=1$. Hence, $\gb>0$ a.s.

Finally, \eqref{lb2b} follows by \eqref{lb2a} and \eqref{eq:meanphik}.
\end{proof}

\subsection{Estimates for $S_{n,k}$}
Below, let $\psi_n$ be a positive function such that
$\psi_n\le n-1$ and 
$\psi_n\to\infty$ as $n\to\infty$ (we later choose $\psi_n=n/\log n$). The next lemma shows that w.h.p., for all $k\ge \psi_n$, the random variables $S_{n,k}$ are close enough to the constants $(k/n)^\chi$. 

\begin{lemma}\label{le:Sest}
    Let $S_{n,k}$ be as in \eqref{de:S} and $\psi_n$ be as above. Define
    $\delta_n = \psi_n^{-\eps}$ for some $\eps\in(0,1/2)$. 
Then, there is a positive constant $C$  such that 
    \begin{align}
        \IP\bigg[\max_{\ceil{\psi_n}\leq k<n}\bigg|S_{n,k}-\bbclr{\frac{k}{n}}^\chi\bigg|\geq 2\delta_n\bigg]\leq {C}{\psi_n}^{2\eps-1}.
    \end{align}
\end{lemma}

The proof of Lemma \ref{le:Sest} is based on a standard
  martingale argument that is similar to \cite{lo2024} (see also
  \cite{Berger2014}), but we present it here for
  completeness. To prepare for the main proof, we start by estimating
$\E(S_{n,k})$. 

\begin{lemma}\label{le:smom}
    Let $S_{n,k}$ be as in \eqref{de:S}. For every $1\leq k\leq n-1$, we have
    \begin{align}\label{eq:meanS}
        \E(S_{n,k}) = \frac{\G\bclr{n-(3m+\rho)/\theta}}{\G\bclr{n-2m/\theta}}\frac{\G\bclr{k+1-2m/\theta}}{\G\bclr{ k+1-(3m+\rho)/\theta}}.
    \end{align}
\end{lemma}

\begin{proof}
    Recalling that $(B_j)^{n-1}_{j=2}$ are independent,
we obtain from \eqref{eq:Bmean}
      \begin{align*}
        \E S_{n,k} &=\prod^{n-1}_{j=k+1} \E(1-B_i) = \prod^{n-1}_{j=k+1} \frac{\theta j -3m-\rho}{\theta j -2m} = \prod^{n-1}_{j=k+1} \frac{j-(3m+\rho)/\theta}{j-2m/\theta} \\&= \frac{\G\bclr{n-(3m+\rho)/\theta}}{\G\bclr{n-2m/\theta}}\frac{\G\bclr{k+1-2m/\theta}}{\G\bclr{ k+1-(3m+\rho)/\theta}},\numberthis
    \end{align*}
    as claimed in the lemma.
\end{proof}

\begin{lemma}\label{le:smean}
    Let $S_{n,k}$ be as in \eqref{de:S}.
    Then,  
there is a positive constant $C$ such that,
for $1\leq k \leq n-1$,
    \begin{align}\label{eq:smean}
        \bigg|\E(S_{n,k})-\bbclr{\frac{k}{n}}^{\chi}\bigg|\leq  {\frac{C}{n^\chi k^{1-\chi}}.}
    \end{align}
\end{lemma}

\begin{proof}
  By \eqref{eq:meanS} and \eqref{gg}, 
we have, recalling \eqref{de:chi},
  \begin{align}
\E(S_{n,k})= n^{-\chi}\bigpar{1+O(n\qw)}k^{\chi}\bigpar{1+O(k\qw)}
=\bbclr{\frac{k}{n}}^{\chi}\bigpar{1+O(k\qw)}
,  \end{align}
which yields \eqref{eq:smean}.
\end{proof}

\begin{proof}[Proof of Lemma \ref{le:Sest}]
 By \refL{le:smean}, for $n$ large enough and any $k\in[\psi_n, n]$,
\begin{align}
  \biggabs{\E(S_{n,k})-\bbclr{\frac{k}{n}}^\chi}
\le \frac{C}{{n^\chi k^{1-\chi}}}
\le \frac{C}{k} \le \frac{C}{\psi_n}<\gd_n.
\end{align}
Hence,
    \begin{align*}
    \IP\bigg[\max_{\ceil{\psi_n}\leq k{<} n}\bigg|S_{n,k}-\bbclr{\frac{k}{n}}^\chi\bigg|\geq 2\delta_n\bigg] 
        &\leq 
       \IP\bigg[\max_{\ceil{\psi_n}\leq k{<} n}\big|S_{n,k}-\E(S_{n,k})\big| \geq \delta_n\bigg] \\
       &\leq \IP\bigg[\max_{\ceil{\psi_n}\leq k{<} n}\big|S_{n,k}/\E(S_{n,k})-1\big| \geq \delta_n\bigg],
\numberthis\label{hw1}
    \end{align*}
    noting that $\E S_{n,k}\leq 1$ for $k\geq 1$. To bound the right-hand side of \eqref{hw1}, we first observe that for $k\geq {0}$,
    \begin{align}
        \MM_k := \prod^{{n-1}}_{{j=n-k}} \frac{1-B_{j}}{\E(1-B_{j})} = \frac{S_{{n,n-1-k}}}{\E S_{n,n-1-k} }
    \end{align}
    is a martingale with respect to {the} $\sigma$-algebras
    generated by {$(B_j)^{n-1}_{j=n-k}$}, with $\E \MM_k=1$. 
Now, by Doob's inequality for the submartingale
    $(\MM_k-1)^2$, see e.g.\ \cite[Theorem 10.9.1]{Gut},
    \begin{align*}\label{eq:mdoob}
        \IP\bigg[\max_{\ceil{\psi_n}\leq k<n}\big|S_{n,k}/\E(S_{n,k})-1\big| \geq \delta_n\bigg] &= {\IP\bigg[\max_{0\leq k\leq n-1-\ceil{\psi_n}}|\MM_k-1| \geq \delta_n\bigg]} \\
        &\leq \delta_n^{-2}\var\big(\MM_{n-1-\ceil{\psi_n}}\big). 
\numberthis
    \end{align*}
Using $\E \MM_{n-1-\ceil{\psi_n}}=1$ and the independence of the beta variables, we have
    \begin{align}
        \var\big(\MM_{n-1-\ceil{\psi_n}}\big) &= \E\big(\MM_{n-1-\ceil{\psi_n}}^2\big) -1 = \prod^{n-1}_{k = \ceil{\psi_n}+1}\frac{\E[(1-B_k)^2]}{(\E[1-B_k])^2}-1,
    \end{align}
    and by \eqref{de:betas}, \eqref{betamom}, and simplifying, we get
    \begin{align*}
        \var\big(\MM_{n-1-\ceil{\psi_n}}\big) &= \prod^{n-1}_{k=\ceil{\psi_n}+1 } {\bbclr{\frac{\theta k-3m-\rho+1}{\theta k-3m-\rho}\cdot \frac{\theta k-2m}{\theta k-2m+1}} }-1\\
        &\leq \prod^{n-1}_{k={\ceil{\psi_n}+1} }(1+Ck^{-2}) -1 \leq \frac{C}{\psi_n}. \numberthis\label{eq:mvar}
    \end{align*} 
    Applying \eqref{eq:mvar} {and $\delta_n=\psi_n^{-\eps}$} to 
\eqref{hw1} and \eqref{eq:mdoob} yields
    \begin{align*}
 \IP\bigg[\max_{\ceil{\psi_n}\leq k<n}\bigg|S_{n,k}-\bbclr{\frac{k}{n}}^\chi\bigg|\geq 2\delta_n\bigg] 
\leq \gd_n^{-2}\frac{C}{\psi_n}
= C \psi_n^{2\eps-1}, 
    \end{align*}
    hence proving the lemma.
\end{proof}

\subsection{Asymptotics of two sums}

Fix a sequence $\gl_n\to\infty$ (we will later choose 
$\gl_n=n^\nu$)
and define, for $y\ge0$ and $0\le k\le \ell<\infty$,
\begin{align}\label{lh1}
  H^y_{k,\ell}&:=\sum_{i=k+1}^\ell\bigsqpar{(1-B_i)^{\gln y}-1+\gln yB_i}, \qquad  I^y_{k,\ell} := \sum_{i=k+1}^\ell \bigsqpar{1-(1-B_i)^{\gln y}}
\end{align}
and, 
for $y\ge0$ and $0\le s\le t<\infty$,
\begin{align}\label{lh2}
  \hH^y_{s,t}:=\gln\qw H^y_{\floor{s\gln},\floor{t\gln}},\qquad  \wh I^y_{s,t}:=\gln\qw I^y_{\floor{s\gln},\floor{t\gln}}.
\end{align}

\begin{lemma}\label{LH}
Let\/ $0<s\le t<\infty$. Then, for every $y\ge0$,
as \ntoo,
\begin{align}\label{lh3}
 \hH^y_{s,t}&\pto 
\int_s^t \Bigpar{\frac{(\theta u)^{m+\rho}}
  {(\theta u+y)^{m+\rho}}-1+\frac{\chi y}{u}}\dd u 
=
\int_s^t \Bigpar{\frac{1}{(1+y/(\theta u))^{m+\rho}}-1+\frac{\chi y}{u}}\dd u 
\end{align}
and
\begin{align}\label{lh34}
 \wh I ^y_{s,t}\pto 
\int_s^t \Bigpar{1-\frac{(\theta u)^{m+\rho}}{(\theta u+y)^{m+\rho}}}\dd u
=
\int_s^t \Bigpar{1-\frac{1}{(1+y/(\theta u))^{m+\rho}}}\dd u.
\end{align}
\end{lemma}

\begin{proof}
Denote the summand in \eqref{lh1} by $\gD H_i$.
Then $-1\le \gD H_i\le \gln y B_i$,
and thus, using \eqref{eq:Ebetasq}, 
for $\floor{s\gln}<i\le\floor{\gln t}$,
\begin{align}\label{lh4}
  \Var (\gD H_i) 
\le   \E (\gD H_i)^2 
\le C + C\gln^2\E B_i^2 \le C_s.
\end{align}
The summands $\gD H_i$ are independent, and thus \eqref{lh1}--\eqref{lh2}
and \eqref{lh4} 
yield
\begin{align}\label{lh5}
  \Var(\hH^y_{s,t})
=\gln^{-2}\sum_{i=\floor{s\gln}+1}^{\floor{t\gln}}\Var\bigpar{\gD H_i}
\le {C_{s,t}} \gln^{-1} =o(1).
\end{align}
Hence, it suffices to show that the expectation $\E \hH^y_{s,t}$ converges
to the limit in \eqref{lh3}.
We have, applying \eqref{betamom} to {$1-B_i\in\Beta(\theta i-3m-\rho,m+\rho)$} and using \eqref{gg},
uniformly for $s\gln < i \le t\gln$,
\begin{align}\label{lh6}
\E(1-B_i)^{\gln y}&
=
\frac{\gG(\theta i-2m)\gG(\theta i-3m-\rho+\gln y)}{\gG(\theta i-3m-\rho)\gG(\theta i-2m+\gln y)}
= \bbclr{\frac{\theta i}{\theta i +\gln y}}^{m+\rho} +o(1)\notag\\
&= \bbclr{\frac{\theta i/\gln}{\theta i/\gln + y}}^{m+\rho} +o(1). 
\end{align}
Hence, using also \eqref{eq:Bmean}, if $i=u\gln$ with $u\in(s,t]$,
\begin{align}\label{lh7}
  \E[\gD H_i]
&=
\E(1-B_i)^{\gln y}-1+\gln y \E B_i
\notag\\ 
&=\bbclr{\frac{\theta i/\gln}{\theta i/\gln + y}}^{m+\rho} 
-1 + \frac{\chi}{i}\gln y +o(1). 
\end{align}
It follows that  $\E\hH^y_{s,t}$ is $o(1)$ plus a Riemann sum of the integral in
\eqref{lh3}. 

The proof of \eqref{lh34} is similar, where we now replace $\Delta H_i$ with 
\begin{align}
    \Delta  I_i = 1- (1-B_i)^{\lambda_n y}
\end{align}
and, for $s\lambda_n <i\leq t\lambda_n$, use the estimates
$0\le\gD I_i\le1$ and thus, using \eqref{lh6},
\begin{align}
   & \var(   \Delta I_i) 
\le \E (\Delta I_i)^2
\le1 
,\\
   & \E[\gD I_i]
=
1-\E(1-B_i)^{\gln y} = 1- \bbclr{\frac{\theta i/\gln}{\theta i/\gln + y}}^{m+\rho} + o(1)
\end{align}
to proceed.
\end{proof}

\section{Basic analysis}\label{se:basicanalysis}

In this and subsequent sections, we follow the framework (and hence
the notation) in \cite{Janson2023}. 
To concentrate on the important aspects of the proof,  
we assume that $m=2$ and $\rho=0$;
note that then $\chi=\frac12$ and $\gth=2m=4$.
The minor modifications for the general case are discussed in \refS{Sgen}.

\subsection{The stochastic recursion.}\label{sse:sr}
Let $D_n$ 
be the  subgraph in $G_n$, consisting of vertex~$n$, all
vertices that can be reached from $n$ via a directed path, and all the edges
between them. 
We think of the vertices and edges in $D_n$ as coloured red.
We use the following stochastic recursion to construct $D_n$. 
It is similar to the recursion used in \cite{Janson2023},
with differences that  stem
  from the difference of the models.
\begin{enumerate}
    \item Sample the beta variables $(B_j)^{n-1}_{j=2}$ defined in
\eqref{de:betas}.
    \item Declare vertex $n$ to be red and all others black. Initiate
        the recursion by setting $k:=n$. 
    \item If vertex $k$ is red, choose the recipients of the two
      outgoing edges from vertex $k$ according to the construction given
        in Definition \ref{de:PUR}. 
    After sampling the recipients, declare them as red.\\
    If vertex $k$ is black, delete $k$ and do nothing else.
    \item If $k=2$ then stop; otherwise let $k:=k-1$ and repeat from {(3)}.
\end{enumerate}

For integers $0\leq k\leq n-1$, let $Y_k$ be the number of edges in
$D_n$ that start from $\clc{k+1,\dots,n}$ and end in
$\clc{1,\dots,k}$. Define $Z_k$ as the number of edges in $Y_k$ that
end in~$k$. Note that we have the boundary conditions $Y_{n-1}=2$ and
$Y_0=0$; as well as $Z_1=Y_1$ and $Z_0=0$.  

For $1\leq k\leq n-1,$ denote the indicator that at least one edge of $Y_k$
ends at $k$ as 
\begin{align}
    J_k=\tone[Z_k\geq 1], 
\end{align}
which is the same as the indicator that $k$ is red.
Thus summing $J_k$ over $k\in[n-1]$ gives the number $X\nn$ of red vertices. 
For $2\leq k\leq n-1$, 
the number of edges that start at $k$ is $2J_k$, 
and we thus have
\begin{align}\label{eq:stocrecur}
    Y_{k-1} = Y_k - Z_k + 2\cdot J_k = Y_k -Z_k + 2\cdot \tone[Z_k\geq 1].
\end{align}

As in \cite{Janson2023}, we use a modified version of the procedure above,
where we use the construction in \refR{Rstop}.
In (3) above, we thus do not choose the recipients of the outgoing edges, we
just note that they have endpoints in $[k-1]$.
We then at the next vertex toss a coin for each edge with unassigned endpoint
to decide whether it ends there or not. This yields the following equivalent
version of the construction.

\begin{enumerate}
    \item Sample the beta variables $(B_j)^{n-1}_{j=2}$ defined in
\eqref{de:betas}.
    \item Declare vertex $n$ to be red and all others black. 
Initiate the recursion by setting $k:=n$. 
    \item If vertex $k$ is red,
add two outgoing edges from vertex $k$,
with as yet undetermined endpoints in $[k-1]$;
mark these edges \emph{incomplete}.
\item  Let $k:=k-1$.
\item For each incomplete edge, toss
a coin with heads probability $B_k$, independently given $B_k$.
If the outcome is heads, the edge ends at $k$ and is marked complete;
furthermore,  vertex $k$ is coloured red.
Otherwise do nothing (so the edge is still incomplete). 
    \item If $k=1$ then stop; otherwise repeat from (3).
\end{enumerate}

Let
$\cF_k$ be the $\sigma$-field generated by 
all beta variables $(B_j)^{n-1}_{j=2}$ and 
the coin tosses at vertices $n-1,\dots,k+1$.
Then $\cF_1,\dots,\cF_{n-1}$ forms a decreasing sequence of
$\sigma$-fields, 
and $Y_{n-1},\dots, Y_k$ are measurable with respect to $\cF_k$.
Moreover,   
conditioned on $\cF_k$, we have
\begin{align}\label{eq:distZk}
    Z_k\mid \cF_k\in \mathrm{Bin}(Y_k, B_k)\quad{\text{for $1\leq k\leq n-1$.} }
\end{align}
Thus, in view of the recursion \eqref{eq:stocrecur}, we have for $2\leq
k\leq n-1$,
\begin{align}
 \E(Y_{k-1}\mid \cF_k)&= Y_k - \E(Z_k\mid \cF_k) + 2\cdot\IP(Z_k\geq 1\mid \cF_k)
\notag\\
    &=Y_k -B_k Y_k +2\bclr{1-(1-B_k)^{Y_k}}.\label{eq:exact}
\end{align}
By Markov's inequality, we also have {for $2\leq k\leq n-1$,}
\begin{align*}
    \E[Y_{k-1}\mid\cF_k]&\leq  Y_k - \E[Z_k\mid \cF_k] + 2\cdot\E[Z_k\mid \cF_k] = Y_k+\E[Z_k\mid \cF_k]\\
    &{=} (1+B_k)Y_k. \numberthis\label{eq:markov}
\end{align*}
Define, recalling \eqref{de:phi},
\begin{align}\label{de:Wk}
    W_k = \Phi_k  Y_k\quad  \text{{for $0\leq k\leq n-1$,}}
\end{align}
noting that $W_0=2Y_0=0$.
Using \eqref{eq:markov} and \eqref{de:phi}, we find {for $2\leq k\leq n-1$,}
\begin{align}\label{eq:Wk}
    \E(W_{k-1}\mid \cF_k) &= \Phi_{k-1}\E(Y_{k-1}\mid \cF_k)
   \leq \Phi_{k-1} (1+B_k) Y_k =\Phi_k Y_k = W_k; 
\end{align}
and so $W_0,\dots,W_{n-1}$ is a reverse supermartingale. The initial value is 
\begin{align}
    W_{n-1} = \Phi_{n-1} Y_{n-1}=2\Phi_{n-1}.
\end{align}
By Doob's decomposition, 
\begin{align}\label{eq:doobdecom}
    W_k = M_k-A_k, \quad 0\leq k \leq n-1,
\end{align}
where
\begin{align}\label{eq:mg}
    M_k := 2\Phi_{n-1} + \sum^{n-1}_{j=k+1} (W_{j-1} - \E(W_{j-1}\mid\cF_j))
\end{align}
is a reverse martingale and 
\begin{align}\label{eq:incp}
    A_k := \sum^{n-1}_{j=k+1} (W_j-\E(W_{j-1}\mid\cF_j))
\end{align}
is positive and reverse increasing. To see these properties of $A_k$, we note $A_{n-1}=0$ and by \eqref{eq:Wk},
\begin{align}\label{eq:Akdiff}
    A_{k-1} - A_k = W_k - \E(W_{k-1}\mid \cF_k) \geq 0\quad \text{for $1\leq k \leq n-1$.}
\end{align}
Hence, for $0\leq k\leq n-1$, we have $0\le W_k\leq M_k$. 

From the exact formula \eqref{eq:exact},
\begin{align}
    \E (W_{k-1}\mid \cF_k) = \Phi_{k-1}\E(Y_{k-1}\mid \cF_k) = \Phi_{k-1} (1-B_k) Y_k  + 2\Phi_{k-1}\bclr{1-(1-B_k)^{Y_k}},
\end{align}
and so \eqref{eq:Akdiff} can be written as
\begin{align}
    A_{k-1}-A_k &= 2B_k\Phi_{k-1}Y_k - 2\Phi_{k-1}\bclr{1-(1-B_k)^{Y_k}}
\notag\\
    &=2\Phi_{k-1} \bclr{(1-B_k)^{Y_k}-1+B_kY_k}.\label{eq:incA}
\end{align}
Following the steps in \cite[equation (2.16)]{Janson2023} for evaluating
$\var(Y_{k-1}\mid \cF_{k})$ in the {uniform attachment case}, 
here we have, for  $1\leq k\leq n-1$,
\begin{align}\label{nov1}
    &\var(Y_{k-1}\mid \cF_{k})\notag\\
    &\quad =\Var(Z_k-2\cdot \mathbf{1}[Z_k\ge 1]\mid \cF_{k})\le 2\Var(Z_k \mid \cF_k) + 2 \Var(2\cdot \mathbf{1}[Z_k\ge 1] \mid \cF_{k})\notag\\ 
   & \quad \le 2 B_k Y_k + 8 \IP(Z_k\ge 1\mid \cF_{k})\le 2 B_k Y_k + 8 \E(Z_k\mid \cF_k) 
    = 10 B_k Y_k.
\end{align}
Thus, 
\begin{align}\label{eq:ubvarW}
    \var(W_{k-1}\mid \cF_k) = \Phi^2_{k-1} \var(Y_{k-1}\mid\cF_k) \leq 10\Phi^2_{k-1} B_k Y_k.
\end{align}

Let $\bB$
be the $\sigma$-field generated by the beta variables
$(B_j)^{n-1}_{j=2}$, and 
let $\E_\bB$ and $\var_\bB$ denote conditional expectation and variance
with respect to $\bB$.
Note that $M_0,\dots,M_{n-1}$ is a reverse
martingale also conditioned on $\bB$,
since $\bB=\cF_{n-1}\subseteq \cF_k$ for every $k$. 
In particular,
 \begin{align}\label{eq:ubmeanW}
     \E_\bB W_k \leq \E_\bB M_k = M_{n-1} = 2\Phi_{n-1}\quad \text{{for $0\leq k\leq n-1$}}. 
 \end{align}
Hence, by applying \eqref{eq:mg}, the reverse   martingale property, 
\eqref{eq:ubvarW}, \eqref{de:phi}, \eqref{de:Wk}, and then
  \eqref{eq:ubmeanW}, 
for $0\le k\le n-1$,
\begin{align}
    \var_\bB (M_k) &= \E_\bB (M_k - 2\Phi_{n-1})^2 = \sum^{n-1}_{j=k+1}\E_\bB \var(W_{j-1}\mid \cF_j)
\notag\\
    &\leq 10 \sum^{n-1}_{j=k+1} \Phi^2_{j-1}B_j \E_\bB(Y_j) = 10 \sum^{n-1}_{j=k+1} \frac{B_j}{1+B_j}\Phi_{j-1}\E_\bB (W_j) 
\notag\\
    & \leq {20}  \sum^{n-1}_{j=k+1} \frac{B_j}{1+B_j}\Phi_{j-1}\Phi_{n-1} = {20} \sum^{n-1}_{j=k+1} \Phi^2_{j-1} B_j \prod^{n-1}_{i=j+1}(1+B_i).\numberthis\label{eq:varM}
\end{align}

\subsection{Some estimates}\label{sse:est}
Below we provide several estimates for $W_k$, $M_k$, $Z_k$ and $A_k$ that we
need later. The results are analogous to \cite[Lemmas 2.1--2.3]{Janson2023}.
Recall that $\chi=1/2$ for $m=2$ and $\rho=0$.  

\begin{lemma}\label{le:doob}
   For $1\leq k\leq n-1$, we have
    \begin{align}\label{eq:meanW2}
      \E_\bB W_k^2 \leq  \E_\bB  M_k^2\leq 20\sum^{n-1}_{j=k+1}\Phi_{j-1}^2  B_j  \prod^{n-1}_{i=j+1}(1+B_i) + 4\Phi_{n-1}^{2}.
    \end{align}
Furthermore,
    \begin{align}\label{mW2}
        \E_\bB \max_{0\le k\le n-1} W_k^2 \le  \E_\bB \max_{0\leq k \leq n-1} M_k^2  \le 4 \E_\bB M_0^2
    \end{align}
    and there is a positive constant $C$
    such that
    \begin{align}\label{eq:maximalineq}
      \E \max_{0\leq k \leq n-1} W_k^2 \le \E \max_{0\leq k \leq n-1} M_k^2 \leq Cn.
    \end{align}
\end{lemma}

\begin{proof}
The first inequalities in \eqref{eq:meanW2} and \eqref{mW2} follow from $0\le W_k\leq M_k$. 
  For the second inequality in \eqref{eq:meanW2}, note that 
    \begin{align}\label{eq:notdoob}
        \E_\bB M_k^2 = \var_\bB(M_k) +  (\E_\bB M_k)^2,
    \end{align}
and the inequality follows from this by the inequality \eqref{eq:varM} and the equalities in \eqref{eq:ubmeanW}.

The second inequality in \eqref{mW2} follows from Doob's inequality.
By \eqref{mW2}, \eqref{eq:maximalineq} follows from showing $\E M_0^2 \le Cn$.
For every $0\leq k\leq n-1$,
we use \eqref{eq:varM} and the independence of the beta variables to obtain
    \begin{align}
        \E\big(\var_\bB(M_k)\big) &\leq 20 \E\sum^{n-1}_{j=k+1} \Phi^2_{j-1} B_j \prod^{n-1}_{i=j+1}(1+B_i)
\notag\\
        &=20\sum^{n-1}_{j=k+1}\E\Phi_{j-1}^2 \E B_j  \prod^{n-1}_{i=j+1}\E(1+B_i).
    \end{align}
    So applying \eqref{eq:Bmean}, \eqref{eq:mprodbeta} and
      \eqref{eq:phibd}, we get
    \begin{align}\label{eq:varMbd}
         \E\big(\var_\bB(M_k)\big) 
         \leq  C \sum^{n-1}_{j=k+1} \frac{j}{j}\cdot\bbclr{\frac{n}{j}}^{1/2}
         \leq C n^{1/2}\sum^{n-1}_{j=k+1} \frac{1}{j^{1/2}} \leq C n .
    \end{align}
    Using the {equalities} in \eqref{eq:ubmeanW} and the estimate in \eqref{eq:phibd}, we also have
    \begin{align}\label{eq:EM2}
        \E\bigsqpar{( \E_\bB M_k)^2} 
= \E\bclr{M_{n-1}^2} = 4 \E(\Phi_{n-1}^2)\leq Cn \quad \text{{for $0\leq k\leq n-1$}}.
    \end{align}
    Applying \eqref{eq:varMbd} and \eqref{eq:EM2} to \eqref{eq:notdoob}, with $k=0$, we thus have 
    \begin{align}
        \E  M_0^2  \leq {\E}\bigsqpar{\var_\bB(M_0)} +
          {\E}\bigsqpar{(\E_\bB M_0)^2}\leq  C n,
    \end{align}
    which together with \eqref{mW2} implies \eqref{eq:maximalineq}.
\end{proof}


\begin{lemma}\label{le:Zk}
There is a positive constant $C$ such that,    
for $1\leq k\leq n-1$,
    \begin{align}
        \IP(Z_k\geq 1) &\leq C\frac{n^{1/2}}{k^{{3/2}}}; \label{eq:Zkg1}\\ 
        \IP(Z_k\geq 2) &\leq C\frac{n}{k^3} \label{eq:Zkg2}.
    \end{align}
\end{lemma}

\begin{proof}
    We start by proving \eqref{eq:Zkg1}. Firstly, it follows from \eqref{eq:distZk} and \eqref{de:Wk} that 
    \begin{equation}
        \E(Z_k\mid \cF_k)=Y_k B_k=\Phi_k^{-1} B_k W_k, 
    \end{equation}
    which, along with \eqref{eq:ubmeanW}, imply that
    \begin{equation}\label{MarZ}
        \E_\bB(Z_k) = \Phi^{-1}_{k} B_k \E_\bB(W_k)\leq {2}\Phi^{-1}_{k} B_k \Phi_{n-1} = {2}B_k\prod^{n-1}_{i=k+1} (1+B_i).
    \end{equation}
    Using the independence of $(B_k)^{n-1}_{k=2}$, \eqref{eq:Bmean} and \eqref{eq:mprodbeta}, we therefore have
    \begin{align}\label{eq:meanZk}
        \E(Z_k)&{\leq 2}\E(B_k)\prod^{n-1}_{i=k+1}\E(1+B_i)
\leq\frac{C}{k}\bbclr{\frac{n}{k}}^{1/2}= C\frac{n^{1/2}}{k^{{3/2}}} ,
    \end{align}
    and \eqref{eq:Zkg1} follows  from Markov's inequality.

    The proof for \eqref{eq:Zkg2} is similar, this time we observe that by
    Markov's inequality, 
    \begin{align}\label{eq:Markov2}
        \IP(Z_k\geq 2\mid \cF_k) \leq \E\bbbclr{\binom{Z_k}{2}\biggm\vert \cF_k} = \binom{Y_k}{2} B_k^2 \leq B_k^2 
 Y^2_k = B_k^2\Phi_k^{-2} W_k^2.
    \end{align}
    By \eqref{eq:meanW2} of Lemma \ref{le:doob} and \eqref{de:phi}, we have 
    \begin{align*}
         &\E\big(B_k^2\Phi_k^{-2} W_k^2\big)
          = \E\big(B_k^2\Phi_k^{-2} \E_\bB(W_k^2)\big)\\
         & \leq 20 \E\bbbclr{B_k^2\Phi_k^{-2}\sum^{n-1}_{j=k+1}\Phi_{j-1}^2  B_j  \prod^{n-1}_{i=j+1}(1+B_i)} + 4\E\bclr{B^2_k\Phi_k^{-2}\Phi_{n-1}^2}\\
         &= 20 \E\bbbclr{B_k^2\sum^{n-1}_{j=k+1} B_j\prod^{j-1}_{i=k+1}(1+B_i)^2\prod^{n-1}_{l=j+1}(1+B_l)} + 4\E\bbclr{B^2_k\prod^{n-1}_{i=k+1}(1+B_i)^2}\\
         &=20 \E(B_k^2)\sum^{n-1}_{j=k+1} \E(B_j)\prod^{j-1}_{i=k+1}\E(1+B_i)^2\prod^{n-1}_{l=j+1}\E(1+B_l) + 4\E(B_k^2) \prod^{n-1}_{i=k+1}\E(1+B_i)^2
         \numberthis \label{eq:EZk2}
    \end{align*}
    Applying the estimates in \eqref{eq:Ebetasq}, \eqref{eq:mprodbeta} and \eqref{eq:trunphi2bd} to \eqref{eq:EZk2}, we find
\begin{align}\label{eq:BphiW2}
    \E\big(B_k^2\Phi_k^{-2} W_k^2\big)\leq \frac{C}{k^2}\sum^{n-1}_{j=k+1}\frac{1}{j} \cdot\frac{j}{k}\cdot\bbclr{\frac{n}{j}}^{1/2} + \frac{Cn}{k^3} \leq \frac{Cn}{k^3}.
\end{align}
Taking the expectation in \eqref{eq:Markov2} and plugging in 
 \eqref{eq:BphiW2} yields \eqref{eq:Zkg2}.
\end{proof}

\begin{lemma}\label{le:Ak}
    For $1\leq k\leq n-1$, 
\begin{equation}\label{eq:incAsim}
    A_{k-1} -A_k \leq (W_kB_k)^2\Phi^{-1}_{k}
\end{equation}
    and there is a positive constant $C$ such that
    \begin{equation}\label{eq:EAk}
        \E A_k \leq \frac{Cn}{k^{3/2}}.
    \end{equation}
\end{lemma}

\begin{proof}
    By \eqref{eq:incA}, Taylor's formula, 
the increasing property of $\Phi_k$, and \eqref{de:Wk},
    \begin{align}\label{eq:Abd}
        A_{k-1}-A_k&=2\Phi_{k-1}\bclr{(1-B_k)^{Y_k}-1+B_kY_k}\leq 2\Phi_{k-1} \binom{Y_k}{2} B_k^2 
\notag\\
        &\leq Y_k^2\Phi_{k-1} B_k^2 
\leq Y_k^2\Phi_{k} B_k^2 
=W_k^2 \Phi^{-1}_{k}B_k^2 ,
    \end{align}
    yielding \eqref{eq:incAsim}. To prove \eqref{eq:EAk}, we note that by a telescoping argument, \eqref{eq:incAsim} implies
    \begin{align}\label{eq:tele}
        A_k \leq \sum^{n-1}_{i=k+1} \bclr{W_i B_i}^{2} \Phi^{-1}_{i}
    \end{align}
    Furthermore, \eqref{eq:meanW2} of Lemma \ref{le:doob} and \eqref{de:phi} together yield
    \begin{align*}
        &\E_\bB\bclr{W^2_i B_i^{2} \Phi^{-1}_{i}}\\
        &\leq  20 B_i^2 \Phi^{-1}_{i} \sum^{n-1}_{j=i+1} \Phi_{j-1}^2 B_j\prod^{n-1}_{h=j+1}(1+B_h) + 4 B_i^2 \Phi^{-1}_{i} \Phi^2_{n-1}\\
        &=20 B_i^2 \Phi_i\sum^{n-1}_{j=i+1} B_j \prod^{j-1}_{l=i+1} (1+B_l)^2  \prod^{n-1}_{h=j+1}(1+B_h) + 4 B_i^2 \Phi_{i} \prod^{n-1}_{h=i+1}(1+B_h)^2\\
        &\leq 40 B_i^2 \Phi_{i-1}\sum^{n-1}_{j=i+1} B_j \prod^{j-1}_{l=i+1} (1+B_l)^2 \prod^{n-1}_{h=j+1}(1+B_h) + 8 B_i^2 \Phi_{i-1} \prod^{n-1}_{h=i+1}(1+B_h)^2
.\numberthis    
\end{align*}
Taking expectation and again using the independence of $(B_k)^{n-1}_{k=2}$, we have 
\begin{align*}
    \E\bclr{W^2_i B_i^{2} \Phi^{-1}_{i}} &\leq 40 \E(B_i^2)\E(\Phi_{i-1}) \sum^{n-1}_{j=i+1} \E(B_j) \prod^{j-1}_{l=i+1} \E(1+B_l)^2 \prod^{n-1}_{h=j+1}\E(1+B_h)\\
    &\qquad + 8 \E(B_i^2) \E(\Phi_{i-1}) \prod^{n-1}_{h=i+1}\E(1+B_h)^2.
\numberthis
\end{align*}
Applying \eqref{eq:Bmean}, \eqref{eq:Ebetasq}, \eqref{eq:mprodbeta}, \eqref{eq:meanphik} and \eqref{eq:trunphi2bd} to the last display gives
\begin{align}\label{eq:Asummand}
    \E\bclr{W^2_i B_i^{2} \Phi^{-1}_{i}}\leq \frac{C}{i^{3/2}} \sum^{n-1}_{j=i+1} \frac{1}{2j-2} \cdot \frac{j}{i}\cdot \bbclr{\frac{n}{j}}^{1/2} + \frac{Cn}{i^{5/2}}\leq \frac{Cn}{i^{5/2}}.
\end{align}
Thus, in view of \eqref{eq:Asummand} and \eqref{eq:tele}, we deduce that for $1\leq k\leq n-1$,
\begin{align}
    \E (A_k) \leq \sum^{n-1}_{i=k+1} \frac{Cn}{i^{5/2}}\leq \frac{Cn}{k^{3/2}},
\end{align}
hence proving \eqref{eq:EAk}.
\end{proof}

\section{The early part and a Yule process}\label{se:Yule}
We continue to study the case $m=2$ and $\rho=0$, and recall that then
$\chi=1/2$. 
We show that the early part
of the growth of $D_n$ can be closely coupled to the same
time-changed Yule process as in \cite{Janson2023}, and use this coupling to study $Y_k$ and $W_k$. We start by presenting its construction and key features, following the description in \cite{Janson2023}.

Let $\mathcal{Y}$ be a Yule process starting with two particles, and let
$\mathcal{Y}_t$ be  the number of (living) particles at time $t$ (thus
$\mathcal{Y}_0=2$). Note that $\mathcal{Y}_t$ has the same distribution as
the sum of two copies of the standard Yule process, which starts with a
single particle. 
(See e.g.\ \cite[Section III.5]{Athreya-Ney1972} for definition and some
basic properties.)
To better compare the process to $D_n$, it is convenient
too to view the Yule process $\mathcal{Y}$ as a tree, where the root
$\gamma_0:=0$ marks the beginning of the process and the vertex $\gamma_i$
marks the time of the $i$-th particle split in the process. Note also these
split times are a.s.\ distinct. In this way, each particle can be
represented by an edge from its time of birth to its time of death. 

The time-changed Yule tree $ \mathcal{\widehat Y}$ appearing in
\cite{Janson2023} is obtained by applying the mapping $t\to e^{-t}$, so that
the vertices in $ \mathcal{\widehat Y}$ have labels $e^{-\gamma_i}\in(0,1]$. 
Hence, the root in $\hcY$ has label 1, and
a particle in the original Yule process that is born at time $\gamma_i$ and has lifetime $\tau\in \mathrm{Exp}(1)$ is now represented by an edge from $x=e^{-\gamma_i}$ to $e^{-{(\gamma_i+\tau)}}=xe^{-\tau}=xU$, where $U:=e^{-\tau}\in\sU(0,1)$. In light of this, as well as that $e^{-\gamma_0}=1$ and all lifetimes in the original Yule process are independent and have the $\mathrm{Exp}(1)$ distribution, any vertex in $\mathcal{\widehat Y}$ that is $d$ generations away from the root therefore take labels of the form $\widehat U_1\cdots \widehat U_d$, where $\widehat U_1,\dots,\widehat U_d\in\sU(0,1)$ and are independent. 

Let $\xD_n$ be
the random red graph $D_n$ with label $k$ replaced with $(k/n)^\chi$,
so that the labels now take values in $(0,1]$.
We regard $D_n$ as rooted at $n$; thus the root of $\xD_n$ has label $1$.
We shall compare the time-changed Yule tree $\mathcal{\widehat Y}$ to
$ \xD_n$, considering only vertices with large enough labels. 
In preparation, let 
\begin{equation}\label{eq:n1}
    n_1=n_1^{(n)}:=\floor{n/\log n}.
\end{equation}
We will use the construction of $G_n$ in \eqref{de:S}--\eqref{eq:sU},
using the variables $S_{n,k}$ defined there;
in particular, recall that $\tU_{k,\ell}$ are independent $\sU(0,1)$ random
variables. 

\begin{lemma}\label{LU}
Let\/ $\kappa_n:=\log n/n^{1/3}$.
  With probability at least $1-C\log n/n^{1/3}$, the following hold:
  \begin{romenumerate}
    
  \item \label{LUa}
For every path in $D_n$ between vertex $n$ (the root)
and a vertex $k>n_1$  consisting of $d+1\ge2$ 
red vertices $n=v_0>v_1>\dots>v_{d-1}>v_d=k$ such that
$v_i\overset{\ell_i}\to v_{i+1}$ for $0\le i< d$, we have
\begin{align}\label{lu1}
\bigabs{\widetilde U_{v_0, \ell_0}\cdots \widetilde U_{v_{d-1}, \ell_{d-1}} -
\bigpar{\tfrac{k}{n}}^\chi }
\le 3d\kk_n.
\end{align}
\item   \label{LUb}
For every such path in $D_n$ between vertex $n$
and a vertex $k\le n_1$, we have
\begin{align}\label{lu2}
\widetilde U_{v_0, \ell_0}\cdots \widetilde U_{v_{d-1}, \ell_{d-1}} 
\le \bigpar{\tfrac{n_1}{n}}^\chi + 3d\kk_n.
\end{align}
\end{romenumerate}
\end{lemma}

\begin{proof}
  We may assume that $n$ is large enough such that $n_1\qw < \kk_n$, since
  the result is trivial for small $n$ by choosing $C$ large enough.

Firstly it follows from
\eqref{kil} 
and \eqref{eq:sU} that if $k\overset\ell\to i$, then
\begin{align}\label{eq:Sinterval}
    S_{n,i-1}\leq \widetilde U_{k,\ell} S_{n,k-1}<S_{n,i}.
\end{align}
Secondly,
again assuming that $n$ is large,
we take $\eps=1/3$ and $\psi_n=n_1$ in Lemma
\ref{le:Sest} and find 
that 
there is a positive constant $C$ such that
with probability at least $1-C\log n/n^{1/3}$, 
\begin{align}\label{eq:Sest}
    \max_{n_1\leq j\leq n} \big| S_{n,j}-\bclr{\tfrac{j}{n}}^\chi \big| 
\leq 2 n_1^{-1/3}
\le 3 \log^{1/3} n/n^{1/3}
\leq \kappa_n.
\end{align}
We assume in the rest of the proof that \eqref{eq:Sest} holds, and show first that
\eqref{lu1} follows by induction on $d$.
Note first that if $j> n_1$, then \eqref{eq:Sest} implies
\begin{align}\label{eq:Sest-}
  S_{n,j-1} \ge \bclr{\tfrac{j-1}{n}}^\chi -\kk_n&
= \bclr{\tfrac{j}{n}}^\chi  \bclr{1-\tfrac{1}{j}}^\chi -\kk_n
\ge \bclr{\tfrac{j}{n}}^\chi  \bclr{1-\tfrac{1}{j}} -\kk_n
\ge \bclr{\tfrac{j}{n}}^\chi -\tfrac{1}{j} -\kk_n
\notag\\&
\ge \bclr{\tfrac{j}{n}}^\chi -2\kk_n.
\end{align}

For the base case $d=1$ we have by \eqref{eq:Sinterval}, \eqref{eq:Sest},
\eqref{eq:Sest-}, 
and recalling $S_{n,n-1}=1$,
\begin{align}
\tU_{n,\ell_0} \le S_{n,k} \le \bclr{\tfrac{k}{n}}^\chi +\kk_n
\label{eq:basel} 
\end{align}
and 
\begin{align}
\tU_{n,\ell_0} \ge S_{n,k-1} \ge \bclr{\tfrac{k}{n}}^\chi -2\kk_n,
\label{eq:base2} 
\end{align}
which show \eqref{lu1} in this case.

For $d\ge2$, we use induction and find, using the induction hypothesis and
\eqref{eq:Sinterval}--\eqref{eq:Sest-},
\begin{align}\label{lu3}
\widetilde U_{v_0, \ell_0}\cdots \widetilde U_{v_{d-1}, \ell_{d-1}} &
\le \Bigpar{\bigpar{\tfrac{v_{d-1}}{n}}^\chi +3(d-1)\kk_n}\tU_{v_{d-1}, \ell_{d-1}} 
\notag\\&
\le \bigpar{S_{n,v_{d-1}-1}+2\kk_n+(3d-3)\kk_n}\tU_{v_{d-1}, \ell_{d-1}} 
\notag\\&
\le S_{n,v_{d-1}-1}\tU_{v_{d-1}, \ell_{d-1}} + (3d-1)\kk_n
\notag\\&
< S_{n,k} + (3d-1)\kk_n
\notag\\&
\le \bigpar{\tfrac{k}{n}}^\chi +3d\kk_n
\end{align}
and similarly,
using also $S_{n,v_{d-1}}\ge S_{n,v_{d-1}-1}$,
\begin{align}\label{lu4}
\widetilde U_{v_0, \ell_0}\cdots \widetilde U_{v_{d-1}, \ell_{d-1}} &
\ge \Bigpar{\bigpar{\tfrac{v_{d-1}}{n}}^\chi -3(d-1)\kk_n}\tU_{v_{d-1}, \ell_{d-1}} 
\notag\\&
\ge \bigpar{S_{n,v_{d-1}}-\kk_n-(3d-3)\kk_n}\tU_{v_{d-1}, \ell_{d-1}} 
\notag\\&
\ge S_{n,v_{d-1}-1}\tU_{v_{d-1}, \ell_{d-1}} - (3d-2)\kk_n
\notag\\&
\ge S_{n,k-1} - (3d-2)\kk_n
\notag\\&
\ge \bigpar{\tfrac{k}{n}}^\chi -3d\kk_n
.\end{align}
These inequalities prove \eqref{lu1}, which completes the proof of \ref{LUa}.

To prove \ref{LUb},
assume first $v_{d-1}> n_1\ge k$.
Then, using \eqref{lu1}, the first lines of \eqref{lu3} still hold and yield 
\begin{align}\label{lu5}
\widetilde U_{v_0, \ell_0}\cdots \widetilde U_{v_{d-1}, \ell_{d-1}} &
< S_{n,k} + (3d-1)\kk_n
.\end{align}
Furthermore, by $n_1\ge k$ and \eqref{eq:Sest},
\begin{align}\label{lu6}
  S_{n,k} \le S_{n,n_1}
\le \bigpar{\tfrac{n_1}{n}}^\chi +\kk_n,
\end{align}
and \eqref{lu2} follows by \eqref{lu5} and \eqref{lu6}.
Finally, in the remaining case $v_{d-1}\le n_1$, we use the trivial
$\tU_{v_0, \ell_0}\cdots\tU_{v_{d-1}, \ell_{d-1}} 
\le \tU_{v_0, \ell_0}\cdots\tU_{v_{d-2}, \ell_{d-2}} $ and induction on $d$.
\end{proof}

Recall that a vertex in $\mathcal{\widehat Y}$ that is $d$ generations
   away from the root has label of the form $\widehat U_1\cdots \widehat
   U_d$, where $\widehat U_i\in \sU[0,1]$ and are independent.
   In view of \eqref{lu1}, we  couple $\mathcal{\widehat Y}$ and
 $\xD_n$  by generating them together as follows, 
where we also  construct
a  mapping $\Psi$ of the vertices of $\hcY$ to the vertices of $\xD_n$. 
In the construction below, 
$\hcY$ and $\xD_n$ will be finite subsets of the
final Yule tree and digraph, and $\Psi$ maps the current $\hcY$ onto the
current $\xD_n$.
Recall that $\hcY$ and $\xD_n$ determine $\cY$ and $D_n$ by 
(deterministic) relabelling.
\begin{enumerate}
    \item Sample the beta variables $(B_j)^{n-1}_{j=2}$ defined in
\eqref{de:betas}.
This defines also $S_{n,j}$ by  \eqref{de:S}.
\item 
Start the construction by letting
$\hcY$ and $\xD_n$ just consist of their roots, both labelled  1.
Let $\Psi$ map the root of $\hcY$ to the root of $\xD_n$.
\item 
Let $x$ be the vertex in the constructed part of $\hcY$ that has the largest label
among all vertices that do not yet have children assigned.
Give $x$ children $xU'_{x,1}$ and $xU'_{x,2}$ 
(which are added to the current $\hcY$),
where $U'_{x,\ell}$ are independent $\sU[0,1]$ variables that are independent
of all other variables.

The vertex $x$ is mapped to  some vertex $\Psi(x)=(k/n)^\chi$ in
$\xD_n$, which thus corresponds to vertex $k$ in $D_n$.
There are three cases:
\begin{enumerate}
\item 
If $k>1$ and $\Psi(x)$ has not yet got any children,
define $\tU_{k,\ell}:=U'_{x,\ell}$ for $\ell=1,2$.
This defines by \eqref{kil}--\eqref{eq:sU} 
the edges from $k$ and thus the children of $k$
in $D_n$; if these children are $k_1$ and $k_2$, 
the corresponding children in $\xD_n$ are $(k_1/n)^\chi$ and $(k_2/n)^\chi$;
we add them to $\xD_n$
and we define $\Psi(x_{\ell}):=(k_\ell/n)^\chi$,
thus mapping the children of $x$ in $\hcY$
to the children of $\Psi(x)$ in $\xD_n$. 

\item 
If $k>1$ and $\Psi(x)$ already has children
(because it equals $\Psi(y)$ for some $y>x$),
then we just extend $\Psi$ 
by mapping the children of $x$ to the children of $\Psi(x)$
(in any order). 

\item 
If $k=1$, so $x$ maps to $v=(1/n)^\chi$ (which has no children in the final
$\xD_n$), 
we extend $\Psi$ by mapping also the children of $x$ to $v$.
\end{enumerate}

\item Repeat from (3) (\emph{ad infinitum}).

\end{enumerate}

It is easy to see that 
running this ``algorithm'' an infinite number of iterations yields 
$\hcY$ and $\xD_n$ with the right distributions, together with 
a map $\Psi$ of the vertices of $\hcY$
onto the vertices of $\xD_n$ such that every path from the root in $\xD_n$ is
the image of a path from the root in $\hcY$. 
$\Psi$ is obviously not injective since $\hcY$ is an infinite
tree. Nevertheless, we show that restricted to rather large labels,
the mapping $\Psi$
is \whp\ a bijection which perturbs the vertex labels with small errors.

\begin{theorem}\label{th:coupling}
    Let $n_1:=\floor{n/\log n}$. We may w.h.p.\ couple the
    $\xD_n$ and the time-changed Yule tree $\mathcal{\widehat Y}$,
    such that considering only vertices with labels in
    $((n_1/n)^\chi,1]$ and edges with the starting points in this set,
    there is a bijection between these sets of vertices in the two models
    which displaces each label by at most $3\log^2n/n^{1/3}$, and
    a corresponding bijection between the edges (preserving the incidence
    relations). 
In particular, \whp{} 
\begin{equation}\label{eq:Yulecoupling}
Y^{(n)}_{n_1} = \mathcal{\widehat Y}_{(n_1/n)^\chi},
\end{equation}
where $\hcY_x=\cY_{-\log x}$ is the number  of edges in $\hcY$ alive at time $x$.
\end{theorem}

\begin{proof}
The proof is similar to that of \cite[Theorem 3.1]{Janson2023}, but with
several technical complications.  
We use the coupling constructed before the theorem.
Let $\gd_n:= 3\log^2n/n^{1/3}=(\log n)\kk_n$, with $\kk_n$ as in \refL{LU}.

\stepp\label{TCO1}
We first note that if some vertex in $\xD_n$ with label in 
$[(n_1/n)^\chi,1]$ is the image of two or more vertices in $\hcY$, then the
corresponding vertex $k\ge n_1$ in $D_n$ can be reached from $n$ by at least
two different paths in $D_n$, and if we let $k$ be maximal with this
property, then its indegree $Z_k\ge2$.
Consequently, the probability that this happens is at most,
using \eqref{eq:Zkg2} of Lemma \ref{le:Zk},
\begin{align}\label{medges}
    \sum^{n-1}_{k=n_1}\IP(Z_k\geq 2) \leq C n
  \sum^{n-1}_{k=n_1}k^{-3}=O(n/n^2_1) = {O(\log^2 n/n)}
=o(1).
\end{align}
Hence, \whp\ the mapping $\Psi$ from $\hcY$ to $\xD_n$ is injective at every 
vertex in $\xD_n\cap[(n_1/n)^\chi,1]$.
We may thus in the sequel assume that this injectivity holds.
Note that this implies that in the construction of the mapping $\Psi$ above, 
we have $\tU_{k,\ell}=U'_{x,\ell}$ for every $k\ge n_1$ 
and vertex $x\in\hcY$ such that $\Psi(x)=(k/n)^\chi$.

\stepp\label{TCO2} 
As in \cite{Janson2023},
$\mathcal{\widehat Y}_x=\mathcal{Y}_{-\log x}$ for every $x\in (0,1]$, so by
standard properties of the Yule process 
(see e.g.\ \cite[Section 3]{Janson2023})
    \begin{align}\label{eq:Yule1}
        \E \mathcal{\widehat Y}_x = \E \mathcal{Y}_{-\log x} = 2e^{-\log x} = 2/x.
    \end{align}
In $\mathcal{\widehat Y}$, there are $\mathcal{\widehat Y}_x-1$ vertices
    with labels in $[x,1]$. 
Thus, \eqref{eq:Yule1} implies that w.h.p.\ there are less than 
$\floor{\log n}$  
such vertices for $x=(n_1/n)^\chi \sim \log^{-\chi} n$.
It follows trivially that w.h.p., in $\hcY$ the number of generations
    from the root to any point in $[(n_1/n)^\chi,1]$ is less than 
$\floor{\log n}$. 
Hence, we may assume this property too.

\stepp\label{TCO3}
 The expected number of vertices in $\mathcal{\widehat Y}$ that are
within $\gd_n$ from $(n_1/n)^\chi$ is, by \eqref{eq:Yule1},
\begin{align*}
 \E\bclr{\hcY_{(n_1/n)^\chi-\gd_n} - \hcY_{(n_1/n)^\chi+\gd_n}}
&= \frac{2}{(n_1/n)^\chi-\gd_n} - \frac{2}{(n_1/n)^\chi+\gd_n} 
\\&
= O\bbclr{\frac{\gd_n}{(n_1/n)^{2\chi}}}
=O\bclr{\gd_n \log^{2\chi}n}
=o(1).
\numberthis
\label{eq:displacement}
    \end{align*}
Hence, \whp\ there are no vertices $x$ in $\hcY$ 
with $|x-(n_1/n)^\chi|\le\gd_n$. 
We may  in the sequel assume this.

\stepp\label{TCO4}
Consequently,  \whp{} the properties in 
\refSteps{TCO1}--\ref{TCO3} hold, and also 
the conclusions \ref{LUa} and \ref{LUb} of
\refL{LU}. We assume this for the rest of the proof.

Suppose that $k> n_1$ is a vertex of $D_n$, and let $v=(k/n)^\chi$ be the
corresponding vertex of $\xD_n$. 
By \refStep{TCO1}, $v=\Psi(x)$ for a unique vertex $x\in\hcY$.
If $x\in((n_1/n)^\chi,1]$, then the number of generations from the root
to $x$ in $\hcY$ is at most $\log n$ by \refStep{TCO2}. The number of
generations to $\Psi(x)$ in $\xD_n$ is the same, and it follows from
\eqref{lu1} and the equality
$\tU_{k,\ell}=U'_{x,\ell}$  in \refStep{TCO1} that,
denoting the path from $n$ to $k$ as in \refL{LU},
\begin{align}\label{tco4}
|x-\Psi(x)|=
  \bigabs{\tU_{v_0, \ell_0}\cdots \tU_{v_{d-1}, \ell_{d-1}} 
-\bigpar{\tfrac{k}{n}}^\chi}
\le 3 (\log n) \kk_n = 3\log^2n/n^{1/3} = \gd_n.
\end{align}
It remains to show only that no vertex $x$ in $\hcY$ is pushed over the
boundary $(n_1/n)^\chi$ (in any direction) by $\Psi$.

\stepp\label{TCO5}
Suppose that there exists a vertex $x\le (n_1/n)^\chi$ in $\hcY$ such that
$\Psi(x)> (n_1/n)^\chi$. Let $y>x$ be the parent of $x$ in $\hcY$, so that
  $\Psi(y)>\Psi(x)$. Assume also $y>(n_1/n)^\chi$. By \refStep{TCO2}, it
follows that the number of generations from the root
to $y$ is less than $\floor{\log n}$, and thus the number of generations to
$x$ is at most $\floor{\log n}$.
Consequently, \eqref{lu1} shows, similarly to \eqref{tco4}, that
\begin{align}
|x-\Psi(x)|
\le 3 (\log n) \kk_n =\gd_n
\end{align}
and hence
\begin{align}
  x\ge \Psi(x)-\gd_n \ge (n_1/n)^\chi-\gd_n.
\end{align}
However, we have also $x\le(n_1/n)^\chi$, so $|x-(n_1/n)^\chi|\le\gd_n$,
and by \refStep{TCO3}, there is no such vertex $x$ in $\hcY$.
If $y\le(n_1/n)^\chi$, we may instead replace $x$ by $y$ (and repeat this if
necessary) until we encounter a vertex $x$ with parent $y$ such that
$x\le(n_1/n)^\chi$, $\Psi(x)>(n_1/n)^\chi$ and $y>(n_1/n)^\chi$. However, we
have shown that such a pair cannot exist.  

\stepp\label{TCO6}
Suppose that there exists a vertex $x>(n_1/n)^\chi$ in $\hcY$ such that
$\Psi(x)\le (n_1/n)^\chi$. 
By \refStep{TCO2}, it follows that the number of generations from the root
to $x$ is less than $\log n$.
This time \eqref{lu2} applies and shows that
\begin{align}
x\le 
\bigpar{\tfrac{n_1}{n}}^\chi +  3 (\log n) \kk_n =
\bigpar{\tfrac{n_1}{n}}^\chi +  \gd_n.
\end{align}
However, by \refStep{TCO3} again, there is no such vertex $x$ in $\hcY$.

The various claims proved in the steps above show that with the coupling
and mapping $\Psi$ constructed before the theorem, \whp\ $\Psi$ yields a bijection
with the stated properties.
In particular, \whp\ $Y_{n_1}$ equals the number of edges in $\hcY$ that
are alive at $(n_1/n)^\chi$,
meaning that they start in $J:=[(n_1/n)^\chi,1]$ and end outside $J$.
(Note that a.s.\ $\hcY$ has no point exactly at $(n_1/n)^\chi$, so it does
not matter whether we include that point in $J$ or not.)
Finally, note that for any $x\in(0,1)$,
the number of edges of $\hcY$ that are alive at $x$
equals the number of edges in $\cY$ that are alive at $-\log x$,
which equals the number $\cY_{-\log x}$ of particles at $-\log x$, since
each edge represents the lifeline of one particle.
\end{proof}

We define
\begin{equation}\label{de:Xi}
    \Xi=\Xi^{(n)}:=\frac{W^{(n)}_{n_1}}{n^{\chi}}.
\end{equation}
This random variable plays the same role as $\Xi$ in \cite{Janson2023}, but
note the different scaling.
Recall also  $\tgb$ defined in \eqref{lb2a}--\eqref{lb2b}.

\begin{lemma}\label{LXi}
    As $n\to\infty$, 
    \begin{equation}\label{lxi}
\Xi^{(n)}\dto \tgb\cdot \xi,
    \end{equation}
where $\tgb$ is given by \eqref{lb2b}, and $\xi\in\mathrm{Gamma}(2,1)$ is
independent of $\tgb$.
\end{lemma}

\begin{proof}
As in \cite[Lemma 3.2]{Janson2023}, 
we first generate the Yule process $\cY$, and then for each $n$ separately,
we construct $D_n$ by the construction given before \refT{th:coupling}.
This yields for each $n$ the coupling in \refT{th:coupling}. 
In particular, \eqref{eq:Yulecoupling} holds w.h.p.
By standard properties of the Yule process 
(see e.g.\ \cite[Section III.5  and Problem III.2]{Athreya-Ney1972}),
\begin{equation}\label{eq:Yulepop}
  x\mathcal{\widehat Y}_{x} = x\mathcal{Y}_{-\log x} \asto \xi 
\quad \text{as $x\to 0$},
\end{equation}
with $\xi\in\mathrm{Gamma}(2,1)$.
Therefore, \eqref{eq:Yulecoupling} and \eqref{eq:Yulepop} together imply that
    \begin{equation}\label{eq:Yconvergence}
        \bbclr{\frac{n_1}{n}}^{\chi} Y^{(n)}_{n_1} \pto \xi.
    \end{equation}
Hence, using also \eqref{de:Xi}, \eqref{de:Wk}, and \eqref{lb2b},
\begin{align}\label{lxi3}
  \Xi\nn 
= n^{-\chi}\Phi_{n_1}Y_{n_1}
= n_1^{-\chi}\Phi_{n_1}\bbclr{\frac{n_1}{n}}^{\chi} Y_{n_1}
\pto \tgb\xi.
\end{align}

Finally,
$Y\nn_{n_1}$ is a function of $(B_i)_{i>n_1}$
and the coin tosses made for $k>n_1$ in the construction.
Hence, for any fixed $K$, $Y\nn_{n_1}$ is independent of $(B_i)_{i=1}^K$ for
large enough $n$. Consequently,
\eqref{eq:Yconvergence} implies that $\xi$ is independent of $(B_i)_{i=1}^K$ for
every $K<\infty$, and thus $\xi$ and $\gb$ are independent.
\end{proof}

\section{The flat middle part}\label{se:flat} 
Let $n_2$ be any sequence of integers satisfying $n^{1/3}\ll n_2 \leq
n_1$. We show that similar to the case of uniform attachment, the
variable $W_k$ does not fluctuate much in the range $n_1\geq k\geq
n_2$. Below we give analogues of \cite[Lemmas 4.1--4.2 and Theorem
4.3]{Janson2023}, where we recall the definitions of $W_k$, $M_k$, $A_k$ and
$\Xi^{(n)}$ in \eqref{de:Wk}, \eqref{eq:doobdecom} and \eqref{de:Xi}. The
results and proofs are again similar, but with a different  scaling.  
\begin{lemma}\label{le:maxAk}
    As $n\to\infty$, 
    \begin{equation}\label{eq:IImaxAk}
        \max_{n_2\leq k\leq n-1} \bigg|\frac{A_k}{n^{1/2}}\bigg| = \frac{A_{n_2}}{n^{1/2}}\pto 0.
    \end{equation}
\end{lemma}

\begin{proof}
    The first equality in \eqref{eq:IImaxAk} follows from the fact that $A_k$ is reverse increasing. By \eqref{eq:EAk} in Lemma \ref{le:Ak} and the choice of $n_2$, 
    \begin{equation}
        \E \frac{A_{n_2}}{n^{1/2}} \leq C\frac{n}{n^{1/2}n_2^{3/2}} =o(1),
    \end{equation}
    implying the convergence in probability in \eqref{eq:IImaxAk}.
\end{proof}

\begin{lemma}\label{le:maxMk}
     As $n\to\infty$, 
     \begin{equation}\label{lfM}
     \max_{0\leq k \leq n_1} \bigg|\frac{M_k}{n^{{1/2}}}-\Xi^{(n)}\bigg| \pto 0.
     \end{equation}
\end{lemma}

\begin{proof}
Recall that $M_k$ is a reverse martingale.
Hence we obtain by Doob's inequality, \eqref{eq:mg}, \eqref{eq:ubvarW},
\eqref{de:Wk}, \eqref{de:phi}, and \eqref{eq:ubmeanW},
    \begin{align*}
        \E  \max_{0\leq k \leq n_1} |M_k-M_{n_1}|^2 &\leq 4\E|M_0-M_{n_1}|^2 = 4 \sum^{n_1}_{i=1} \E \var(W_{i-1}\mid \cF_i)\\
        &\leq 40 \sum^{n_1}_{i=1} \E\bclr{\Phi_{i-1}^{2}B_iY_i}=40 \sum^{n_1}_{i=1} \E\bclr{\Phi_{i-1}W_iB_i(1+B_i)^{-1}} \\
     &= 40 \sum^{n_1}_{i=1}\E\bclr{\Phi_{i-1}B_i(1+B_i)^{-1}\E_\bB(W_i)}\\
        &\le 80 \sum^{n_1}_{i=1} \E\bclr{\Phi_{i-1}B_i(1+B_i)^{-1}\Phi_{n-1}} \\
   & = 80  \sum^{n_1}_{i=1} \E\bbclr{\Phi^2_{i-1}B_i\prod^{n-1}_{j=i+1}(1+B_j)}.
\numberthis\label{hw2}
    \end{align*}
    Using the independence of $(B_i)^{n-1}_{i=2}$, 
    \eqref{eq:Bmean}, \eqref{eq:mprodbeta} and \eqref{eq:phibd}, it
    follows from \eqref{hw2} that 
    \begin{align*}
         \E \big( \max_{0\leq k \leq n_1} |M_k-M_{n_1}|\big)^2
&=  \E  \max_{0\leq k \leq n_1} |M_k-M_{n_1}|^2 
\leq C \sum^{n_1}_{i=1} \frac{i}{i}\bbclr{\frac{n}{{i}}}^{1/2}\\
          &\leq C(n n_1)^{1/2} = o(n).\numberthis \label{mmm}
    \end{align*}
 Thus, by \eqref{de:Xi}, \eqref{eq:doobdecom},
 the triangle inequality, Lemma \ref{le:maxAk}, and \eqref{mmm},
    \begin{align}
        \max_{0\leq k \leq n_1} \bigg|\frac{M_k}{n^{1/2}}-\Xi^{(n)}\bigg| &\leq  \max_{0\leq k \leq n_1} \bigg|\frac{M_k}{n^{1/2}}-\frac{M_{n_1}}{n^{1/2}}\bigg| + \frac{A_{n_1}}{n^{1/2}}  \pto 0,
    \end{align}
    as required.
\end{proof}

\begin{theorem}\label{TF1}
    As $n\to\infty$, 
    \begin{equation}\label{tf1}
        \max_{n_2\leq k\leq n_1} \bigg|\frac{W_k}{n^{1/2}}-\Xi^{(n)}\bigg|\pto 0.
    \end{equation}
\end{theorem}

\begin{proof}
    The result is a direct consequence of \eqref{eq:doobdecom}, the
    triangle inequality, and Lemmas \ref{le:maxAk} and~\ref{le:maxMk}.
\end{proof}

\section{The final part: tightness}\label{Stig}
Most vertices in $D_n$ turn out to have labels of the order $n^{1/3}$. To study this region in detail, we extend the processes $W_k$, $M_k$ and $A_k$ to real arguments $t\in[0,n-1]$
by linear interpolation. We  for convenience extend them further to
$t\in[0,\infty)$ by defining them to be constant on $[n-1,\infty)$.

The aim of this and the next section is to prove convergence of
$W_{tn^{1/3}}$ and $Y_{tn^{1/3}}$ as $n\to\infty$, after suitable rescaling. 
A key ingredient
is the tightness of the random function
\begin{align}\label{hA}
\hA\nn_t:=n^{-1/2} A_{t n^{1/3}}\nn,
\qquad t\ge0.
\end{align}

Recall that $C[a,b]$ is the space of continuous functions on
$[a,b]$.

\begin{lemma}\label{LA2}
  Let $0<a <b<\infty$.
Then the stochastic processes $\hA\nn_t$, $n\ge1$,
are tight in $C[a,b]$.
\end{lemma}

The proof of \refL{LA2} is more complicated than for the corresponding
Lemma 5.2 in \cite{Janson2023}, and we show first two other lemmas.
We begin by stating a simple general lemma on tightness in the space $C[a,b]$.
(See e.g.\ \cite{Billingsley} for a background.)

\begin{lemma}\label{LC}
Let $-\infty<a<b<\infty$.
  Let\/ $(X_n(t))\nxoo$ and\/ $(Y_n(t))\nxoo$  be two sequences of random
  continuous functions on $[a,b]$.
Suppose that there exists a sequence
$(Z_n)\nxoo$ of random variables such that
for every $n$ and $s,t\in[a,b]$, we have
\begin{align}\label{lc}
  |X_n(t)-X_n(s)| \le Z_n |Y_n(t)-Y_n(s)|.
\end{align}
If the sequences 
$(X_n(a))\nxoo$ and $(Z_n)\nxoo$ are tight, and
$(Y_n(t))\nxoo$ is tight in $C[a,b]$, 
then 
the sequence $(X_n(t))\nxoo$ is tight in $C[a,b]$.
\end{lemma}

\begin{proof}
We may for convenience assume $[a,b]=\oi$.
We define 
for any function $f$ on \oi{} its modulus of continuity
\begin{align}
  \go(f;\gd):=\sup_{s,t\in\oi;\; |s-t|<\gd}|f(s)-f(t)|,
\qquad \gd>0.
\end{align}
Then \cite[Theorem 8.2]{Billingsley} says that a sequence $(X_n(t))\nxoo$
in $C\oi$
is tight  if and only if 
\begin{romenumerate}
\item\label{bill1} the sequence $(X_n(0))_n$ is tight, and
\item\label{bill2} 
for each positive $\eps$ and $\eta$, there exists $\gd>0$ such that,
for every $n$,
  \begin{align}
    \P\bigpar{\go(X_n;\gd)\ge\eps}\le\eta.
  \end{align}
\end{romenumerate}
We have already assumed \ref{bill1}.
Moreover, the assumption \eqref{lc} implies
\begin{align}\label{lc4}
\go(X_n;\gd)\le Z_n\go(Y_n;\gd)   
\end{align}
for every $\gd$.
Let $\eps,\eta>0$ be given.
Since the sequence $(Z_n)\nxoo$ is tight,
there exists a number $K>0$ such that $\P(|Z_n|>K) <\eta/2$ for every $n$.
Hence, \eqref{lc4} implies
\begin{align}\label{lc5}
\P\bigpar{\go(X_n;\gd)\ge\eps}&
\le \P(|Z_n|>K) + \P\bigpar{K\go(Y_n;\gd)\ge \eps}   
\notag\\&
\le \eta/2 + \P\bigpar{\go(Y_n;\gd)\ge \eps/K}.   
\end{align}
Since $(Y_n(t))\nxoo$ is tight, conditions \ref{bill1}--\ref{bill2} hold for
$Y_n(t)$; in particular, there exists $\gd>0$ such that
$\P\bigpar{\go(Y_n;\gd)\ge\eps/K}\le\eta/2$ for every $n$.
Then \eqref{lc5} shows that \ref{bill2} holds.
Consequently, $(X_n(t))\nxoo$ is tight.
\end{proof}

Recall that \refL{le:Ak} shows that 
\begin{align}\label{ha2}
0\le A_{k-1}-A_k
\le W_k^2\Phi_k\qw B_k^2
\le M_k^2\Phi_k\qw B_k^2,
\qquad 1\le k\le n-1. 
\end{align}
We  define the simpler
\begin{align}\label{ha3}
  V_k:=\sum_{j=1}^k B_j^2,\quad T_k :=\sum_{j=1}^k B_j,
\qquad 0\le k\le n-1; 
\end{align}
we extend also $V_k$ and $ T_k$ by linear interpolation to real arguments,
and define 
\begin{align}\label{hV}
\hV\nn_t:=n^{1/3} V_{t n^{1/3}}\nn, \quad\wh T\nn_t := T\nn_{tn^{1/3}},
\qquad t\ge0.
\end{align}
The proof of \refL{LA2} only uses $\hV\nn_t$, but $\wh T\nn_t$ is needed later when we prove \eqref{nov7}.

\begin{lemma}\label{LV2}
  Let $0<a<b<\infty$. Then the stochastic processes $\hV\nn_t-\hV\nn_a$ and  
$\wh T\nn_t- \wh T\nn_a$, $n\ge1$,
are tight in $C[a,b]$.
\end{lemma}

\begin{proof}
We start with $\hV\nn_t-\hV\nn_a$, $n\geq 1$.
If $1\le k\le \ell\le n-1$, then
\begin{align}\label{hdj1}
  \E|V_\ell- V_k|^2
=\E\Bigpar{ \sum_{i,j=k+1}^\ell B_i^2 B_j^2}
= \sum_{i,j=k+1}^\ell\E\bigsqpar{ B_i^2 B_j^2}
.\end{align}
If $i\neq j$, then $B_i$ and $B_j$ are independent, and thus, by
\eqref{eq:Ebetasq}, 
$\E\bigsqpar{ B_i^2 B_j^2}=\E[B_i^2]\,\E[B_j^2] = O\bigpar{i^{-2}j^{-2}}$.
On the other hand, if $i=j$, then, by \eqref{betamom}, 
recalling that $B_i\in \Beta(2,4i-6)$,
\begin{align}\label{hdj2}
  \E B_i^4 = \frac{2\cdot3\cdot4\cdot5}{(4i-4)(4i-3)(4i-2)(4i-1)}
= O\bigpar{i^{-4}}.
\end{align}
Consequently, \eqref{hdj1} yields
\begin{align}\label{hdj3}
  \E|V_\ell- V_k|^2
\le  \sum_{i,j=k+1}^\ell C i^{-2}j^{-2}
\le C (\ell-k)^2 k^{-4}
.\end{align}
This trivially holds for $\ell>n-1$ too, since then $V_\ell=V_{n-1}$ by
definition. Furthermore, 
writing \eqref{hdj3} as $\norm{V_\ell-V_k}_{L^2}\le C (\ell-k)k^{-2}$,
it follows by Minkowski's inequality that we can interpolate between integer
arguments and conclude that \eqref{hdj3} holds for all real $k$ and $\ell$
with $1\le k\le \ell$.

Let $s$ and $t$ be real numbers with $0<s\le t$.
Then \eqref{hV} and \eqref{hdj3} yield,
with $k:=sn^{1/3}$ and $\ell:=tn^{1/3}$, 
\begin{align}\label{hdj4}
 \E\bigabs{\hV\nn_t-\hV\nn_s}^2
= n^{2/3} \E |V_\ell-V_k|^2
\le C n^{2/3}|\ell-k|^2 k^{-4}
= C|t-s|^2 s^{-4}.
\end{align}
For the restriction to $[a,b]$ we thus have 
\begin{align}\label{hdj4b}
 \E\bigabs{\hV\nn_t-\hV\nn_s}^2
\le C_a|t-s|^2 ,
\qquad a\le s\le t\le b,
\end{align}
which shows the tightness of $\hV\nn_t-\hV\nn_a$ by 
\cite[Theorem 12.3]{Billingsley}.

Tightness of $\wh T\nn_t-\wh T\nn_a $ can be shown by using 
\begin{align}
    \E| T_k- T_\ell|^2 =\sum^\ell_{i,j=k+1} \E(B_iB_j) \leq C(\ell-k)^2 k^{-2}, \quad 1\leq k\leq \ell\leq n-1,
\end{align}
instead of \eqref{hdj3}
and proceeding similarly.
\end{proof}

\begin{proof}[Proof of \refL{LA2}]
Note first that \eqref{eq:EAk} implies that 
\begin{align}
  \E \hA\nn_a \le  n^{-1/2}\frac{Cn}{(an^{1/3})^{3/2}} = Ca^{-3/2},
\end{align}
and thus the sequence $\hA\nn_a$ is tight.

Let 
\begin{align}\label{hdj0}
  \xM_n:=n^{-1/2}\max_{k\ge1} M_k
\qquad\text{and}\qquad
\Phix_n:= n^{1/6}\Phi\qw_{\floor{a n^{1/3}}}. 
\end{align}
Then \refL{le:doob} 
shows that
$\E\xM_n^2\le C$, 
and \eqref{eq:phibd} yields $\E\Phix_n=O(1)$;
hence the sequences $(\xM_n)\nxoo$ and $(\Phix_n)\nxoo$ are tight.

By \eqref{ha2} and \eqref{ha3}, we have for any integer
$k$ with $k\ge a n^{1/3}$,
since $\Phi_k$ is increasing by the definition \eqref{de:phi},
using \eqref{hdj0},
  \begin{align}\label{hdj01}
   | A_k- A_{k-1}| &
\le M_k^2\Phi_k\qw (V_k-V_{k-1})
\le n\xM_n^2 \Phi\qw_{\floor{a n^{1/3}}} (V_k-V_{k-1})
\notag\\&
= n^{5/6}\xM_n^2 \Phix_n (V_k-V_{k-1}).
  \end{align}
Thus, if $\floor{a n^{1/3}}\le k\le \ell$,
  \begin{align}\label{hdj5}
   | A_\ell- A_{k}| 
\le n^{5/6}\xM_n^2 \Phix_n (V_\ell-V_{k}),
  \end{align}
We can interpolate between integer arguments and conclude that \eqref{hdj5}
holds for all real $k$ and $\ell$ with $\floor{a n^{1/3}}\le k\le \ell$.
By \eqref{hA} and \eqref{hV}, this shows that if $a\le s\le t$, then 
  \begin{align}\label{hdj6}
\bigabs{\hA\nn_t-\hA\nn_s}
\le \xM_n^2 \Phix_n \bigpar{\hV\nn_t-\hV\nn_s}.
  \end{align}
The result now follows from \refLs{LC} and \ref{LV2}, 
taking $X_n(t):=\hA\nn_t$, $Y_n(t):=\hV\nn_t-\hV\nn_a$, and
$Z_n:=\xM_n^2\Phix_n$, 
noting that $(Z_n)\nxoo$ is tight since
$(\xM_n)\nxoo$ and
$(\Phix_n)\nxoo$ are.
\end{proof}

\section{The final part: convergence}\label{Sconv}

Recall the definition of the spaces $C(0,\infty)$ and $C[0,\infty)$ in
Section \ref{Snot}. 

\begin{theorem}\label{th:cvg}
As $n\to\infty$ we have
    \begin{align}\label{cvg3}
  \frac{W_{tn^{1/3}}}{n^{1/2}}&\dto 
4\tilde \beta
   \bclr{\bclr{t^{9/2}+\tfrac{3}{4}\xi t^{3}}^{1/3}-t^{3/2}}
\qquad\text{in $C[0,\infty)$},
\intertext{and}
\label{cvg4}
  \frac{Y_{tn^{1/3}}}{n^{1/3}}&\dto 
4 \bclr{\bclr{t^{3}+\tfrac{3}{4}\xi t^{3/2}}^{1/3}-t}
\qquad\text{in $C(0,\infty)$},
    \end{align}
with $\tilde\beta$ as in \eqref{lb2b} and\/ $\xi\in\mathrm{Gamma}(2,1)$ 
independent of $\tilde \beta$.
\end{theorem}
\begin{remark}
We believe that \eqref{cvg4} holds also in $C[0,\infty)$,
but we see no simple proof so we leave this as an open problem.
\end{remark}

\begin{proof}
The proof is similar to that of \cite[Theorem 5.3]{Janson2023}, but with
several technical complications. 

\resetsteps
\steppx{Convergence in $C(0,\infty)$ for a subsequence.}
  As in \cite{Janson2023},
\refL{LA2} implies that by considering a subsequence,
we may assume that the processes $\hA\nn_t$ converge in distribution in
every space $C[a,b]$ with $0<a<b$ to some continuous stochastic process
$\cA(t)$ on $(0,\infty)$; in other words, as \ntoo,
\begin{align}
  \label{aw1}
\hA\nn_t \to \cA(t)
\end{align}
holds in distribution in the space $C(0,\infty)$.
Furthermore, also as in \cite{Janson2023}, we may
by the Skorohod coupling theorem \cite[Theorem~4.30]{Kallenberg}, 
assume
that this convergence holds a.s.;
in other words (as $\ntoo$ along the subsequence)
a.s.\ \eqref{aw1} holds
uniformly on every interval $[a,b]$.
We use such a coupling until further notice.
Note that this means that we consider all random variables as defined separately
for each $n$ (with some unknown coupling); in particular, this means that we
have potentially a different sequence $B_i\nn$ ($i\ge1$) for each $n$, 
and thus different limits $\gb\nn$ and $\tgb\nn$.
(The variables $B_i\nn$ for a fixed $n$ are independent, but we do not know
how they are coupled for different $n$.)
Hence, we cannot directly use the a.s.\ convergence results in \refL{LB2}.
Instead we note that \eqref{lb2b} (which holds for each $n$, with the
distributions of the variables the same for all $n$) implies,
for any coupling,
\begin{align}\label{jup1}
\sup_{k\ge\log n}\bigabs{k^{-1/2}\Phi\nn_k- \tgb\nn}\pto0.
\end{align}
Furthermore, trivially (since the distributions are the same)
\begin{align}\label{jup2}
  \tgb\nn\dto\tgb
\end{align}
for some random variable $\tgb$ with, by \refL{LB2}, $\tgb>0$ a.s.;
note also that \eqref{jup2} holds jointly with \eqref{lxi}, since this is
true for the coupling used in the proof of \refL{LXi},
where $B_i$ does not depend on $n$ and thus trivially
$\tgb\nn\to\tgb$ holds together with 
\eqref{lxi3}. 

We may select the subsequence above such that \eqref{aw1} (in distribution),
\eqref{jup2} and \eqref{lxi} hold jointly 
(with some joint distribution of the limits).
We may then assume that \eqref{aw1}, \eqref{jup1}, \eqref{jup2}, \eqref{lxi},
\eqref{lfM},  and \eqref{tf1} all hold a.s., by
redoing the application of the Skorohod coupling theorem and
including all these limits.
It then follows from
\eqref{eq:doobdecom}, \eqref{hA}, \eqref{lfM}, \eqref{aw1},
and \eqref{lxi} that a.s.
\begin{align}\label{aw2}
  n^{-1/2}W_{tn^{1/3}}=
n^{-1/2}M_{tn^{1/3}}-\hA\nn_t
=\Xi\nn -\cA(t) + o(1)
\to \cB(t):=\tgb\xi-\cA(t)
\end{align}
uniformly on every compact interval in $(0,\infty)$.
In other words, \eqref{aw2} holds a.s.\ in $C(0,\infty)$.
From \eqref{aw2} we obtain by \eqref{de:Wk}, \eqref{jup1}, and \eqref{jup2},
letting $k:=\floor{tn^{1/3}}$, a.s.
\begin{align}\label{aw3}
 Y_{\floor{tn^{1/3}}}=\frac{W_k}{\GF_k}=\frac{W_k}{k^{1/2}}(\tgb\nn+o(1))\qw
=\frac{n^{1/2}}{k^{1/2}}\bigpar{\cB(t)+o(1)}\bigpar{\tgb+o(1)}\qw
\end{align}
and thus
\begin{align}\label{aw4}
n^{-1/3} Y_{\floor{tn^{1/3}}}
\to t^{-1/2}\tgb\qw\cB(t),
\end{align}
again uniformly on every compact interval in $(0,\infty)$,
and thus in $C(0,\infty)$.

\steppx{Identifying the limit.}
Fix $0<s<t$ and let $k:=\floor{sn^{1/3}}$ and $\ell:=\floor{tn^{1/3}}$.
It follows from \eqref{aw1} that $\hA\nn_s-\hA\nn_{n^{-1/3}k}\to0$ and
$\hA\nn_t-\hA\nn_{n^{-1/3}\ell}\to0$ a.s.
Hence, \eqref{hA} and \eqref{eq:incA} imply that
\begin{align}\label{aw5}
  \hA\nn_s-\hA\nn_t 
&= \sum_{i=k+1}^\ell n^{-1/2}(A_i-A_{i-1})+o(1)
\notag\\&
=n^{-1/2}\sum_{i=k+1}^\ell2\GF_{i-1}\bigsqpar{(1-B_i)^{Y_i}-1+Y_iB_i}+o(1)
.\end{align}
For any $B\in(0,1)$, with $x:=-\log(1-B)$,
and any $y\ge1$,
\begin{align}\label{aw6}
  \frac{\ddx}{\ddx y}\bigpar{(1-B)^{y}-1+yB}
&=
  \frac{\ddx}{\ddx y}\bigpar{e^{-xy}-1+y(1-e^{-x})}
=-xe^{-xy}+1-e^{-x}
\notag\\&
\ge 1-(1+x)e^{-x}>0.
\end{align}
Hence $(1-B)^{y}-1+yB$ is an increasing function of $y\ge1$.
Let $\eps>0$, and let
\begin{align}\label{aw7}
  y_+:=\max\set{u^{-1/2}\tgb\qw\cB(u): u\in[s,t]}+\eps.
\end{align}
Then \eqref{aw4} implies that, for large enough $n$, $Y_{i}\le y_+n^{1/3}$
when $k\le i\le\ell$, and hence \eqref{aw5} implies,
noting that $\GF_i$ is increasing in $i$,
\begin{align}\label{aw8}
  \hA\nn_s-\hA\nn_t &
\le2 n^{-1/2}\GF_{\ell}\sum_{i=k+1}^\ell
 \bigsqpar{(1-B_i)^{n^{1/3}y_+}-1+n^{1/3}y_+B_i}+o(1)
.\end{align}
Using the notation \eqref{lh1}--\eqref{lh2}, 
we thus have by \eqref{lh3} and \eqref{lb2b}
\begin{align}\label{aw9}
  \hA\nn_s-\hA\nn_t &
\le2 n^{-1/6}\GF_{\ell}\hH^{y_+}_{s,t}+o(1)
\notag\\&
= 2 t^{1/2}\tgb
\int_s^t \Bigpar{\frac{1}{(1+y_+/(4u))^2}-1+\frac{y_+}{2u}}\dd u +\op(1)
.\end{align}
Similarly,
defining
\begin{align}\label{aw7-}
  y_-:=\min\set{u^{-1/2}\tgb\qw\cB(u): u\in[s,t]}-\eps
\end{align}
(adjusted to 0 if \eqref{aw7-} becomes negative), we obtain a lower bound
\begin{align}\label{aw9-}
  \hA\nn_s-\hA\nn_t &
\ge 2 s^{1/2}\tgb
\int_s^t \Bigpar{\frac{1}{(1+y_-/4u)^2}-1+\frac{y_-}{2u}}\dd u +\op(1)
.\end{align}
We now subdivide $[s,t]$ into a large number $N$ 
of small subintervals of equal length and apply 
\eqref{aw9} and \eqref{aw9-} for each subinterval $[s_i,t_i]$.
Since $u^{-1/2}\tgb\qw\cB(u)$ is continuous, we may choose $N$ such that with
probability $>1-\eps$, for each subinterval, the corresponding values of
$y_+$ and $y_-$ differ by at most $3\eps$, and also that $t_i/s_i<1+\eps$.
We may then, for each subinterval, replace $y_+$ and $y_-$ in \eqref{aw9}
and \eqref{aw9-} by $u^{-1/2}\tgb\qw\cB(u)$ with a small error, and by summing
over all subintervals it finally follows  by letting
$\eps\to0$ (we omit the routine details) that
\begin{align}\label{aw10}
  \hA\nn_s-\hA\nn_t &
= 2 \tgb
\int_s^t u^{1/2}\Bigpar{\bigpar{1+\tfrac{1}{4}u^{-3/2}\tgb\qw\cB(u)}^{-2}
-1+\tfrac12u^{-3/2}\tgb\qw\cB(u)}\dd u +\op(1)
.\end{align}
Since we also assume \eqref{aw1}, the left-hand side converges a.s.\ to
$\cA(s)-\cA(t)$, and thus we have, a.s.,
\begin{align}\label{aw11}
  \cA(s)-\cA(t) &
= 2 \tgb
\int_s^t u^{1/2}\Bigpar{\bigpar{1+\tfrac{1}{4}u^{-3/2}\tgb\qw\cB(u)}^{-2}
-1+\tfrac12u^{-3/2}\tgb\qw\cB(u)}\dd u
.\end{align}
This holds a.s.\ simultaneously for every pair of rational  $s$ and $t$, and
thus by continuity a.s.\ for every real $s$ and $t$ with $0< s\le t$. 
Consequently, $\cA(t)$ is a.s.\ continuously differentiable on $(0,\infty)$,
with 
\begin{align}\label{aw12}
\cA'(t) &
= -2 \tgb
 t^{1/2}\Bigpar{\bigpar{1+\tfrac{1}{4}t^{-3/2}\tgb\qw\cB(t)}^{-2}
-1+\tfrac12t^{-3/2}\tgb\qw\cB(t)}
.\end{align}
By the definition of $\cB(t)$ in \eqref{aw2}, this implies that also
$\cB(t)$ is continuously differentiable and
\begin{align}\label{aw12B}
\cB'(t) &
=-\cA'(t)
= 2 \tgb
 t^{1/2}\Bigpar{\bigpar{1+\tfrac{1}{4}t^{-3/2}\tgb\qw\cB(t)}^{-2}
-1+\tfrac12t^{-3/2}\tgb\qw\cB(t)}
.\end{align}
We may simplify a little by defining
\begin{align}\label{tcB}
  \tcB(t):=\tgb\qw\cB(t).
\end{align}
Then \eqref{aw12B} becomes
\begin{align}\label{aw13}
\tcB'(t) &
= 2 t^{1/2}\Bigpar{\bigpar{1+\tfrac{1}{4}t^{-3/2}\tcB(t)}^{-2}
-1+\tfrac12t^{-3/2}\tcB(t)}
.\end{align}

By definition, $\hA\nn_t$ is decreasing on $[0,\infty)$, and thus 
\eqref{aw1} shows that $\cA(t)$ is decreasing, and thus $\cB(t)$ is
increasing by \eqref{aw2}. (This also follows from \eqref{aw13}, since the
right-hand side is positive.)
Furthermore, \eqref{aw1}, \eqref{hA}, \eqref{eq:EAk}, and Fatou's inequality
yield, for every $t>0$,
\begin{align}
  \E \cA(t) \le\liminf_\ntoo\E\hA\nn_t
\le\liminf_\ntoo n^{-1/2}\frac{Cn}{(tn^{1/3})^{3/2}}
=\frac{C}{t^{3/2}}.
\end{align}
Hence, by dominated convergence,
\begin{align}
  \E\lim_{t\to\infty}\cA(t)=\lim_{t\to\infty}\E\cA(t)=0.
\end{align}
Consequently, a.s.\ $\cA(t)\to0$ as \ttoo, and thus by \eqref{tcB} and
\eqref{aw2} 
\begin{align}\label{aw16}
  \tcB(t)\upto\xi,
\qquad\text{as \ttoo}.
\end{align}

We show in \refApp{Sdiff} below, see in particular \eqref{eq:de3} and
\eqref{eq:gc0}, 
that the differential equation 
\eqref{aw13} has a unique solution satisfying 
the boundary condition \eqref{aw16},
viz.\
\begin{align}\label{tBt}
    \tcB(t) = 4t^{3/2} \bclr{\bclr{1+\tfrac{3}{4} \xi t^{-3/2}}^{1/3}-1},
\qquad t>0.
\end{align}
(It can easily be verified by differentiation that this is a solution; 
\refApp{Sdiff} shows how the solution may be found, and that it is unique.)
Hence, by \eqref{tcB}, 
\begin{align*}
    \cB(t) &= \tgb\tcB(t) 
= 4\tilde \beta t^{3/2} \bclr{\bclr{1+\tfrac{3}{4} \xi t^{-3/2}}^{1/3}-1}\\
    &=4\tilde \beta\bclr{\bclr{t^{9/2}+\tfrac{3}{4} \xi t^{3}}^{1/3}-t^{3/2}}.\numberthis\label{Bt}
\end{align*}

\steppx{Convergence in $C[0,\infty)$.}
Note that \eqref{Bt} shows that $\cB(t)$ extends to a continuous function on
$[0,\infty)$ with $\cB(0)=0$; hence it follows from \eqref{aw2} that
also $\cA(t)$ extends to a continuous function on $[0,\infty)$ with 
$\cA(0)=\tilde \beta \xi$. Using $A\nn_0 =M\nn_0 - W\nn_0\leq M\nn_0$,
and the assumed a.s.\ versions of \eqref{lfM} and \eqref{lxi}, 
we have that, a.s., 
\begin{align}
    \limsup_{n\to\infty} n^{-1/2} A\nn_0 \leq  \limsup_{n\to\infty} n^{-1/2} M\nn_0 =\tilde\beta \xi = \cA(0).
\end{align}
By the reverse increasing property of $A\nn_k$, we also have,
for every $t>0$, 
\begin{align}
    \liminf_{n\to\infty} n^{-1/2} A\nn_0\geq \liminf_{n\to\infty} n^{-1/2} A\nn_{tn^{1/3}}  =\cA(t).
\end{align}
Sending $t\downto 0$ and thus $\cA(t)\upto \cA(0)$ then yields
\begin{align}
    \liminf_{n\to\infty} n^{-1/2} A\nn_0\geq \cA(0).
\end{align}
It follows that, $\wh A\nn_0 \to \cA(0)$ a.s., 
and thus \eqref{aw1} holds a.s.\ for every fixed $t\ge0$.
Since $\hA\nn_t$ and $\cA(t)$ are decreasing in $t$, and $\cA(t)$ is
continuous, this implies that \eqref{aw1} holds a.s.\ uniformly for every
interval $[0,b]$ with $0<b<\infty$. 
It then follows that also \eqref{aw2} holds a.s.\ uniformly on every compact
interval in $[0,\infty)$;
in other words, \eqref{aw2} holds a.s.\ in $C[0,\infty)$,
with $\cB(t)$ as in \eqref{Bt}.

\steppx{Conclusion.}
We have so far considered a subsequence,
and a special coupling, and have then shown \eqref{aw2} in $C[0,\infty)$ and
\eqref{aw4} in $C(0,\infty)$, which by \eqref{Bt} 
yields \eqref{cvg3} and \eqref{cvg4}.
Since \eqref{cvg3} and \eqref{cvg4} use convergence
in distribution, they hold in general along the subsequence, also without
the special coupling used in the proof. 
Moreover, the limits in \eqref{cvg3} and \eqref{cvg4} do not depend on the
subsequence, and the proof shows that every subsequence has a subsequence 
such that the limits in distribution \eqref{cvg3} and \eqref{cvg4} hold.
As is well known, this implies that the full sequences converge in
distribution, see e.g.\ \cite[Section 5.7]{Gut}).
\end{proof}

\begin{remark}
  The argument above using possibly different $B_i\nn$ is rather
  technical. A more elegant, and perhaps more intuitive version, 
of the argument would be to assume
  that $B_i$ is the same for all $n$, and then condition on $(B_i)_{i=1}^\infty$
  before applying the Skorohod coupling theorem.  However, while
  intuitively clear, this seems technically more difficult to justify, and
  it seems to require that we prove that earlier convergence results hold
  also conditioned on $(B_i)_{i=1}^\infty$, a.s. We therefore prefer the 
somewhat clumsy argument above.
\end{remark}

\section{The number of descendants}\label{se:desc}
Let $X=X^{(n)}$ be the number of red vertices in the preferential attachment
graph. Vertex $n$ is red by definition, and $J_k=\tone[Z_k\geq 1]$ is the
indicator that takes value 1 if vertex $k$ is red for $k< n$; thus 
\begin{equation}
    X=1+\sum^{n-1}_{k=1} J_k,
\end{equation}
noting that $J_k$ is $\cF_{k-1}$-measurable. We now set out to prove \eqref{nov7} (a special case of \refT{Tmain} with $m=2$ and $\rho=0$), which we for convenience restate below as a separate theorem.
\begin{theorem}\label{th:X}
    As $n\to\infty$, 
    \begin{equation}\label{eX}
        \frac{X^{(n)}}{n^{1/3}}\dto 
2^{-\xfrac43}3^{-\xfrac13}\frac{\gG(\frac13)^2}{\gG(\frac23)}
\xi^{2/3}  ,
    \end{equation}
    where $\xi\in \mathrm{Gamma}(2,1)$.
\end{theorem}

\begin{remark}
    Unlike in \eqref{cvg3} and \eqref{cvg4}, $\tgb$ does not appear in the distributional limit of $X^{(n)}/n^{1/3}$. This is because $\beta$ in \eqref{lb2a} is essentially
determined by $B_k$ corresponding to the first few vertices; and in the red subgraph $D_n$, the number of these vertices is insignificant, since most vertices have labels of the order $n^{1/3}$.
\end{remark}

As in \cite{Janson2023}, we use the Doob decomposition 
\begin{equation}\label{eq:X}
    X=1+L_0+P_0,
\end{equation}
where for $ 0\leq k\leq n-1$, 
\begin{equation}\label{eq:Xm}
    L_k := \sum^{n-1}_{i=k+1}\bigpar{J_i-\E(J_i\mid \cF_i)}
\end{equation}
is a reverse martingale, and by \eqref{eq:distZk}, 
\begin{equation}\label{eq:Xpred}
    P_k := \sum^{n-1}_{i=k+1} \E(J_i\mid \cF_i) = \sum^{n-1}_{i=k+1} \IP(Z_i\geq 1\mid \cF_i) = \sum^{n-1}_{i=k+1} \bbclr{1-(1-B_i)^{Y_i}}  
\end{equation}
is positive and reverse increasing. Furthermore, by Markov's inequality and \eqref{eq:distZk}, 
\begin{equation}\label{eq:Xpred1}
    P_{k-1}-P_k =\IP(Z_k\geq 1\mid \cF_k) \leq  B_kY_k,  \quad 1\leq k\leq n-1. 
\end{equation}

By \eqref{eq:Xpred} and Lemma \ref{le:Zk}, for $1\leq k\leq n-1$,
\begin{equation}\label{eq:EPk}
    \E P_k = \sum^{n-1}_{i=k+1} \IP(Z_i\geq 1) \leq  \sum^{n-1}_{i=k+1}  \frac{Cn^{1/2}}{i^{3/2}} \leq \frac{Cn^{1/2}}{k^{1/2}}.
\end{equation}
Using also the crude bound $0\leq J_i\leq 1$, we have $P_0-P_k\leq k$. Choosing $k=\floor{n^{1/3}}$ and applying \eqref{eq:EPk} thus yield
\begin{equation}\label{eq:EP0}
    \E P_0  \leq \E P_{\floor{n^{1/3}}}+ \floor{n^{1/3}} \leq Cn^{1/3}.
\end{equation}
Moreover, it follows from the reverse martingale property of $L_k$, $\var(J_i\mid \cF_i)\leq \E(J_i\mid \cF_i)$ and \eqref{eq:EP0} that
\begin{equation}\label{eq:L02}
    \E L_0^2 = \sum^{n-1}_{i=1} \E[\var(J_i\mid \cF_i)] 
\le\sum^{n-1}_{i=1} \E J_i  =\E P_0 \leq Cn^{1/3}.
\end{equation}
This in turn implies that as $n\to\infty$, 
\begin{equation}\label{eq:L0}
    \frac{L_0}{n^{1/3}} \overset{p}{\longrightarrow} 0,
\end{equation}
and thus by \eqref{eq:X}, it is enough to show that as $n\to\infty$, $n^{-1/3}P_0$ converges in distribution to the \rhs{} of \eqref{eX}. To this end, we extend also $P_k$ to real arguments by linear interpolation and define 
    \begin{equation}\label{de:hP}
        \wh P^{(n)}_t=n^{-1/3} P^{(n)}_{tn^{1/3}},\qquad t\ge 0.
    \end{equation}

\begin{lemma}\label{le:Pt}
    Let $0<a<b<\infty$. Then the stochastic processes $\wh P^{(n)}_t$, $n\geq 1$, are tight in $C[a,b]$.
\end{lemma}

\begin{proof}
The proof is very similar to that of Lemma \ref{LA2}. 
First,
by \eqref{de:hP} and \eqref{eq:EP0},
\begin{align}
  \E \hP\nn_a =  n^{-1/3}\E P\nn_{an^{1/3}}
\le n^{-1/3}\E P\nn_0
\le C,
\end{align}
and thus the sequence $(\hP\nn_a)\nxoo$ is tight.

Let $\xM_n$ and $\Phix_n$ be as in \eqref{hdj0}, $T_k$ and $\wh T\nn_k$ be as in \eqref{ha3} and \eqref{hV}. From \eqref{eq:Xpred1} and \eqref{de:Wk}, $W_k\leq M_k$, the increasing property of $\Phi_k$ and also \eqref{ha3}, for any integer $k\geq an^{1/3}$, 
    \begin{align}\label{eq:hP}
        |P_k-P_{k-1}| \leq  B_k Y_k = B_k\Phi_k^{-1} W_k 
\leq B_k\Phi_{\floor{an^{1/3}}}^{-1}\max_{k\geq 1} M_k
        = n^{1/3}  \xM_n \Phix_n ( T_k- T_{k-1}). 
    \end{align}
    Extending \eqref{eq:hP} to real arguments and using  \eqref{de:hP}, we thus have 
\begin{align}
     |\wh P\nn_t-\wh P\nn_s|\leq  \xM_n \Phix_n (\wh T\nn_t-\wh T\nn_s)\quad
  \text{if $a\leq s\leq t$.}
\end{align}
The result then follows from Lemmas \ref{LC} and \ref{LV2}, this time taking $X_n(t):=\wh P\nn_t$, $Y_n(t):=\wh T\nn_t-\wh T\nn_a$ and $Z_n:=\xM_n \Phix_n$; tightness of $\xM_n \Phix_n$ follows from that of $(\xM_n)^\infty_0$ and $ (\Phix_n)^\infty_0$. 
\end{proof}

\begin{proof}[Proof of \refT{th:X}]
In view of Lemma \ref{le:Pt}, by considering a subsequence, we may assume
that the processes $\wh P\nn_t$ converge in distribution in every space
$C[a,b]$ for $0<a<b<\infty$, 
and thus in $C(0,\infty)$, 
to some stochastic process $\mathcal{P}(t)$ on
$(0,\infty)$. Again using the Skorohod coupling theorem, we can assume that
all a.s.\ convergence results in the proof of Theorem \ref{th:cvg} hold and
also  
\begin{equation}\label{eq:Pc}
    \wh P\nn_t \to \mathcal{P}(t)
\end{equation}
a.s.\ uniformly on every interval $[a,b]$. From \eqref{aw4}, 
    \begin{equation}
        n^{-1/3} Y_{\floor{tn^{1/3}}} 
\to t^{-1/2}\tilde \beta\qw \mathcal{B}(t)
= t^{-1/2}\tcB(t)
    \end{equation}
a.s.\ uniformly on each compact interval in $(0,\infty)$.

Let $s,t$ be real numbers with $0<s<t$, and let $k=\floor{sn^{1/3}}$ and
$\ell=\floor{tn^{1/3}}$. By the same argument leading to \eqref{aw5}, now
using \eqref{eq:Pc} and \eqref{eq:Xpred}, we have 
\begin{align}
    \wh P\nn_s -  \wh P\nn_t =n^{-1/3} \sum^\ell_{i=k+1}\bclr{1-(1-B_i)^{Y_i}} + o(1)
\end{align}
Let $y_+$ and $y_-$ be as in \eqref{aw7} and \eqref{aw7-}. By \eqref{lh34} of Lemma \ref{LH}, we can therefore conclude that 
\begin{align}\label{eq:P+}
    \wh P\nn_s -  \wh P\nn_t \leq \int^t_s \bbclr{1-\frac{1}{(1+y_+/(4u))^2}}\dd u +\op(1)
\end{align}
and 
\begin{align}\label{eq:P-}
     \wh P\nn_s -  \wh P\nn_t \geq \int^t_s \bbclr{1-\frac{1}{(1+y_-/(4u))^2}}\dd u +\op(1).
\end{align}
The sandwich argument in the proof of \eqref{aw10}, the bounds in
\eqref{eq:P+} and \eqref{eq:P-}, together with \eqref{tBt}, imply that for $0<s<t<\infty$,
\begin{align}
     \wh P\nn_s -  \wh P\nn_t 
&=\int_s^t\bbclr{1-\frac{1}{(1+u^{-1/2}\tcB(u)/(4u))^2}}\dd u +\op(1)
\notag\\&
=\int^t_s \bbclr{1-\frac{1}{\bclr{1+\tfrac{3}{4} \xi u^{-3/2}}^{2/3}}}\dd u 
+\op(1).
\end{align}
In light of \eqref{eq:Pc}, 
we thus have a.s.,
\begin{align}\label{eq:Pc1}
    \mathcal{P}(t)- \mathcal{P}(s)
=\int^t_s \bbclr{1-\frac{1}{\bclr{1+\tfrac{3}{4} \xi u^{-3/2}}^{2/3}}}\dd u. 
\end{align}

Furthermore,
\begin{equation}
    \frac{\dd}{\dd s} \wh P\nn_s = -\E(J_k\mid \cF_k) = (1-B_k)^{Y_k} -1,
\end{equation}
where $k=\ceil{sn^{1/3}}$. Thus
$ \big| \tfrac{\dd}{\dd s} \wh P\nn_s\big|\le 1$, which in turn implies that
\begin{align}\label{eq:P1}
  \big| \wh P\nn_0-  \wh P\nn_s\big| \leq s.
\end{align}
For $t\geq 1$, it follows from the reverse increasing property of $P_k$ and \eqref{eq:EPk} that 
\begin{equation}\label{eq:P2}
    \E \wh P\nn_t = n^{-1/3} \E  P\nn_{tn^{1/3}} \leq n^{-1/3} \E  P\nn_{\floor{tn^{1/3}}} \leq  \frac{Cn^{1/6}}{\floor{tn^{1/3}}^{1/2}} \leq Ct^{-1/2}.
\end{equation}
Sending $s\to 0$ and $t\to\infty$, we deduce from \eqref{eq:P1} and \eqref{eq:P2} that
$\wh P\nn_0-(\wh P\nn_s - \wh P\nn_t)\pto 0$, uniformly in $n$. Combining this and \eqref{eq:Pc1} with a standard argument gives
\begin{equation}\label{eq:Pc2}
    \wh P\nn_0 \pto 
\int^\infty_0  \bbclr{1-\frac{1}{\bclr{1+\tfrac{3}{4} \xi u^{-3/2}}^{2/3}}}\dd u 
. \end{equation}
The change of variable $x=3\xi u^{-3/2}/4$ yields, using \refL{Lnora}, 
\begin{align}\label{px}
\int^\infty_0  \bbclr{1-\frac{1}{\bclr{1+\tfrac{3}{4} \xi u^{-3/2}}^{2/3}}}\dd u 
&=
\frac{2}{3}\Bigpar{\frac{3\xi}{4}}^{2/3}
\int^\infty_0  \bbclr{1-\frac{1}{(1+x)^{2/3}}}x^{-5/3}\dd x 
\notag\\&
=\frac{2}{3}\Bigpar{\frac{3}{4}}^{2/3}\xi^{2/3}\frac{-\gG(-\frac23)\gG(\frac43)}{\gG(\frac23)}
\notag\\&
=2^{-\xfrac43}3^{-\xfrac13}\frac{\gG(\frac13)^2}{\gG(\frac23)}\xi^{2/3}.
\end{align}
The result thus follows from 
\eqref{eq:X},
\eqref{eq:L0},
\eqref{de:hP},
\eqref{eq:Pc2},
and
\eqref{px}.
\end{proof}

\section{The general case}\label{Sgen}
Here we consider the general case. The argument is similar to the case where $m=2$ and $\rho=0$, so we give only the main changes here. 

\subsection{The stochastic recursions and new estimates}
Define $Y_k, J_k, Z_k, \cF_k$ as in Section \ref{se:basicanalysis}, but now
with the boundary condition $Y_{n-1}=m$. 
We use the stochastic recursions in \refS{sse:sr} to obtain the subgraph $D_n$, where we now sample $m$ outgoing edges instead of two. The recursion in \eqref{eq:stocrecur} now becomes
\begin{align}\label{eq:sr}
    Y_{k-1} = Y_k - Z_k + m\cdot J_k=Y_k - Z_k + m\cdot\mathbf{1}[Z_k\ge 1],
\qquad 2 \le k \le n-1. 
\end{align}
As \eqref{eq:distZk} still holds, we have
\begin{align}
    \E(Y_{k-1}\mid \cF_{k}) &= Y_k - \E(Z_k\mid \cF_k) + m \IP(Z_k\ge 1\mid \cF_k)\notag \\
    &= Y_k - B_k Y_k + m \big(1-(1-B_k)^{Y_k}\big),
\end{align}
and, again by Markov's inequality, 
\begin{align}\label{eq:gmy}
     \E(Y_{k-1}\mid \cF_{k}) \le Y_k  - B_k Y_k + m B_k Y_k = (1+(m-1)B_k) Y_k.
\end{align}
 Thus, with  $\Phi_k$ as in \eqref{de:phi}, we can define 
\begin{align}\label{de:gmw}
    W_k := \Phi_k Y_k,\qquad 0\le k\le n-1.
\end{align}
It follows from \eqref{eq:gmy} and \eqref{de:gmw} that $W_0,\dots W_{n-1}$ is a reverse supermartingale with 
\begin{align}
    W_{n-1} = \Phi_{n-1} Y_{n-1} = m\Phi_{n-1.}
\end{align}
We again consider the Doob decomposition $W_k=M_k-A_k$, with 
\begin{align}\label{yngve}
    M_k := m\Phi_{n-1} + \sum^{n-1}_{j=k+1} (W_{j-1}-\E(W_{j-1}\mid \cF_j) ),
\end{align}
and $A_k$ as in \eqref{eq:incp}. Analogous to \eqref{eq:incA} and \eqref{eq:ubvarW}, 
\begin{align}\label{eq:AA}
    A_{k-1}-A_k = m\Phi_{k-1} \big((1-B_k)^{Y_k}-1+B_k Y_k \big) 
\end{align}
and 
\begin{align}\label{eq:gvw}
\var(W_{k-1}\mid \cF_k ) \le C\Phi_{k-1}^2 B_k Y_k.     
\end{align}
Using \eqref{eq:gvw} and arguing as in \eqref{eq:varM}, we obtain
\begin{align} \label{eq:vm}
    \var_\bB (M_k) 
\leq C \sum^{n-1}_{j=k+1} \Phi_{j-1}^2 B_j \prod^{n-1}_{i=j+1} \bclr{1+(m-1)B_i}.
\end{align}

By the same proofs as for Lemmas \ref{le:doob}--\ref{le:Ak}, again using  \eqref{eq:Bmean}, \eqref{eq:Ebetasq}, \eqref{eq:mprodbeta}, \eqref{eq:meanphik} and \eqref{eq:phibd}, we get
\begin{gather}
    \E\max_{0\le k\le n-1} W_k^2 \le \E\max_{0\le k\le n-1} M_k^2 \le Cn^{2(m-1)\chi}; \label{eq:e1}\\
   \IP(Z_k\ge 1) \le \frac{Cn^{(m-1)\chi}}{k^{1+(m-1)\chi}},\qquad 
    \IP(Z_k\ge 2) \le \frac{Cn^{2(m-1)\chi}}{k^{2+2(m-1)\chi}};\label{eq:e2}\\
    A_{k-1}-A_k \le C \Phi_k^{-1} W_k^2 B_k^2, \qquad 
     \E A_{k} \le \frac{Cn^{2(m-1)\chi}}{k^{1+(m-1)\chi}}.\label{eq:e3}
\end{gather}

\subsection{The branching process}\label{se:Yule1}
 As in Section \ref{se:Yule}, we can couple the early part of $D_n$ to a
 suitable time-changed branching process. Let $\cY$ be a branching process
 that starts with $m$ particles at time 0, and each particle has an
 independent $\mathrm{Exp}(1)$ lifetime, before splitting into $m$ new
 particles. Let $\wh \cY$ be the time-changed counterpart of $\cY$, again by
 the mapping $t\mapsto e^{-t}$; 
thus $\wh \cY_x=\cY_{-\log x}$ is the number of particles in $\wh \cY$ alive at time $x$. By standard properties of branching 
 processes (see \cite[Section 8]{Janson2023} and e.g.\ 
   \cite[Chapter III]{Athreya-Ney1972})
 \begin{align}\label{eq:NB}
 \E \cY_t = me^{(m-1)t}
 \end{align}
and as $\ttoo$ and thus $x=e^{-t}\to 0$,
\begin{align}\label{eq:GG}
    x^{m-1}\hcY_{x}=e^{-(m-1)t}\cY_t 
\asto \xi\in \mathrm{Gamma}\bbclr{\frac{m}{m-1},m-1}.
\end{align}

The statements of \refL{LU} and \refT{th:coupling} hold with the same $n_1$
and $\kappa_n$, but $\chi$ is now as in \eqref{de:chi}, and $\delta_n=3\log^m n/n^{1/3}$. The analogue of Lemma \ref{LU} can be proved using entirely the same argument, but several straightforward modifications are needed to obtain the analogue of \refT{th:coupling}. For instance, in Step \ref{TCO1}, we use \eqref{eq:e2} and $n_1=\floor{n/\log n}$ to show that 
\begin{align}
    \sum^{n-1}_{k=n_1} \IP(Z_k\ge 2) = O\bclr{n^{-1}(\log n)^{1+2(m-1)\chi}}=o(1).
\end{align}
In Step \ref{TCO2}, it follows from \eqref{eq:NB} that
\begin{align}
    \E \wh \cY_x = \E \cY_{-\log x} = \frac{m}{x^{m-1}}, \quad 0<x\le1,
\end{align}
and so for $x=(n_1/n)^\chi\sim \log^\chi n$, w.h.p.\ there are at most $\log^{m-1} n$ generations from the root $1$ to any point in $[(n_1/n)^\chi,1]$. Adjusting the remaining steps accordingly then yield the desired conclusion.  

Redefine 
\begin{align}\label{Xi+}
    \Xi^{(n)} := \frac{W^{(n)}_{n_1}}{n^{(m-1)\chi}}.
\end{align}
It follows from
\eqref{Xi+}, \eqref{de:Wk}. \eqref{de:phi},
\eqref{lb2b},
\eqref{eq:Yulecoupling},
and \eqref{eq:GG}, that
the statement of \refL{LXi} holds, with $\tilde \beta$ as in \eqref{lb2b},
and $\xi\in \mathrm{Gamma}(m/(m-1),m-1)$, independent of $\tilde \beta$.

\subsection{The flat middle part}\label{se:flat1}
We first note that $\nu$ defined in \eqref{de:nu} also satisfies,
by \eqref{de:chi} and a simple calculation, 
\begin{align}\label{de:nu2}
    \nu = \frac{(m-1)\chi}{1+(m-1)\chi}
\end{align}
and thus
\begin{align}\label{de:nu3}
(1-\nu)(m-1)\chi = \nu.
\end{align}
We now choose $n^{\nu}\ll n_2\le n_1:=\floor{n/\log n}$.
Then, as in Section \ref{se:flat}, 
 \eqref{eq:e3} and \eqref{de:nu2} yield
    \begin{equation}\label{eq:IImaxAk+}
\E  \max_{n_2\leq k\leq n-1} \bigg|\frac{A_k}{n^{(m-1)\chi}}\bigg| 
= \frac{\E A_{n_2}}{n^{(m-1)\chi}} \le C\frac{n^{(m-1)\chi}}{n_2^{1+(m-1)\chi}}
= C\Bigpar{\frac{n^\nu}{n_2}}^{1+(m-1)\chi} = o(1).
    \end{equation}
Hence \refL{le:maxAk} holds, with denominators $n^{(m-1)\chi}$.
Similarly,
the proofs in Section \ref{se:flat}, 
now using \eqref{eq:vm}, \eqref{eq:mprodbeta}, \eqref{eq:phibd}, and 
\eqref{Xi+}, 
yield the conclusion that as $n\to\infty$, 
\begin{align}
    \max_{0\le k\le n_1} \bigg|\frac{M_k}{n^{(m-1)\chi}}-\Xi^{(n)}\bigg| &\pto 0; \label{eq:cm} \\
     \max_{n_2\le k\le n_1}\bigg|\frac{W_k}{n^{(m-1)\chi}}-\Xi^{(n)}\bigg| &\pto 0. \label{eq:cw} 
\end{align}

\subsection{The final part}
Let $V_k$ and $T_k$ be as in \eqref{ha3}.
As before, we extend $W_k, M_k, A_k, V_k, T_k$ to real arguments by linear
interpolation. 
Now let, for $t\ge0$,
\begin{align}
    \wh A\nn_t &:= n^{-(m-1)\chi} A\nn_{tn^\nu}, \label{hA+} 
\\
    \wh V\nn_t&:= n^\nu V\nn_{tn^\nu},\label{hV+}
\\
    \wh T\nn_t&:= T\nn_{tn^\nu}.\label{hT+}
\end{align}
Note that \refL{LV2} holds in this more general setting
(with the exponent $1/3$ replaced by $\nu$ in the proof).
Moreover, fix $a>0$ and define also
\begin{align}
  \xM_n:=n^{-(m-1)\chi}\max_{k\ge1} M_k
\qquad\text{and}\qquad
\Phix_n:= n^{\nu (m-1)\chi}\Phi\qw_{\floor{a n^{\nu}}}. 
\end{align}
In view of \eqref{eq:e3}, we have, for all real $\floor{an^\nu}\le k\le \ell$, 
\begin{align}\label{eq:A1}
    |A_\ell-A_k| \le C n^{(m-1)\chi(2-\nu)}\xM_n^2 \Phix_n (V_\ell-V_{k}).
\end{align}
Since $(2-\nu)(m-1)\chi - (m-1)\chi=\nu$ by \eqref{de:nu3}, 
it follows from
\eqref{eq:A1}  that if $a\le s\le t$, then
\begin{align}
    |\wh A\nn_t-\wh A\nn_s| \le C\xM_n^2 \Phix_n (\wh V\nn_t-\wh V\nn_s).
\end{align}
By \eqref{hA+}, \eqref{eq:e3}, and \eqref{de:nu2},
$\E \wh A\nn_a\le C_a$, which implies that the sequence 
$(\wh A\nn_a)\nxoo$
is tight. 
 From \eqref{eq:e1} and \eqref{eq:phibd}, $\E \xM^2_n\le C$ and $\E \Phix_n=O(1)$, implying that $( \xM_n)^\infty_{n=1}$ and $(\Phix_n)^\infty_{n=1}$ are tight also. Following the proof of \refL{LA2}, with the ingredients above, we then conclude that the stochastic processes $\wh A\nn_t$, $n\ge 1$, are tight in $C[a,b]$ for $0<a<b<\infty$.

 Therefore, arguing as in the beginning of the proof of \refT{th:cvg}, we
 may assume (by considering a subsequence and a special coupling)
that \eqref{aw1} holds a.s.\ in $C(0,\infty)$ together with 
\eqref{jup2}, \eqref{lxi},
\eqref{eq:cm},  \eqref{eq:cw},
and, instead of \eqref{jup1},
\begin{align}\label{jup1+}
\sup_{k\ge\log n}\bigabs{k^{-(m-1)\chi}\Phi\nn_k- \tgb\nn}\pto0.
\end{align}
Then, 
analogously to \eqref{aw2} and \eqref{aw4},
 \begin{align}\label{gaw1}
     n^{-(m-1)\chi} W_{tn^\nu} \to \cB(t):=\tilde\beta \xi -\cA(t) 
 \end{align}
and, again using \eqref{de:nu3},
\begin{align}\label{gaw2}
    n^{-\nu} Y_{\floor{tn^\nu}} \to t^{-(m-1)\chi}\tilde \beta^{-1} \cB(t)
\end{align}
a.s.\ uniformly on every compact interval in $(0,\infty)$. 

Let $k:=\floor{sn^\nu}$ and $\ell:=\floor{tn^\nu}$ for some $0<s<t$. 
Similarly to \eqref{aw5}, we have by \eqref{aw1} and \eqref{eq:AA}
\begin{align}\label{gaii1}
    \wh A\nn_s-\wh A\nn_t = mn^{-(m-1)\chi}\sum^\ell_{i=k+1}\Phi_{i-1} [(1-B_i)^{Y_i}-1+Y_iB_i]+ o(1).
\end{align}
Following \eqref{aw6}--\eqref{aw10} with minor adjustments (in particular,
choosing $\lambda_n=n^\nu$ in \refL{LH} and using again \eqref{de:nu3})
leads to, a.s.\ for every real $0<s\le t$, 
\begin{multline}\label{gaii2}
    \cA(s)-\cA(t) 
= m\tilde \beta \int^t_s u^{(m-1)\chi}\Bigl(\bclr{1+\tfrac{1}{\theta} u^{-(1+(m-1)\chi)}\tilde\beta^{-1}\cB(u)}^{-(m+\rho)} \\
    \qquad  -1+\chi u^{-(1+(m-1)\chi)}\tilde\beta^{-1}\cB(u)\Bigr)\dd u.
\end{multline}
Let, for convenience,  recalling \eqref{de:nu2},
\begin{align}\label{de:al}
     \al := 1+(m-1)\chi
=\frac{1}{1-\nu}
.\end{align}
Then, by \eqref{gaii2}, a.s.\ $\cA(t)$ is differentiable on $(0,\infty)$ and
\begin{equation}\label{gaii3}
    \cA'(t) = -m\tilde \beta t^{(m-1)\chi} 
\Bigl(
\bclr{1+\tfrac{1}{\theta} t^{-\ga}\tilde\beta^{-1}\cB(t)}^{-(m+\rho)}
 -1+\chi t^{-\ga}\tilde\beta^{-1}\cB(t)\Bigr),
\end{equation}
and $\cB'(t)=-\cA'(t)$ by \eqref{gaw1}. Define again
\begin{align}\label{eq:bb}
    \wt \cB(t) = \tilde\beta^{-1}\cB(t),
\end{align}
so that 
\begin{align}\label{eq:de}
     \wt \cB'(t)  
     = mt^{(m-1)\chi} \bbclr{\bclr{1+\tfrac{1}{\theta}t^{-\ga}\wt \cB(t)}^{-(m+\rho)}-1+\chi t^{-\ga}\wt \cB(t)}.
\end{align}
Moreover, $\E \cA(t)\le Ct^{-\ga}$ for $t\ge1$, say,
by \eqref{hA+}, \eqref{eq:e3}, and Fatou's lemma, 
and dominated convergence further implies that $\E\lim_{t\to\infty}\cA(t)=0$. 
Hence $\cA(t)\to 0$ a.s.\ as $t\to\infty$, 
and thus we have we have from \eqref{gaw1} and \eqref{eq:bb} that 
\begin{align}\label{eq:aw16}
  \tcB(t)\upto\xi
\qquad\text{as \ttoo}.
\end{align}
As shown in detail in \refApp{Sdiff}, see \eqref{eq:de3} and \eqref{eq:gc0}, 
the unique solution to \eqref{eq:de} satisfying \eqref{eq:aw16} is given by
\begin{align}\label{eq:tB}
    \wt \cB(t) = \theta t^\al \bbclr{\bclr{1+\tfrac{m+\rho+1}{\theta} \xi t^{-\al}}^{\frac{1}{m+\rho+1}}-1}.
\end{align}
Hence, by \eqref{eq:bb} and \eqref{eq:tB}, 
\begin{align}\label{eq:b}
    \cB(t)=\tgb \wt \cB(t)=\theta \tgb t^\al \bbclr{\bclr{1+\tfrac{m+\rho+1}{\theta}\xi t^{-\al}}^{\frac{1}{m+\rho+1}}-1}. 
\end{align}
Now, proceeding as in the remaining steps of the proof of \refT{th:cvg}
yields the conclusion that as $\ntoo$, 
\begin{equation}\label{cvg3+}
       n^{-(m-1)\chi} W_{tn^{\nu}}
\dto \theta \tgb t^\al \bbclr{\bclr{1+\tfrac{m+\rho+1}{\theta}\xi
    t^{-\al}}^{\frac{1}{m+\rho+1}}-1}
\qquad\text{in $C[0,\infty)$},
    \end{equation}
and
\begin{equation}\label{cvg4+}
    n^{-\nu}Y_{{tn^\nu}} \dto\theta t\bbclr{\bclr{1+\tfrac{m+\rho+1}{\theta}\xi t^{-\al}}^{\frac{1}{m+\rho+1}}-1}
\qquad\text{in $C(0,\infty)$}.
\end{equation}

\subsection{The number of descendants}
As in Section \ref{se:desc}, let $X=X^{(n)}$ be the number of red vertices,
and define $L_k$ and $P_k$ as in \eqref{eq:Xm} and \eqref{eq:Xpred}. 
As in \eqref{eq:EPk}, it follows from \eqref{eq:e2} that
\begin{align}
    \E P_k \le C \bbclr{\frac{n}{k}}^{(m-1)\chi},\qquad 1\le k\le n-1.
\end{align}
Furthermore, arguing as in \eqref{eq:EP0} with the cutoff $n^\nu$ yields,
recalling \eqref{de:nu3},
\begin{align}\label{EP0}
  \E P_0 \le C n^\nu.
\end{align}
The argument  for \eqref{eq:L02} now yields
\begin{align}\label{EL0}
\E L_0^2 \le C n^\nu,  
\end{align}
which implies that 
\begin{align}\label{eq:cL}
    n^{-\nu}L_0 \pto 0\qquad \text{as $\ntoo$}.
\end{align}
As before, we extend $P_k$ to real arguments by linear interpolation, but now let 
\begin{align}\label{nov2}
    \wh P\nn_t = n^{-\nu} P\nn_{tn^\nu}, \qquad t\ge 0.
\end{align}
The same proof as for \refL{le:Pt} then shows that for $0<a<b<\infty$, 
the sequences $\wh P\nn_t $, $n\ge 1$ are tight in $C[a,b]$. Proceeding as
in the proof of \refT{th:X}, where we use the Skorohod coupling theorem again,
we get
\begin{align}\label{eq:cP}
    \wh P\nn_0 &\pto   
\int^\infty_0 \bbclr{1-\bclr{1+\gth\qw u^{-\ga}\tcB(u)}^{-(m+\rho)} }\dd u.
\notag\\&\qquad
=    \int^\infty_0 \bbclr{1-\bclr{1+\tfrac{m+\rho+1}{\theta} \xi u^{-\al}}^{-\frac{m+\rho}{m+\rho+1} } }\dd u.
\end{align}
By the change of variable $v=\theta^{-1}(m+\rho+1)\xi u^{-\al}$,
\begin{align}\label{nov5}
&\int^\infty_0 \bbclr{1-\bclr{1+\tfrac{m+\rho+1}{\theta} \xi
  u^{-\al}}^{-\frac{m+\rho}{m+\rho+1} } }\dd u
\notag\\&\qquad
=\frac{1}{\al} \bbclr{\frac{m+\rho+1}{\theta}\xi}^{1/\al} \int^\infty_0 \bigg(1-{(1+v)^{-\tfrac{m+\rho}{m+\rho+1}}}\bigg) v^{-(1+1/\al)} \dd v.
\end{align}

We take $a=-1/\al$ and $b={(m+\rho)}/({m+\rho+1})$ in \refL{Lnora}, 
and note that $1/\ga=1-\nu$ by \eqref{de:al}, and thus by \eqref{de:nu},
\begin{align}
  b-a=
\frac{m+\rho}{m+\rho+1}+\frac{1}{\ga}
=\frac{m+\rho}{m+\rho+1}-\nu+1
=\frac{m+\rho}{m(m+\rho+1)}+1
.\end{align}
Hence, \eqref{nov5} and \refL{Lnora} yield
\begin{align}
    &\int^\infty_0 \bbclr{1-\bclr{1+\tfrac{m+\rho+1}{\theta} \xi u^{-\al}}^{-\frac{m+\rho}{m+\rho+1}} }\dd u \notag\\
    &\qquad = -\frac{1}{\al} \cdot \frac{\G(-\tfrac{1}{\al})\G\big(\tfrac{m+\rho}{m(m+\rho+1)}+1\big)}{\G\big(\frac{m+\rho}{m+\rho+1}\big)} \bbclr{\frac{m+\rho+1}{\theta}\xi}^{1/\al} \notag\\
    &\qquad =\frac{\G(1-\frac{1}{\al})\G\big(\frac{m+\rho}{m(m+\rho+1)}+1\big)}{\G\big(\tfrac{m+\rho}{m+\rho+1}\big)} \bbclr{\frac{m+\rho+1}{\theta}\xi }^{1/\al} . \label{eq:ir}
\end{align}
Finally, \eqref{eq:X},  \eqref{eq:cL}, \eqref{nov2},
\eqref{eq:cP}, and \eqref{eq:ir} together imply that as $\ntoo$,
\begin{align}\label{nov3}
  n^{-\nu} X\dto \frac{\G(1-\tfrac{1}{\al})\G\big(\frac{m+\rho}{m(m+\rho+1)}+1\big)}{\G\big(\tfrac{m+\rho}{m+\rho+1}\big)} \bbclr{\frac{m+\rho+1}{\theta}\xi }^{1/\al}. 
\end{align}
We here note that, by \eqref{de:al} and \eqref{de:nu}, 
\begin{align}\label{nova}
  1-\frac{1}{\ga}&=\nu=
\frac{(m-1)(m+\rho)}{m(m+\rho+1)}
.\end{align}
We write also $\xi=(m-1)\xi_1$, with $\xi_1\in\GAMMA(m/(m-1),1)$,
and recall that $\gth=2m+\rho$.
Hence, \eqref{nov3} can be written as \eqref{tmain}.
\qed

\section{Moment convergence}\label{Smom}
In this section, we prove \refT{Tmom} on moment convergence;
we use the standard method of proving uniform moment estimates and 
thus uniform integrability. 
This time we choose to treat general $m$ and $\rho$ from the beginning.

We consider first the reverse martingale $M_k$, recalling that $M_k\ge W_k\ge 0$. We denote
the maximal function by
\begin{align}\label{mx}
  \Mx:=\max_{n-1\ge k\ge 0} M_k,
\end{align}
and define 
the martingale differences, for $n-1\ge k\ge 1$, 
recalling \eqref{yngve}, \eqref{de:gmw}, \eqref{eq:sr},
and that $Y_k$ is $\cF_k$-measurable,
\begin{align}\label{gdm}
  \gD M_k&:= M_{k-1}-M_k 
=W_{k-1}-\E\bigpar{W_{k-1}\mid\cF_k}
\notag\\&\phantom:
=\GF_{k-1} \bigpar{Y_{k-1}-\E\bigpar{Y_{k-1}\mid\cF_k}}
\notag\\&\phantom:
=-\GF_{k-1} \bigpar{Z_k-\E\bigpar{Z_k\mid\cF_k}}
+m\GF_{k-1} \bigpar{J_k-\E\bigpar{J_k\mid\cF_k}}
.\end{align}
We define also the
conditional square function
\begin{align}\label{ssm}
 s(M):=
\lrpar{\sum_{i=1}^{n-1}\E\bigpar{\xpar{\gD M_i}^2\mid\cF_i}}\qq.
\end{align}
Let for convenience
\begin{align}\label{gk}
  \gk:=(m-1)\chi.
\end{align}
(Thus, in the case $m=2$, $\rho=0$, we have $\gk=\chi=\frac12$.)
We use also the standard notation, for any random variable $\cX$,
\begin{align}
  \norm{\cX}_p:=\bigpar{\E[|\cX|^p]}^{1/p}.
\end{align}
Note that for any  $p>0$, \eqref{de:betas}, \eqref{betamom}, and \eqref{gg}
yield, cf.\ \eqref{eq:Bmean}--\eqref{eq:Ebetasq},
\begin{align}\label{brage}
  \E[B_k^p]\le C_p k^{-p}.
\end{align}

\begin{lemma}\label{LpM}
  For every $p>0$,
\begin{align}\label{lpm}
 \E[\Mx^p]  \le C_p n^{p\gk}.  
\end{align}
\end{lemma}
\begin{proof}
We assume in the proof for simplicity that $p\ge2$ is an integer;
the general case follows   by Lyapunov's inequality. 

We use as in \cite{Janson2023} 
one of Burkholder's martingale inequalities 
\cite[Theorem 21.1]{Burkholder1973},
\cite[Corollary 10.9.1]{Gut}
on the reverse martingale $M_k-M_{n-1}=M_k-m\GF_{n-1}$,
which yields 
\begin{align}\label{burk}
  \E [\Mxx^p]&
\le C_p\E[\GF_{n-1}^p] + C_p\E \bigsqpar{\bigpar{\max_k|M_k-M_{n-1}|}^p}
\notag\\&
\le C_p \E[\GF_{n-1}^p]
+C_p \E [s(M)^p] + C_p \E \bigsqpar{\max_k|\gD M_k|^p}
\notag\\&
\le C_p\E[\GF_{n-1}^p]+C_p \E \bigsqpar{s(M)^p} 
+ C_p \sum_{k=1}^{n-1}\E \bigsqpar{|\gD M_k|^p}
.\end{align}
We estimate the three terms on the \rhs{} separately.

First, we have 
by the independence of $B_i$, \eqref{freja}, and \eqref{eq:Bmean},
similarly to \eqref{eq:trunphi2bd}--\eqref{eq:phibd}, 
\begin{align}\label{frej}
  \E [\GF_k^p]&
=\prod_{i=1}^k \E{(1+(m-1)B_i)^p}
=\prod_{i=1}^k \Bigpar{1+p(m-1)\frac{\chi}{i}+O\bigpar{i^{-2}}}
\notag\\&
=\exp\Bigpar{\sum_{i=1}^k \Bigpar{\frac{p\gk}{i}+O\bigpar{i^{-2}}}}
=\exp\Bigpar{p\gk \log k+O\bigpar{1}}
\notag\\&
\le C_pk^{p\gk}
.\end{align}

Next, by \eqref{gdm}, \eqref{eq:gvw}, and \eqref{de:gmw},
\begin{align}\label{sw1}
\E\bigsqpar{\xpar{\gD M_k}^2\mid\cF_k}
&=\Var\bigsqpar{W_{k-1}\mid\cF_k}
\le C\GF_{k-1}^2 B_kY_k 
\le C \GF_{k-1}B_kW_k
\notag\\&
\le C \GF_{k-1}B_kM_k
\le C \GF_{k-1}B_k\Mx.
\end{align}
Note that $\GF_k-\GF_{k-1}=(1+(m-1)B_k-1)\GF_{k-1}=(m-1)B_k\GF_{k-1}$.
Hence, \eqref{ssm} and \eqref{sw1} yield
\begin{align}\label{sw2}
  s(M)^2\le C\sum_{k=1}^{n-1}(\GF_k-\GF_{k-1})\Mx\le C\GF_{n-1}\Mx.
\end{align}
H\"older's inequality (or Cauchy--Schwarz's) and \eqref{frej} thus yield
\begin{align}\label{sw3}
 \E[s(M)^p]
\le C_p\E\bigsqpar{\GF_{n-1}^{p/2}\Mxx^{p/2}}
\le C_p\Bigpar{\E\sqpar{\GF_{n-1}^{p}}\E[\Mxx^{p}]}\qq
\le C_p n^{p\gk/2}\norm{\Mx}_p^{p/2}.
\end{align}

For the final term in \eqref{burk}, we use the decomposition of $\gD M_k$
in \eqref{gdm} and treat the two terms on the last line there separately.
We use as in \cite[(7.9)]{Janson2023} the well-known general estimate
for a binomial random variable $\zeta\in\Bin(N,q)$:
\begin{align}\label{binp}
  \E |\zeta-\E\zeta|^p \le C_p (Nq)^{p/2} + C_pNq
.\end{align}
Conditioned on $\cF_k$, we have $Z_k\in\Bin(Y_k,B_k)$ by \eqref{eq:distZk},
and thus \eqref{binp} yields
\begin{align}\label{sw4}
  \E \bigpar{\bigabs{Z_k-\E(Z_k\mid\cF_k)}^p\mid\cF_k}&
\le C_p (Y_kB_k)^{p/2} + C_p Y_kB_k.
\end{align}
Similarly, since $J_k=\tone[Z_k\geq 1]$ has a conditional Bernoulli
distribution, 
\begin{align}\label{sw7}
  \E \bigpar{\bigabs{J_k-\E(J_k\mid\cF_k)}^p\mid\cF_k}&
\le C_p   \E \bigpar{\abs{J_k}^p\mid\cF_k}
= C_p  \E \bigpar{J_k\mid\cF_k}
\notag\\&
\le C_p \E \bigpar{Z_k\mid\cF_k}
= C_p Y_kB_k.
\end{align}
Hence, \eqref{gdm}, \eqref{sw4}, and \eqref{sw7} yield,
\begin{align}\label{qw4}
\E \bigsqpar{|\gD M_k|^p\mid\cF_k}
&\le 
C_p\GF_{k-1}^p\bigsqpar{
  \E \bigpar{\bigabs{Z_k-\E(Z_k\mid\cF_k)}^p\mid\cF_k}
+
  \E \bigpar{\bigabs{J_k-\E(J_k\mid\cF_k)}^p\mid\cF_k}
}
\notag\\&
\le C_p \GF_{k-1}^pY_k^{p/2}B_k^{p/2} + C_p\GF_{k-1}^p Y_kB_k
\notag\\&
\le C_p \GF_{k-1}^{p/2}W_k^{p/2}B_k^{p/2} + C_p\GF_{k-1}^{p-1} W_kB_k
\notag\\&
\le  C_p \GF_{k-1}^{p/2}\Mx^{p/2}B_k^{p/2} + C_p\GF_{k-1}^{p-1} \Mx B_k.
\end{align}
Hence, using H\"older's inequality, the independence of $\GF_{k-1}$ and $B_k$,
\eqref{frej}, and \eqref{brage},
\begin{align}\label{sw5}
\E \bigsqpar{|\gD M_k|^p}&
\le C_p \E\bigsqpar{\GF_{k-1}^{p/2}B_k^{p/2}\Mxx^{p/2}}
+C_p \E\bigsqpar{\GF_{k-1}^{p-1}B_k\Mxx}
\notag\\&
\le C_p \Bigpar{\E\bigsqpar{\GF_{k-1}^{p}B_k^{p}}\E \bigsqpar{\Mxx^{p}}}\qq
+C_p \Bigpar{\E\bigsqpar{\GF_{k-1}^{2p-2}B_k^{2}}\E \bigsqpar{\Mxx^{2}}}\qq
\notag\\&
\le C_p k^{p\gk/2-p/2}\norm{\Mx}_p^{p/2}
+C_p k^{(p-1)\gk-1}\norm{\Mx}_2
\notag\\&
\le C_p k^{p\gk/2-1}\norm{\Mx}_p^{p/2}
+C_p k^{(p-1)\gk-1}\norm{\Mx}_p
.\end{align}
Consequently, 
\begin{align}\label{sw6}
\sum_{k=1}^{n-1} \E \bigsqpar{|\gD M_k|^p}&
\le C_p n^{p\gk/2}\norm{\Mx}_p^{p/2}
+C_p n^{(p-1)\gk}\norm{\Mx}_p
.\end{align}

Finally, \eqref{burk} yields, collecting the estimates  \eqref{frej},
\eqref{sw3}, and \eqref{sw6},
\begin{align}\label{sw9}
\E[\Mx^p]&
\le C_p n^{p\gk}+C_p n^{p\gk/2}\norm{\Mx}_p^{p/2}+C_p n^{(p-1)\gk}\norm{\Mx}_p
.\end{align}
It follows trivially from the definitions 
that for every $n$, $\Mx$ is deterministically bounded by some constant
(depending on $n$), and thus $\norm{\Mx}_p<\infty$.
Let $x:=\norm{\Mx}_p/n^{\gk}\in(0,\infty)$; 
then \eqref{sw9} can be written as
\begin{align}\label{sw99}
  x^p \le C_p + C_p x^{p/2} + C_p x.
\end{align} 
Since $p>1$, it follows that $x\le C_p$,
which is the same as \eqref{lpm}. 
Alternatively, we can proceed as in \cite{Janson2023} to consider only $p=2^j$, with $j$ being positive integers. The conclusion \eqref{lpm} then follows from an induction over $j$, \eqref{sw9} and the base case $(p=2)$ proved in \eqref{eq:e1}.
\end{proof}

We use the decomposition $X=1+L_0+P_0$ in \eqref{eq:X}, and estimate the
terms $L_0$ and $P_0$ separately.

\begin{lemma}\label{LpP}
  For every $p>0$,
\begin{align}\label{lpp}
 \E[P_0^p]  \le C_p n^{p\nu}.  
\end{align}
\end{lemma}
\begin{proof}
We may by Lyapunov's  inequality assume that $p\ge1$ is an integer.
  By \eqref{eq:Xpred} and \eqref{eq:distZk},
\begin{align}\label{dk1}
    P_k = \sum^{n-1}_{i=k+1} \E(J_i\mid \cF_i)
\le \sum^{n-1}_{i=k+1}Y_i B_i
= \sum^{n-1}_{i=k+1}\GF_i\qw W_i B_i
\le \Mx \sum^{n-1}_{i=k+1}\GF_{i-1}\qw B_i.
\end{align}
Hence, by H\"older's and Minkowski's inequalities,
\begin{align}\label{dk2}
\norm{P_k}_p
\le \norm{\Mx}_{2p} \lrnorm{\sum^{n-1}_{i=k+1}\GF_{i-1}\qw B_i}_{2p}
\le \norm{\Mx}_{2p} \sum^{n-1}_{i=k+1}\bignorm{\GF_{i-1}\qw B_i}_{2p}.
\end{align}
We have
$(1+x)^{-p}\le 1-px+C_px^2$
for all $x\ge0$, 
and thus
by \eqref{eq:Bmean}--\eqref{eq:Ebetasq} and \eqref{gk},
generalizing \eqref{dx0},
\begin{align}\label{dx1}
  \E\bigsqpar{(1+(m-1)B_i)^{-p}} &
\leq 1-p(m-1)\E [B_i]+C_p\E [B_i^2 ]
\notag\\&
= 1- p\gk i^{-1} + O(i^{-2}).
\end{align}
Hence,
by the same argument as for \eqref{frej}, for any integers $p\ge1$ and $k\ge1$,
\begin{align}\label{dx2}
  \E [\GF_k^{-p}]&
=\prod_{i=1}^k \E\bigsqpar{(1+(m-1)B_i)^{-p}}
=\prod_{i=1}^k \Bigpar{1-\frac{p\gk}{i}+O\bigpar{i^{-2}}}
\notag\\&
=\exp\Bigpar{-\sum_{i=1}^k \Bigpar{\frac{p\gk}{i}+O\bigpar{i^{-2}}}}
=\exp\Bigpar{-p\gk \log k+O\bigpar{1}}
\notag\\&
\le C_pk^{-p\gk}
.\end{align}
In other words,
$\norm{\GF_k\qw}_p\le C_p k^{-\gk}$. 
Furthermore, $\norm{B_k}_p\le C_p k\qw$ by \eqref{brage}.
Since $\GF_{i-1}$ and $B_i$ are independent, it follows that, for $i\ge2$,
\begin{align}\label{dx3}
\bignorm{\GF_{i-1}\qw B_i}_{p}= \norm{\GF_{i-1}\qw}_p\norm{B_i}_{p}
\le C_p i^{-\gk-1}.
\end{align}
We may here replace $p$ by $2p$, and it follows from \eqref{dk2} 
and \eqref{lpm}
that, 
for $k\ge1$,
\begin{align}\label{dx4}
\norm{P_k}_p
\le \norm{\Mx}_{2p}\sum^{n-1}_{i=k+1}\bignorm{\GF_{i-1}\qw B_i}_{2p}
\le  C_pn^{\gk}\sum^{n-1}_{i=k+1} i^{-\gk-1}
\le  C_pn^{\gk} k^{-\gk}
.\end{align}

Furthermore, as in \refS{se:desc}, we have $P_0-P_k \le k$ for any
$k\ge0$, and thus 
Minkowski's inequality and \eqref{dx4} yield,
choosing
$k:=\floor{n^\nu}$ and noting that 
\eqref{de:nu2} and \eqref{gk}
imply $\gk(1-\nu)=\nu$,
\begin{align}
\norm{P_0}_p
\le
\norm{P_k}_p+k
\le C_p n^{\gk-\gk\nu} +n^{\nu}
\le C_p n^{\nu}
,\end{align}
which completes the proof.
\end{proof}

\begin{lemma}\label{LpL}
  For every $p>0$,
\begin{align}\label{lpl}
 \E[|L_0|^p]  \le C_p n^{p\nu/2}.  
\end{align}
\end{lemma}

\begin{proof}
  Recall that $(L_k)_{k=0}^{n-1}$ is a reverse martingale.
By \eqref{eq:Xm}, its conditional square function is given by
\begin{align}
s(L)^2:=
\sum^{n-1}_{i=1}\E\bigsqpar{\bigpar{J_i-\E(J_i\mid \cF_i)}^2\mid\cF_i}
=\sum^{n-1}_{i=1}\Var\sqpar{J_i\mid\cF_i}
\le \sum^{n-1}_{i=1}\E\sqpar{J_i\mid\cF_i}
=P_0,
\end{align}
where the inequality follows because $J_i$ has a conditional Bernoulli
distribution. Furthermore, again using \eqref{eq:Xm},
the martingale differences $\gD L_k:=L_{k-1}-L_k$
are bounded by 
\begin{align}
  |\gD L_k| = \bigabs{J_k-\E(J_k\mid \cF_k)}
\le 1.
\end{align}
Hence, Burkholder's inequality yields, similarly to \eqref{burk},
using also \refL{LpP},
\begin{align}
  \E [L_0^p]
\le C_p \E [s(L)^p] +  C_p \E \bigsqpar{\max_k|\gD L_k|^p}
\le C_p \E [P_0^{p/2}] + C_p
\le C_p n^{p\nu/2},
\end{align}
which completes the proof.
\end{proof}

\begin{proof}[Proof of \refT{Tmom}]
  It follows from \eqref{eq:X} and \refLs{LpP} and \ref{LpL} that,
for any $p>0$,
  \begin{align}
    \E [X^p] \le C_p + C_p \E[L_0^p]+C_p\E[P_0^p]
\le C_p n^{p\nu}.
  \end{align}
In other words, 
$\E[(X\nn/n^\nu)^p] \le C_p$ for every $p>0$.
By a standard argument, see{} e.g.\ \cite[Theorems 5.4.2 and 5.5.9]{Gut},
this implies uniform integrability of the sequence
$|X\nn/n^\nu|^p$  for every $p>0$ and thus 
the convergence in distribution in \eqref{tmain} implies
convergence of all moments.

Since 
$\xi_1\in\GAMMA\bigpar{\frac{m}{m-1},1}$,
\begin{align}
\E \bigsqpar{\xi_1^{p(1-\nu)}}
=\frac{\gG(p(1-\nu)+\frac{m}{m-1})}{\gG(\frac{m}{m-1})}      
,\end{align}
and thus the explicit formula \eqref{tmom} follows.
\end{proof}

\section{The model with self-loops} \label{Sloop}
In this section, we consider a variation of the preferential attachment
graph in \refD{de:pa}, where self-loops are possible. 
We use the version in \cite[Section 8.2]{vdh2017} 
(see also \cite{Bollobas2004,Bollobas2001})
and start with a single vertex 1 with $m$ self-loops.
For $n\ge2$, 
each outgoing edge of vertex $n$ is now attached to a vertex
$j\in[n]$, again with probability proportional to 
$\rho$ + the current degree of vertex~$j$, where
we define the current degree of vertex $n$ when we add the $(k+1)$th edge
from it 
to be $k$ + $1$ + the number of loops attached to $n$ so far.
(We thus count all outgoing edges up to the $(k+1)$th;
a loop contributes 2 to the degree.)
Hence, recalling that $d_j(n)$
is the degree of vertex $j$ in $G_n$,
when adding
vertex $n\ge 2$ to $G_{n-1}$, the $(k+1)$-th outgoing edge of vertex $n$
attaches to vertex $j\in  [n]$ with probability 
\begin{align}\label{eq:pa3}
  \begin{cases}
\frac{d_j(n-1)+\sum^k_{\ell=1} \mathbf{1}[n\overset{\ell}{\rightarrow}j] +\rho}
{2(n-1)m+2k+1+n\rho},
&j<n,
\\
\frac{k+1+\sum^k_{\ell=1} \mathbf{1}[n\overset{\ell}{\rightarrow}j] +\rho}
{2(n-1)m+2k+1+n\rho},
&j=n.
  \end{cases}
\end{align}
\begin{remark}
The details of the model can be  modified  without affecting  the following
asymptotic result, with only straightforward changes to its proof.
For example, we may again start with $m$ edges between vertices 1 and 2, and
thus no loops there, 
or we may include all $m$ outgoing edges in the weight
of vertex $n$ when we add edges from it.
We leave the details to the reader.
\end{remark}

\begin{theorem}\label{Tloop}
    Let $X\nn$ be the number of descendants of vertex $n$ in the model
    above. Then, the statements of \refTs{Tmain} and \ref{Tmom}
    hold.
\end{theorem}

The proof of \refT{Tloop} is largely similar to those of \refTs{Tmain} and \ref{Tmom} so we
only sketch the main differences here. 

First, let $N_i$ be the number of self-loops at vertex $i$.
When we add the $m$ edges from a new vertex $i$, the weight of vertex $i$
and the total weight of the first $i-1$ vertices evolve like a P\'olya urn 
$\cU'_i$
with initially $1+\rho$ red and $(2i-2)m+(i-1)\rho$ black balls, where we
add 2 new balls at each draw: 2 red balls when a red ball is drawn, 
and one ball of each colour when a black ball is drawn;
$N_i$ is the number of times a red ball is drawn.
Note that this urn does not depend on what has happened when the edges from
earlier vertices were added, and in particular not on $N_1,\dots,N_{i-1}$.
Consequently, the random numbers $(N_i)^\infty_{i=1}$ are independent.
Furthermore, if we condition on the entire sequence $(N_i)^\infty_{i=1}$,
then the non-loop edges are added from each new vertex $n\ge2$ to $[n-1]$
by the same random procedure as in \refD{de:pa}, 
except that now we add $m-N_n$ new edges from $n$, and that the degrees of the
vertices include also any existing loops.
This means that after we have added vertex $j\ge2$, 
the weight of vertex $j$ and the total weight of the first $j-1$ vertices
evolve like a standard P\'olya urn $\cU''_j$ with initially $m+N_j+\rho$ red and
$(2j-1)m-N_j+(j-1)\rho$ black balls, after each draw adding one ball of the
same colour as the drawn ball.
As a consequence, the proportion of red balls converges a.s.\ to a random
number $B_j$ with the (conditional) beta distribution
\begin{align}\label{nB}
    B_j\mid(N_i)_{i=1}^\infty
\in\mathrm{Beta}(m+N_j+\rho, (2j-1)m-N_j+(j-1)\rho),
\qquad j\ge2.
\end{align}
Moreover,
conditioned on $(N_i)_{i=1}^\infty$,
we can again construct the preferential attachment graph by the
P\'olya urn representation in \refD{de:PUR}--\refR{RBerger},
using (conditionally) independent $B_j$ with the distributions
\eqref{nB}. (As before, we also let $B_1:=1$.)
In particular, note that since $(N_i)_{i=2}^\infty$ are independent,
the random variables $(B_i)_{i=2}^\infty$ are independent.
and so are the pairs of random variables $(N_i,B_i)$, $i\ge2$.

The distribution of each $B_j$ is thus a mixed beta distribution,
but we do not need exact expressions.
We will show that all estimates in \refS{Sprel} still hold (possibly with
different constants $C$).
Note first that in the urn $\cU'_i$ used to determine $N_i$, we make $m$
draws and thus the number of red balls is at most $m+\rho=O(1)$; hence the
probability of drawing a red ball is $O(1/i)$ for each draw, and thus
\begin{align}\label{ENi}
  \P(N_i>0) \le \E N_i = O(1/i).
\end{align}
Recall $\theta$ and $\chi$  in \eqref{de:theta} and \eqref{de:chi}.
Using $0\le N_i\le m$, \eqref{nB}, \eqref{ENi}, and \eqref{betamom}, it is
easy to show that 
\begin{align}
\E B_i&= \frac{m+\rho+\E N_i}{\theta i} = \frac{\chi}{i} + O(i^{-2})\label{bn},
  \\
\E B_i^r &\le \prod^{r-1}_{j=0}\frac{2m+\rho+j}{\theta i +j} 
\le C_ri^{-r}
,\quad r\ge 2.\label{bn1} 
\end{align}
Similarly, we have, by first conditioning on $N_i$,
\begin{align}\label{enb1}
    \E [N_iB_i] = \frac{\E[N_i(m+\rho+N_i)] }{\theta i} 
\le \frac{\E[N_i(2m+\rho)] }{\theta i} = O(i^{-2}),
\end{align}
and the bound
\begin{align}\label{enb2}
    \E \bigsqpar{N^r_i B_i^r}  =O\big(i^{-(r+1)}\big),
\qquad \text{for each }r\ge 2 .
\end{align}


Define $\Phi_i$ and $S_{n,i}$ by \eqref{de:phi} and \eqref{de:S} as before.
Then
\eqref{bn} and a little calculation using \eqref{gg} shows that, 
for $2\le j\le k<\infty$,
    \begin{align*}
\prod^k_{i=j}\E(1+(m-1)B_i)
&=
\prod^k_{i=j}\frac{i+(m-1)\chi}{i}
\cdot\prod^k_{i=j}\frac{\E(1+(m-1)B_i)}{1+(m-1)\chi/i}
\\
&= \frac{\G\bclr{k+1+(m-1)\chi}\G(j)}{\G\bclr{j+(m-1)\chi}\G(k+1)}
\prod^k_{i=j}\Bigpar{1+O\bigpar{i^{-2}}}
        \notag\\
        &= \bbclr{\frac{k}{j}}^{(m-1)\chi} \bclr{1+O(j^{-1})}\numberthis \label{eq:mprodbetaN}
,    \end{align*}
as in \eqref{eq:mprodbeta},
and, recalling that $B_i$ are independent and 
taking $j=1$ in \eqref{eq:mprodbetaN}, 
    \begin{align*}
\E \Phi_k
&=\prod^k_{i=1}\E(1+(m-1)B_i)
=
\frac{\G\bclr{k+1+(m-1)\chi}}{\G\bclr{1+(m-1)\chi}\G(k+1)}
\prod^k_{i=1}\frac{\E(1+(m-1)B_i)}{1+(m-1)\chi/i}
\\
&=Q k^{(m-1)\chi} \bclr{1+O(k^{-1})}\numberthis \label{pb0}
,    \end{align*}
where
\begin{align}\label{asa}
  Q:=\frac{\prod^\infty_{i=1}\frac{\E(1+(m-1)B_i)}{1+(m-1)\chi/i}}{\G\bclr{1+(m-1)\chi}};
\end{align}
note that the infinite product in \eqref{asa} converges as a consequence of
\eqref{bn}.

Using \eqref{bn} and \eqref{bn1},
the upper bounds \eqref{eq:trunphi2bd} and
\eqref{eq:phibd} follow by the same proof as before.
The statements in Lemmas \ref{LB2} and   \ref{le:Sest}
hold exactly, except for \eqref{lb2b}, which in view of \eqref{pb0}, is now replaced with 
\begin{align}\label{ntb}
    \tilde\beta :=Q\beta.
\end{align}

Let $Y_k, Z_k, J_k, W_k$ be as in \refS{se:basicanalysis}. To streamline the
arguments, from here onwards we concentrate on the $m=2$, $\rho=0$ case, and
leave the general case (with modifications as in \refS{Sgen})
to the reader.
Once we have sampled the self-loops at
every vertex, the stochastic recursions for obtaining $D_n$ are similar to
the ones in \refS{sse:sr}: we sample $(B_i)^{n-1}_{i=1}$ according to
\eqref{nB}, and for each red vertex $k$, we add $2-N_k$ outgoing edges and
proceed as before. The boundary conditions are the same, except now we have
$Y_{n-1} = 2 - N_n$. For $2\le k\le n-1$, the recursion takes the form 
\begin{align}\label{nY}
    Y_{k-1} = Y_k - Z_k + (2 - N_k) J_k , 
\end{align}
and because $0\le N_k \le 2$, we also have
\begin{align}\label{sr1}
    Y_{k-1} \le  Y_k - Z_k + 2 J_k.
\end{align}
Let $\cF_k$ be the $\sigma$-algebra generated by $(N_i)^n_{i=2}$,
$(B_i)^{n-1}_{i=2}$ and the coin tosses at vertices $n-1,\dots,k+1$ in the
stochastic recursion. Note that \eqref{eq:distZk} holds, and in view of
\eqref{sr1}, 
also \eqref{eq:markov}--\eqref{eq:Akdiff}
and \eqref{nov1}--\eqref{eq:varM} hold, with the number $2$ in
\eqref{eq:mg} and \eqref{eq:varM}  replaced with $2-N_n$,
and with the last equality in \eqref{eq:ubmeanW} replaced with $\le$. 
Instead of \eqref{eq:incA}, from
\eqref{nY} we have
\begin{align}\label{iA}
    A_{k-1} -A_k &= W_k - \E (W_{k-1}\mid \cF_k)\notag\\
    &= 2 \Phi_{k-1} ((1- B_k)^{Y_k}-1+B_k Y_k) + \Phi_{k-1} N_k \E(J_k\mid \cF_k).
\end{align}
Now, let $\bB$ be the $\sigma$-field generated by $(B_i)^{n-1}_{i=2}$ and
$(N_i)^{n}_{i=2}$. 
As the upper bounds in Section
\ref{Sprel} still hold, 
and $(B_i)^{n-1}_{i=2}$ are
independent, \refLs{le:doob} and \ref{le:Zk} hold. The probability that
vertex $k\ge 2$ is red and has at least one self-loop is 
\begin{align}
    \IP(Z_k\ge 1, N_k\ge 1)
= \E\bigsqpar{\mathbf{1}\set{N_k\ge 1} \IP_{\bB}(Z_k\ge 1)};
\end{align}
and so by Markov's inequality and \eqref{MarZ}, 
\begin{align}
     \IP(Z_k\ge 1, N_k\ge 1)\le 2 \E \bbclr{N_kB_k \prod^{n-1}_{i=k+1}(1+B_i) }.
\end{align}
By the independence of 
the pairs
$(B_i,N_i)$, 
\eqref{enb1}, and \eqref{eq:mprodbeta},
this yields
\begin{align}\label{eq:zn}
    \IP(Z_k\ge 1, N_k\ge 1)&\le 2 \E (N_kB_k) \prod^{n-1}_{i=k+1}\E(1+B_i) 
\le C\frac{n^{1/2}}{k^{5/2}}.
\end{align}

In view of \eqref{iA} and Markov's inequality,
\eqref{eq:incAsim} in
\refL{le:Ak} is replaced with 
\begin{align}\label{iA1}
    A_{k-1} -A_k \le (W_kB_k)^{2} \Phi_k^{-1} + \Phi_{k} N_k  B_k Y_k
= (W_kB_k)^{2} \Phi_k^{-1} +  N_k  B_k W_k
.\end{align}
Using \eqref{eq:ubmeanW}, we have
\begin{align}\label{iA3}
   \E_\bB \sqpar{N_k B_k W_k} 
=   N_k B_k \E_\bB \sqpar{W_k} \le 2 N_k B_k \Phi_{n-1}.
\end{align}
Since the pairs $(B_i,N_i)$ are independent, 
it follows from \eqref{iA3} and \eqref{de:phi}  that
\begin{align}
   \E \sqpar{N_k B_k W_k} 
\le 2 \E \bigsqpar{N_k B_k(1+B_k)} \prod^{n-1}_{\substack{i=1\\i\neq k}} \E (1+B_i)
\le 4 \E [N_k B_k] \E\Phi_{n-1}
,\end{align}
and applying
 \eqref{enb1} and \eqref{pb0}, we get 
\begin{align}\label{iA2}
   \E \bigsqpar{N_k B_k W_k} 
\le C\frac{n^{1/2}}{k^{2}} \le C\frac{n}{k^{5/2}}.
\end{align}
With \eqref{iA1} and \eqref{iA2}, we may proceed as in the proof of \refL{le:Ak} to show that \eqref{eq:EAk} holds.

\refL{LU} follows from \refL{le:Sest}. 
Thus, the early part of the growth of $D_n$ can be coupled to the same time-changed Yule process $\wh \cY$ with some extra  modifications. Recall that $\Psi(x)$ is the mapping of vertex $x$ in $\wh \cY$ to a vertex $k$ in $D_n$ (or vertex $(k/n)^\chi$ in $\xD_n$). In Step (1) of the coupling, we sample $(N_i)^n_{i=1}$ and then $(B_i)^{n-1}_{i=1}$ as in \eqref{bn}.  
If $\Psi$ maps $x$ to some $k$ that has at least one self-loop, we extend $\Psi$ in Section \ref{se:Yule} by mapping all children of $x$ to $k$ (so all other descendants of $x$ are also mapped to $k$). To prove that \refT{th:coupling} also holds in this case, we need to show that the extended mapping above is w.h.p.\ injective at every vertex in $\xD_n\cap [(n_1/n)^\chi,1]$. By \eqref{eq:zn} and \eqref{ENi}, the probability that a vertex in $\xD_n\cap [(n_1/n)^\chi,1]$ has at least one self-loop is at most  
\begin{multline}\label{loopp}
    \IP(N_n\ge 1) + \sum^{n-1}_{k=n_1}  \IP(Z_k\ge 1, N_k\ge 1)  
    \le   \frac{C}{n} +  \sum^{n-1}_{k=n_1} \frac{Cn^{1/2}}{k^{5/2}} = O(\log^{3/2} n/n) = o(1).
\end{multline}
 The same argument as in Step \ref{TCO1} in the proof of \refT{th:coupling} and \eqref{loopp} then give the desired claim. The remaining steps of the proof can be applied without any  changes.

\refL{le:maxAk}, \refL{le:maxMk} and \refT{TF1} hold with the same proofs as
before,
since we have shown that \eqref{eq:EAk} and the various other estimates that we
use there still hold.

Let $\wh A\nn_t$ be as in \eqref{hA}. When proving tightness of $\wh A\nn_t$
in $C[a,b]$ for $0<a<b<\infty$, we have to use
\eqref{iA1} instead of \eqref{eq:incAsim}. 
Let $V_k$, $T_k$, $\wh V\nn_t$, $\hT\nn_t$,
$\xM_n$, and $\Phix_n$ be as in \eqref{ha3},
\eqref{hV} and \eqref{hdj0}. Using \eqref{iA1}
we obtain instead of \eqref{hdj01}, using the crude bound $N_k\le m$,
  \begin{align}\label{tyr}
   | A_k- A_{k-1}| &
\le M_k^2\Phi_k\qw (V_k-V_{k-1})+ m M_k(T_k-T_{k-1})
\notag\\&
\le n^{5/6}\xM_n^2 \Phix_n (V_k-V_{k-1}) + m n^{1/2}\xM_n(T_k-T_{k-1})
  \end{align}
and thus, arguing as for \eqref{hdj6}, for real numbers $s,t$ such that 
$a\le s\le t$,
\begin{align}\label{balder}
    |\wh A\nn_t - \wh A\nn_s |
&\le \Phix_n \xM_n^2  (\hV\nn_t-\hV\nn_s) + m\xM_n (\hT\nn_t- \hT\nn_s).
\end{align}
We have already shown in Section \ref{Stig} 
that the processes $\wh V\nn_t-\wh V\nn_a$ and 
$\hT\nn_t- \hT\nn_a$ are tight in $C[a,b]$ (\refL{LV2})
and that the sequences $(\xM_n)^\infty_{n=1}$ and $(\Phix_n)^\infty_{n=1}$
are tight. 
Hence, by simple applications of \refL{LC}, the processes
$\Phix_n \xM_n^2  (\hV\nn_t-\hV\nn_a)$ and $m\xM_n (\hT\nn_t- \hT\nn_a)$,
$n\ge1$,
are  tight in $C[a,b]$.
If 
$(X_n(t))^\infty_{n=1}$ and $(Y_n(t))^\infty_{n=1}$ are any two sequences of
random continuous functions on $[a,b]$ that both  are tight in
$C[a,b]$, then so is the sequence $((X_n(t)+Y_n(t)))^\infty_{n=1}$. 
Hence, the sequence
$\Phix_n \xM_n^2  (\hV\nn_t-\hV\nn_a)+m\xM_n (\hT\nn_t- \hT\nn_a)$,
$n\ge1$, is tight in $C[a,b]$;
finally \eqref{balder}
and another application of \refL{LC} (now with $Z_n=1$)
show that 
$\wh A\nn_t$, $n\ge1$, are tight in $C[a,b]$, so \refL{LA2} still
holds.

Some minor adjustments are also required to yield the same result as in
\refT{th:cvg}. When applying the Skorohod coupling theorem, $N\nn_i$, $n\ge
1$, are potentially different for each $n$. 
Let $0<s<t$ and define $k:=\floor{sn^{1/3}}$ and $\ell:=\floor{tn^{1/3}}$. By \eqref{hA} and \eqref{iA}, 
\begin{multline}
  \hA\nn_s-\hA\nn_t 
=n^{-1/2}\sum_{i=k+1}^\ell2\GF_{i-1}\bigsqpar{(1-B_i)^{Y_i}-1+Y_iB_i} \\
 + n^{-1/2}\sum_{i=k+1}^\ell \GF_{i-1} N_i \E(J_i\mid \cF_i) +o(1)
.\end{multline}
However, by Markov's inequality and \eqref{iA2},
\begin{align}\label{at1}
\E \sum_{i=k+1}^\ell \GF_{i-1} N_i \E(J_i\mid \cF_i)&
\le \E \sum_{i=k+1}^\ell \GF_{i-1} N_iB_iY_i
\le \E \sum_{i=k+1}^\ell  N_iB_iW_i
\notag\\&
\le  \sum_{i=k+1}^\ell C\frac{n^{1/2}}{i^2} 
\le C\frac{n^{1/2}}{k} 
\le C_sn^{1/6},
\end{align}
implying that 
\begin{align}
n^{-1/2}\sum_{i=k+1}^\ell \GF_{i-1} N_i \E(J_i\mid \cF_i)\pto 0.
\end{align}
The remainder of the proof of \refT{th:cvg} is then the same as before. 

\begin{proof}[Proof of \refT{Tloop}]
With the preparations above, the same argument as in \refS{se:desc} yields 
\refT{th:X} for this model too; with modifications as in \refS{Sgen} we obtain
\refT{Tmain}. 
Similarly, the arguments in \refS{Smom} still hold, and thus
\refT{Tmom} holds.
\end{proof}

\appendix
\section{The differential equations in \eqref{aw13} and \eqref{eq:de}}
\label{Sdiff}
We rewrite the equation in \eqref{eq:de} as
\begin{align}\label{eq:de1}
     f'(t) 
     = mt^{\ga-1} \bbclr{\bclr{1+\tfrac{1}{\theta}t^{-\ga}f(t)}^{-(m+\rho)}-1+\chi t^{-\ga}f(t)}.
\end{align}
where, as above, 
\begin{align}\label{ga2}
\ga:=1+(m-1)\chi.  
\end{align}
Note that in the special case $m=2$ and $\rho=0$, 
we have $\chi=1/2$, $\theta=4$, and $\ga=3/2$,
so the above yields the differential
equation in \eqref{aw13}.
We define
\begin{align}\label{eq:gsg1}
    g(t) := \theta^{-1} t^{-\ga} f(t) 
\end{align}
so that \eqref{eq:de1} simplifies to,
recalling $\chi\gth=m+\rho$, see \eqref{de:chi},
\begin{align}\label{eq:de2}
    g'(t) = -\frac{\ga}{t}g(t) + \frac{m}{\theta t} \bbclr{(1+g(t))^{-(m+\rho)}-1+(m+\rho)g(t)}.
\end{align}
Letting 
\begin{align}\label{eq:gsh1}
    h(x):=g(e^{(m+\rho)x})
\end{align}
then yields, using \eqref{ga2},
\begin{align}\label{eq:gsh2}
    h'(x)&=-(m+\rho)(1+(m-1)\chi) h(x) + \chi m\bclr{(1+h(x))^{-(m+\rho)}-1+(m+\rho)h(x)}\notag\\
    &= \chi m (1+h(x))^{-(m+\rho)} -\chi m+ (m+\rho)(\chi -1)h(x)\notag\\
    &=\chi m\bclr{(1+h(x))^{-(m+\rho)} -1 - h(x)}, 
\end{align}
where the last equality follows from $(m+\rho)(1-\chi)=\chi m$,
see again \eqref{de:chi}. 
The autonomous differential equation in \eqref{eq:gsh2} can be integrated to 
\begin{align}\label{eq:gsh3}
   \frac{1}{\chi m}\int \frac{1}{(1+h)^{-(m+\rho)}-1-h} \dd h = \int 1 \dd x
\end{align}
Furthermore, 
with $v:=1+h$,
\begin{align}
     \frac{1}{\chi m}\int \frac{1}{(1+h)^{-(m+\rho)}-1-h} \dd h
    & = \frac{1}{\chi m} \int \frac{(1+h)^{m+\rho}}{1-(1+h)^{m+\rho+1}}\dd h 
\notag\\&
=- \frac{1}{\chi m} \int \frac{v^{m+\rho}}{v^{m+\rho+1}-1}\dd v.
\end{align}
The change of variable $u=v^{m+\rho+1}-1$ then gives
\begin{align}
     \frac{1}{\chi m} \int \frac{v^{m+\rho}}{v^{m+\rho+1}-1}\dd v 
= \frac{1}{\chi m(m+\rho+1)} \int \frac{1}{u}\dd u 
= \frac{1}{\chi m(m+\rho+1)} \log u + C.
\end{align}
Thus, reverting back to the original variable $h$, \eqref{eq:gsh3} is equivalent to
\begin{align}
    -\frac{1}{\chi m(m+\rho+1)} \log\bclr{(1+h)^{m+\rho+1}-1} = x+C,
\end{align}
which yields the solution  
\begin{align}\label{eq:gsh4}
    h(x) = \bclr{1+ce^{-\chi m(m+\rho+1)x} }^{\frac{1}{m+\rho+1}} -1, \qquad \text{for some $c\in \bbR$.} 
\end{align}
From \eqref{eq:gsg1} and \eqref{eq:gsh1},
\begin{align}\label{eq:gsh5}
    f(t) = \theta t^{1+(m-1)\chi} h\bclr{\tfrac{1}{m+\rho}\log t}.
\end{align}
so plugging in \eqref{eq:gsh4} into \eqref{eq:gsh5}, and using
\begin{align}\label{eq:amr}
   \al = 1+(m-1)\chi = 1+\frac{(m-1)(m+\rho)}{2m+\rho} = \frac{\chi m(m+\rho+1)}{m+\rho},
\end{align}
we get
\begin{align}\label{eq:de3}
    f(t) = \theta t^\al \bbclr{\bclr{1+c t^{-\al}}^{\frac{1}{m+\rho+1}}-1}.
\end{align}
Using L'H\^opital's rule (or a Taylor expansion) and \eqref{eq:de3}, we obtain
\begin{align}
    f(\infty) 
:=\lim_{t\to\infty}f(t)
    = \frac{\theta c}{m+\rho+1}\lim_{t\to\infty}\bclr{1+c t^{-\al}}^{\frac{1}{m+\rho+1}-1} = \frac{\theta c}{m+\rho+1}  .
\end{align}
Hence, the unique solution $f$ to \eqref{eq:de1} 
with a given $f(\infty)$ is given by \eqref{eq:de3} with
\begin{align}\label{eq:gc0}
    c = \frac{m+\rho+1}{\theta}f(\infty).
\end{align}

\section{{A beta integral}}

Recall the standard beta integral \cite[5.12.3]{NIST}
  \begin{align}\label{erika}
    \intoo\frac{x^{a-1}}{(1+x)^b}\dd x
&=
\frac{\gG(a)\gG(b-a)}{\gG(b)}
  \end{align}
when $0<\Re a<\Re b$.
We use the following less well-known  extension;  
it is not new but we give a proof for completeness.
\begin{lemma}\label{Lnora}
  If\/ $-1<\Re a<0$ and\/ $\Re b>0$, then
  \begin{align}\label{nore}
    \intoo\Bigpar{\frac{1}{(1+x)^b}-1}x^{a-1}\dd x
=\frac{\gG(a)\gG(b-a)}{\gG(b)}.
  \end{align}
\end{lemma}

\begin{proof}
  We consider a more general integral.
  Assume first $\Re a>0$ and $\Re b>0$, and let $\Re c>\Re a$.
Then, by using \eqref{erika} twice,
  \begin{align}\label{nora}
    \intoo\Bigpar{\frac{1}{(1+x)^{b+c}}-\frac{1}{(1+x)^{c}}}x^{a-1}\dd x
&=
\frac{\gG(a)\gG(b+c-a)}{\gG(b+c)}
-\frac{\gG(a)\gG(c-a)}{\gG(c)}.
  \end{align}
For fixed $b$ and $c$ with $\Re b,\Re c>0$, the \lhs{} converges for
$-1<\Re a<\Re c$, and defines an analytic function of $a$ in this strip.
Hence, by analytic continuation, \eqref{nora} holds  throughout this range.
Similarly, if $\Re a>-1$ and $\Re b>0$, then the \lhs{} of \eqref{nora} is an
analytic function of $c$ in the domain $\Re c>\Re a$, and thus
\eqref{nora} holds whenever $-1<\Re a<\Re c$ and $\Re b>0$.

For $-1<\Re a<0$ we thus may take $c=0$ in \eqref{nora} which yields
\eqref{nore}. 
(Recall that $1/\gG(0)=0$.)
\end{proof}

\begin{remark}
  Note that \eqref{erika} and \eqref{nore} give the same formula, but for
different ranges of $a$. The integrals can be interpreted as the Mellin
transforms of $(1+x)^{-b}$ and $(1+x)^{-b}-1$, respectively,
and thus this is an instance of a general phenomenon when considering the Mellin
transforms of a function $f(x)$ and of the difference $f(x)-p(x)$
where, for example, $p(x)$ is a finite Taylor polynomial at $0$,
see \cite[p.~19]{FGD}.
\end{remark}


\end{document}